\documentclass[reqno,11pt]{amsart}
\usepackage{amsmath,amssymb,mathrsfs,amsthm,amsfonts}
\usepackage[inline]{enumitem}
\usepackage[usenames,dvipsnames]{xcolor}

\usepackage{hyperref}
\usepackage[percent]{overpic}
\usepackage{comment}
\usepackage{stmaryrd}
\usepackage{algorithm}
\usepackage{algorithmic}
\usepackage{pgfplots}
\usepackage{subcaption}
\usepackage{tikz,todonotes}
\usetikzlibrary{calc}
\usetikzlibrary{arrows.meta,backgrounds}
\usepackage{multirow,array,longtable,booktabs}
\usetikzlibrary{arrows}
\usepackage{bm}
\usepackage{graphicx}

\hypersetup{
	colorlinks=true, linkcolor=blue,
	citecolor=ForestGreen
}
\usepackage[paper=letterpaper,margin=1in]{geometry}
\DeclareMathAlphabet{\mathpzc}{OT1}{pzc}{m}{it}

\newtheorem{theorem}{Theorem}[section]

\newtheorem{lemma}{Lemma}[section]
\newtheorem{proposition}{Proposition}[section]
\newtheorem{corollary}{Corollary}[section]
\newtheorem{definition}{Definition}[section]

\newtheorem{remark}{Remark}[section]
\theoremstyle{definition}

\numberwithin{equation}{section}
\allowdisplaybreaks
\usepackage{acronym}

\newcommand{\be}{\begin{equation}}
\newcommand{\ee}{\end{equation}}
\renewcommand{\Pr}{\mathbb{P}}

\newcommand{\Ex}{\mathbb{E}}

\newcommand{\ind}{\mathbf{1}}

\makeatletter

\newcommand{\R}{\mathbb{R}}
\newcommand{\Z}{\mathbb{Z}}

\newcommand{\e}{\varepsilon}
\newcommand{\op}{\mathrm{op}}

\newcommand{\T}{\mathbb{T}}

\renewcommand{\bar}{\overline}

\usepackage{comment}

\newcounter{hyp}
\usepackage{etoolbox}

\title[Turnpike Estimates for N-Player Games under Displacement Monotonicity]{Turnpike Estimates and Limiting Behavior for
$N$-Player Differential Games under Displacement Monotonicity}

\author[A.\ Cohen]{Asaf Cohen}
\address{A.\ Cohen,
	Department of Mathematics, University of Michigan,
	\newline\hphantom{\quad \ \ A. Cohen}
	530 Church St, Ann Arbor, MI 48109, USA}
\email{asafc@umich.edu}
\thanks{A. Cohen acknowledges the support of the U.S. National Science Foundation, Grant DMS-2505998.}

\author[A.\ Zitridis]{Antonios Zitridis}
\address{A.\ Zitridis,
	Department of Mathematics, University of Michigan,
	\newline\hphantom{\quad \ \ A. Zitridis}
	530 Church St, Ann Arbor, MI 48109, USA}
\email{zitridis@umich.edu}
\thanks{The authors wish to thank Marco Cirant and Alpár Mészáros for fruitful suggestions and discussions.}

\begin{document}
    \maketitle
\vspace{-10mm}
\begin{abstract}
We study the long-time behavior of open-loop Nash equilibria of $N$-player stochastic differential games with nondegenerate idiosyncratic noise, under displacement monotonicity of the costs and joint strong convexity of the Lagrangians. We prove an exponential turnpike estimate: as the remaining horizon grows, the decoupling fields of the finite-horizon Pontryagin systems converge exponentially fast to the solution of a stationary Pontryagin system, which is unique in a natural stability class, and on bounded time intervals the equilibrium states and controls converge to those generated by this solution. The limiting control is a strong open-loop equilibrium of an associated ergodic $N$-player game, and the ergodic constant $\lambda^{i,N}$ of each player is the same at every such equilibrium. For the finite-horizon game, the law of the equilibrium state, at times away from the initial and terminal times, is exponentially close in the Wasserstein distance $W_2$ to the invariant measure of the limiting dynamics.
Finally, the equilibrium cost of player $i$ on $[t,T]$ equals $(T-t)\lambda^{i,N}$ plus a uniquely determined corrector, up to an error exponentially small in $T-t$.
\end{abstract}

\section{Introduction}

This paper analyzes the long-time behavior of open-loop $N$-player stochastic differential games, as the horizon $T\to +\infty$. We consider the limiting behavior of their Nash equilibria and, in particular, establish the turnpike property and convergence to an equilibrium of an average-cost stochastic differential game.

\subsection{Overview} Finite-horizon and average-cost (also called ergodic) stochastic differential games describe different aspects of strategic interaction. In the former, players optimize over a prescribed time interval $[t_0,T]$ and attach a cost to their terminal states. In the latter, they compare the costs accumulated per unit time over an indefinitely long interaction. A basic question is whether the equilibria of a long finite-horizon ($T$ large) game select an average-cost equilibrium, and whether this selection can be quantified. One expects that the influence of the terminal condition decays as the remaining horizon grows, and that the state distribution approaches a stationary regime away from both ends ($t=t_0$ and $t=T$) of the time interval. This behavior is called the turnpike property
(see e.g. \cite{TrelatZuazua2015} for optimal control problems and the references therein for a historical overview), and an inequality that quantifies it is called a turnpike estimate.
\vspace{1mm}

We study this question directly for a finite population of $N\in \mathbb{N}$ players whose states take values in $\R^d$. Each player $i\in \{1,\ldots,N\}$ controls his/her own drift $\alpha^i$, and his/her state $X^{i,N}_t$ obeys the dynamics $dX_t^{i,N}=\alpha_t^i\,dt+\sqrt2\sigma^i\,dW_t^i,$ for a $d$-dimensional Brownian motion $W^i$ and a $d\times d$ invertible matrix $\sigma^i$. His/her goal is to minimize the following cost, which depends on the configuration of the entire population:
\begin{equation}\label{eq:intro-game}
J_T^i(\bm\alpha)=
\Ex\!\left[\int_{t_0}^T
\bigl(L^i(X_t^{i,N},\alpha_t^i)+F^i(\bm X_t)\bigr)dt
+G^i(\bm X_T)\right].
\end{equation}
Here $\bm X=(X^{1,N},\ldots,X^{N,N})$, $\bm\alpha=(\alpha^1,\ldots,\alpha^N)$, the Brownian motions $W^1,\ldots,W^N$ are independent, and $L^i\colon \R^d\times \R^d\to \R,$ $F^i, G^i\colon (\R^d)^N\to \R$ are given running and terminal costs. The controls (or strategies) $\alpha^i$ are open-loop $\R^d$-valued processes, that is, they are adapted to the common filtration generated by the initial configuration $\bm X_{t_0}$ and all the Brownian motions. A feature of this choice of strategies is that, given a unilateral deviation of a player, the opponents' control processes remain fixed. A Nash equilibrium is a profile $\bm\alpha$ such that no player can lower his/her cost by a unilateral deviation. Unlike other related works, we do not require the players to be symmetric or to interact only through the empirical distribution of their states (mean field interaction).
\vspace{1mm}

Our assumptions on the costs combine joint strong convexity of $L^i$ in the state and control variables with monotonicity of the finite-player own-coordinate gradient maps $D\bm F(\bm x)=\bigl(D_{x^i}F^i(\bm x)\bigr)_{i=1}^N$, $D\bm G(\bm x)=\bigl(D_{x^i}G^i(\bm x)\bigr)_{i=1}^N$. This is the so-called displacement-monotone structure. It yields a quantitative comparison of finite-horizon equilibrium states and controls and allows us to expect that they exhibit long-time behavior on the noncompact configuration space $(\R^d)^N$.
\vspace{1mm}

The main objective of this paper is to show that these Nash equilibria satisfy the turnpike property with respect to a Nash equilibrium of an associated infinite-horizon average-cost (ergodic) $N$-player game. We then establish the long-time behavior of the finite-horizon $N$-player game: the cost of each player at equilibrium grows linearly with the horizon, at a rate given by his/her cost in the average-cost game, up to a correction that depends only on the initial configuration and an error that vanishes as the horizon grows. Our results are described below in Subsection \ref{subsec:intro-results} and stated precisely in Section \ref{not+as}.

\subsection{Background and related work}
Most of the existing results on the long-time behavior of stochastic differential games concern mean field games (MFGs for short), introduced by Lasry and Lions~\cite{LasryLions2007} and by Huang, Malham\'e, and Caines \cite{HuangMalhameCaines2006}. These are models for games with a very large number of identical players, each of whom interacts with the others only through the empirical distribution of their states\footnote{This would be the case in \eqref{eq:intro-game} if $F^i,G^i$ were functions that depend on the population state $\bm X_t$ only through the empirical measure $m^N_{\bm X_t}=\frac{1}{N}\sum_{i=1}^N\delta_{X^i_t}$.}. In the limit $N\to+\infty$ a single player has no influence on this distribution, so one considers a representative player who optimizes against a given evolution of the distribution of the population; at equilibrium, this evolution coincides with the distribution of the optimal state of the representative player. The equilibrium is described by two coupled PDEs: a backward Hamilton--Jacobi--Bellman equation for the value function of the representative player and a forward Fokker--Planck equation for the distribution of the population.
\vspace{1mm}

On the $d$-dimensional torus $\T^d:=\R^d/\mathbb{Z}^d$, Cardaliaguet, Lasry, Lions, and Porretta \cite{CardaliaguetLasryLionsPorretta2012,CardaliaguetLasryLionsPorretta2013} proved that the solution of this system converges, as $T\to+\infty$, to the solution of a stationary system, and Cardaliaguet and Porretta \cite{CardaliaguetPorretta2019} studied the long-time behavior of the so-called master equation, which describes the value of a player for every distribution of the population. These works rely on the monotonicity condition of Lasry and Lions \cite{LasryLions2007} on the dependence of the costs on the distribution; see \cite{CirantPorretta2021,CirantDeBernardi2025LocalTurnpike} for long-time results beyond this condition. On the whole space $\R^d$, where compactness is lost, exponential turnpike estimates are proved in \cite{CecchinConfortiDurmusEichinger2026} by probabilistic coupling, when the interaction is weak and the dynamics are confined, and confined MFGs with fractional diffusion are treated in \cite{ErslandJakobsenPorretta2025Fractional}.
\vspace{1mm}

The displacement monotonicity used in this paper is a different condition; in general, neither of the two implies the other \cite{GraberMeszaros2023}. For MFGs, it gives well-posedness on arbitrary time horizons \cite{GangboMeszarosMouZhang2022,meszaros2024mean,BansilMeszarosMou2025}. The work closest to ours is that of Cirant and M\'esz\'aros~\cite{cirant2024long}, who prove exponential long-time stability of the equilibria and describe the long-time behavior of the value function for displacement monotone MFGs on $\R^d$. Instead of the two PDEs above, they work with the first-order optimality conditions given by the Pontryagin maximum principle. We follow the same strategy for a fixed number of players, who may be different from each other, without passing to the limit $N\to+\infty$; see Subsection \ref{subsec:intro-proofs}.
\vspace{1mm}

For finitely many players, the information used by the players matters. Besides open-loop controls, one often considers closed-loop (or feedback) controls, which are functions of the current states of all the players, so that the opponents react to a deviation. The two notions of equilibrium differ in general, although not in an MFG, where a deviation of the representative player does not affect the population. Under displacement monotonicity, the convergence of the $N$-player equilibria to the MFG equilibrium as $N\to+\infty$ is proved in \cite{JacksonTangpi2024} for open-loop controls and in \cite{jackson2025quantitative} for both types of controls. Other works study games with many players directly, without a limiting model and without assuming that the players are identical: open-loop, closed-loop, and distributed equilibria (where each player uses only his/her own state) are compared in \cite{cirant2025non} when the interaction is weak, and it is shown in \cite{CirantRedaelli2026Distant} that, in open-loop games on sparse networks, a player is almost not affected by distant players, uniformly in the horizon. Our results are in the same spirit: $N$ is fixed, the players are heterogeneous, and the interaction is not assumed to be weak.
\vspace{1mm}

Average-cost games with finitely many players are classical: closed-loop equilibria are studied through systems of Hamilton--Jacobi--Bellman equations in \cite{BensoussanFrehse1995}, and for linear dynamics and quadratic costs in \cite{BardiPriuli2014}, while equilibria are constructed in \cite{CohenFedyashov2017} by means of the ergodic backward stochastic differential equations of \cite{FuhrmanHuTessitore2009}. The relation between long finite-horizon games and average-cost games for a fixed $N$ is less studied. Cohen and Jian~\cite{CohenJian2025Turnpike} prove turnpike estimates that are uniform in $N$ for games with linear dynamics and quadratic costs, with controls that use individual information, whereas we treat nonlinear costs and open-loop controls. The choice of the controls is important: Cardaliaguet and Rainer~\cite{CardaliaguetRainer2020} give an example of an ergodic game, played with closed-loop strategies, with many equilibrium payoffs in the limit $N\to+\infty$, although the MFG equilibrium is unique. Related infinite-horizon problems with discounted costs are studied in \cite{GraberMatterRosengartner2026Asymptotic,ContinelliEtAl2026DiscountTransform}.

\subsection{Description of the main results}\label{subsec:intro-results}
In the estimates below, $C$ and $c$ denote positive constants that do not depend on the horizon $T$ and may change from one estimate to another.
\vspace{1mm}

Our starting point is the Pontryagin maximum principle, that is, the first-order optimality conditions of the players. Let $H^i(x,p)=\sup_{a\in\R^d}\{-a\cdot p-L^i(x,a)\}$ be the Hamiltonian of player $i$, with gradients $D_xH^i$ and $D_pH^i$ in the two variables. We show that the game on $[t_0,T]$ has a unique Nash equilibrium, in which player $i$ uses the control $\alpha^i_t=-D_pH^i\bigl(X^{i,N}_t,v_T^i(t,\bm X_t)\bigr)$. Here $v_T=(v_T^1,\ldots,v_T^N)\colon[t_0,T]\times(\R^d)^N\to(\R^d)^N$ is a classical solution, unique in a suitable class, of the Pontryagin system of PDEs
\begin{equation}\label{eq:intro-PDE}
\begin{cases}
 -\partial_t v_T^i-\sum_{j=1}^N\operatorname{tr}\bigl(\sigma^j(\sigma^j)^\top D^2_{x^jx^j}v_T^i\bigr)+D_xH^i(x^i,v_T^i)\\[2pt]
 \hspace{2.5cm}+\sum_{j=1}^N D_{x^j}v_T^i\,D_pH^j(x^j, v_T^j)= D_{x^i}F^i(\bm x),&(t,\bm x)\in [t_0,T)\times(\R^d)^N,\\[2pt]
 v_T^i(T,\bm x)=D_{x^i}G^i(\bm x), &\bm x\in (\R^d)^N,
\end{cases}
\end{equation}
for $i\in\{1,\ldots,N\}$ (Propositions~\ref{prop:N-wellposed} and~\ref{N-PDE-classical}).
\vspace{1mm}

Our first result is that $v_T$ converges exponentially fast to a function $v_\infty\colon(\R^d)^N\to(\R^d)^N$:
\begin{equation}\label{eq:intro-field-turnpike}
|v_T(t,\bm x)-v_\infty(\bm x)|\leq C(1+|\bm x|)e^{-c(T-t)},\qquad\text{for all }t\in[t_0,T)\text{ and }\bm x\in(\R^d)^N
\end{equation}
(Theorems~\ref{thm:horizon-comparison} and~\ref{thm:stationary-limit}). The function $v_\infty$ does not depend on the terminal costs $G^i$, and it solves the stationary version of \eqref{eq:intro-PDE}, in which the time derivative and the terminal condition are dropped (see \eqref{N-PDEinf}). We associate with $v_\infty$ the state process $\bm X^\infty$ on $[t_0,+\infty)$ in which player $i$ uses the control $\alpha^i_\infty(t)=-D_pH^i(X^{i,\infty}_t,v^i_\infty(\bm X^\infty_t))$, and we write $\bm\alpha_\infty=(\alpha^1_\infty,\ldots,\alpha^N_\infty)$. The stationary system may have several solutions, and $v_\infty$ is the only Lipschitz one whose associated state process has second moments that are bounded uniformly in time (see Theorem~\ref{thm:stationary-limit} and the remark following it). As a consequence of \eqref{eq:intro-field-turnpike}, we obtain a turnpike estimate for the state $\bm X^T$ and the control $\bm\alpha_T$ at the equilibrium of the game on $[t_0,T]$. If $\bm X^T$ and $\bm X^\infty$ have the same initial configuration and the same Brownian motions, then, for every $R\in(t_0,T)$,
\begin{equation}\label{eq:intro-turnpike}
\Ex\Bigl[\sup_{t\in[t_0,R]}\bigl\{|\bm X_t^T-\bm X_t^\infty|^2+|\bm\alpha_T(t)-\bm\alpha_\infty(t)|^2\bigr\}\Bigr]\leq C_R\,e^{-2c(T-R)},
\end{equation}
where $C_R$ depends on $R$ and on the second moment of the initial configuration, but not on $T$ (Theorem~\ref{thm:turnpike}). If the bounds on the data are uniform in $N$ (Assumption \textup{(U)} in Section~\ref{not+as}), then \eqref{eq:intro-field-turnpike} and \eqref{eq:intro-turnpike} hold with constants independent of $N$, after a suitable normalization (see \eqref{goal:estimate2-uniform} and \eqref{eq:turnpike-fixed-window-N}).
\vspace{1mm}

Our second result gives a game-theoretic meaning to $\bm\alpha_\infty$. In the average-cost game, the players use open-loop controls on $[t_0,+\infty)$ and player $i$ minimizes the long-run average, per unit of time, of the running cost in \eqref{eq:intro-game}. We show that $\bm\alpha_\infty$ is a Nash equilibrium of this game, in the sense of Definition~\ref{erg-strong} (Theorem~\ref{prop:average-Nash}). Such equilibria are not unique (Remark~\ref{rem:erg-nonunique}), but the average cost of player $i$ is the same at all of them (Theorem~\ref{thm:main-erg}). We denote it by $\lambda^{i,N}$ and call it the ergodic constant of player $i$.
\vspace{1mm}

Our third result concerns the whole time interval $[t_0,T]$. The process $\bm X^\infty$ has a unique invariant probability measure $\mu^N_\infty$ on $(\R^d)^N$ with finite second moment, that is, a distribution that is preserved by its dynamics, and the law $\mathcal L(\bm X^T_s)$ of the equilibrium state started from a point $\bm x\in(\R^d)^N$ satisfies
\begin{equation}\label{eq:intro-law-turnpike}
W_2\bigl(\mathcal L(\bm X_s^{T}),\mu_\infty^N\bigr)\leq C(1+|\bm x|)\bigl(e^{-c(s-t_0)}+e^{-c(T-s)}\bigr),\qquad s\in[t_0,T],
\end{equation}
where $W_2$ is the quadratic Wasserstein distance between probability measures (Theorem~\ref{thm:mainLTMdistr}). The two terms come from the initial configuration and from the terminal cost, respectively. Moreover, if $V^{i,N}_T(t,\bm x)$ denotes the cost of player $i$ at the equilibrium of the game on $[t,T]$ started from $\bm x$, then
\begin{equation}\label{eq:intro-cost-expansion}
\bigl|V_T^{i,N}(t,\bm x)-(T-t)\lambda^{i,N}-w_G^{i,N}(\bm x)\bigr|\leq C(1+|\bm x|^2)e^{-c(T-t)}
\end{equation}
(Theorem~\ref{thm:mainLTB}). The function $w^{i,N}_G\colon(\R^d)^N\to\R$, which we call the corrector, is the unique solution, with at most quadratic growth and with a prescribed average with respect to $\mu^N_\infty$, of a linear elliptic equation associated with $\bm X^\infty$. It depends on the terminal costs only through an additive constant, as the subscript $G$ recalls.

\subsection{Proof techniques}\label{subsec:intro-proofs}
As in the work of Cirant and M\'esz\'aros \cite{cirant2024long}, who study MFGs under displacement monotonicity, we work with the first-order optimality conditions (Pontryagin systems of FBSDEs), but the main estimates are obtained differently. Our basic tool is to compare two solutions of the first-order optimality conditions (see \eqref{N-FBSDE}), with the same Brownian motions but different horizons, terminal costs, or initial configurations. By the convexity of $L^i$ and the monotonicity of $D\bm F$, the expected inner product of the difference of the two states with the difference of the corresponding processes~$v_T(t,\bm X_t)$ is nonincreasing in time, and its derivative controls the differences of the states and of the controls (see Lemma \ref{lem:pairing}). A computation of the same kind, for fully coupled forward--backward SDEs under monotonicity conditions, is also the core of the method developed in \cite{HuPeng1995,PengWu1999}. Adapting it to the present setting raises three difficulties.
\vspace{1mm}

We may compare solutions only if they are bounded in $L^2$ uniformly in the horizon. This is not automatically true, because the game is set on the unbounded space~$(\R^d)^N$, and the costs have quadratic growth. To resolve this, we prove a Lipschitz bound for $v_T$ and a bound for the second moments of the equilibrium, both independent of $T$ (Propositions~\ref{prop:uniform-Z} and~\ref{prop:uniform-moments}). The moment bound is obtained by comparing the equilibrium state with an Ornstein--Uhlenbeck process (see Proposition \ref{prop:uniform-moments}).
We expect that this approach can also work in the case of the Pontryagin system of FBSDEs for MFGs under displacement monotonicity.
Comparing two games with different horizons and sending the longer horizon to $+\infty$ then gives \eqref{eq:intro-field-turnpike}, and $v_\infty$ is obtained as a limit, without solving the stationary system in advance.
\vspace{1mm}

In an MFG or in a closed-loop game, the stationary problem is an ergodic Hamilton--Jacobi--Bellman equation for the value function of a player of an ergodic game, which contains the ergodic constant, and the equilibrium control is given by the gradient of its solution. In our open-loop setting, the stationary problem is the Pontryagin system \eqref{N-PDEinf}, and $v^i_\infty$ is in general not the gradient of the value of player $i$ (see Remark \ref{rem:erg-corrector-adjoint}), since this gradient also accounts for the dependence of the opponents' controls on the initial configuration; moreover, the equation for $v_\infty$ contains neither the ergodic constants nor a cost criterion. We therefore define the average-cost game directly, and we verify that $\bm\alpha_\infty$ is an equilibrium by a convexity argument that relies on the uniform bounds on the second moments. The ergodic constants and the correctors are constructed only afterwards, from the process~$\bm X^\infty$ and its invariant measure (see Propositions \ref{thm:erg-value} and \ref{prop:erg-poisson}).
\vspace{1mm}

Finally, the estimate~\eqref{eq:intro-turnpike} holds on fixed time windows, whereas \eqref{eq:intro-law-turnpike} and~\eqref{eq:intro-cost-expansion} require a constant that is independent of the observation time. The key point is that $v_\infty$ is strongly monotone, that is, $(v_\infty(\bm x)-v_\infty(\bm y))\cdot(\bm x-\bm y)\geq c|\bm x-\bm y|^2$ (Lemma \ref{lem:erg-strong-monotonicity}). Hence two copies of~$\bm X^\infty$ with the same Brownian motions approach each other exponentially fast, which gives the invariant measure (Proposition \ref{prop:erg-contraction}), and $\bm X^T$ can be viewed as a perturbation of~$\bm X^\infty$ by a drift that is exponentially small in the remaining time (Lemma \ref{lem:erg-pointwise-turnpike}). The expansion \eqref{eq:intro-cost-expansion} follows by splitting the horizon at an intermediate time $S$. On $[t_0,S]$, the equilibrium cost is close to the cost along $\bm X^\infty$, which equals $(S-t_0)\lambda^{i,N}$ plus an increment of the corrector. At time $S$, the law of the equilibrium state is close to $\mu^N_\infty$, so the remaining cost is close to its $\mu^N_\infty$-average. Choosing $S$ at the midpoint gives the exponential rate.

\subsection{Organization of the paper} The rest of the paper is organized as follows. Section~\ref{sec:model} gives the rigorous model, the Nash equilibrium
definition, and its Pontryagin formulation.
Section~\ref{not+as} states the assumptions and main results.
Section~\ref{3} establishes the well-posedness of the finite-horizon Pontryagin systems of FBSDEs and PDEs. In Section \ref{5}, we show the turnpike estimates and we construct the solution of the infinite-horizon Pontryagin system of PDEs. In Section \ref{6}, we prove that the infinite-horizon Pontryagin system of PDEs is associated with an equilibrium of an average-cost ergodic $N$-player game. In Section \ref{7}, we establish the long-time behavior of the finite-horizon $N$-player game. Technical proofs and further characterizations of equilibria for the average-cost game, which might be of independent interest, are collected in the Appendices \ref{Appendix2} and \ref{char}.

\section{Model formulation and notation}\label{sec:model}

\subsection{Setup of the problem} Consider a complete filtered probability space $(\Omega,\mathcal{F},\Pr)$, which satisfies the usual conditions and supports $N$ independent $d$-dimensional standard Brownian motions $\{W^1_t,...,W^N_t\}_{t\ge 0}$. We denote the expectation on this space by $\Ex$.
 Fix an initial time $t_0\geq0$. We use the usual augmentation of the joint filtration
$$
\mathcal F_t=\mathcal I\vee\sigma(W_s^j-W_{t_0}^j:t_0\leq s\leq t,\ j\in[N]),\qquad t\geq t_0,
$$
where the initial sigma-algebra $\mathcal I$ is independent of all Brownian increments after $t_0$. The square-integrable initial vector $\bm\xi_0$ is $\mathcal I$-measurable. All players observe this common filtration, including the initial vector; all controls and solutions below use this filtration. We consider the following system of Stochastic Differential Equations (SDEs for short):
\be\label{N-pl}
\begin{cases}
dX^{i,N}_t=\alpha^i_t\;dt+\sqrt{2}\sigma^idW^i_t,\;\;\;t\ge t_0,\\
    X^{i,N}_{t_0}=\xi_0^i\in L^2(\Omega;\R^d),
\end{cases}
\ee
for $i\in [N]$, where $\sigma^i\in \R^{d\times d}$ is assumed to be invertible. In \eqref{N-pl}, the process $X^i=(X^i_t)_{t\ge t_0}$ represents the state of player $i$ and the process $\alpha^i=(\alpha^i_t)_{t\ge t_0}$ is the corresponding control process, also referred to as the strategy of player $i$. We use the notation $\bm X_t=(X^{1,N}_t,\ldots,X^{N,N}_t)\in (\R^{d})^N$, for $t\ge t_0$, therefore we also denote $\bm X_{t_0}=\bm \xi_0=(\xi^1_0,\ldots,\xi^N_0)\in L^2(\Omega;(\R^d)^N)$. Similarly, we denote by $\bm \alpha=(\alpha^1,\ldots,\alpha^N)$ the strategy profile of the whole population of players.
\vspace{2mm}

Let $T\ge t_0$. We consider the finite-horizon $N$-player game ($t\in [t_0,T]$), where players minimize their expected total cost. In particular, we look at the problem where the $i$-th player aims to minimize the functional
\be\label{cost2}
\begin{aligned}
J^i_T(t_0,\bm\xi_0;\bm\alpha):=\Ex\left[\int_{t_0}^T\bigl(L^i(X_t^{i,N},\alpha_t^i)+F^i(\bm X_t)\bigr)dt+G^i(\bm X_T)\right],
\end{aligned}
\ee
where $L^i\colon \R^d\times\R^d\to \R,\; F^i\colon (\R^d)^N\to \R$ and $G^i\colon (\R^d)^N\to \R$, for $i=1,\ldots,N$, are given Lagrangians, running costs and terminal costs, respectively. We assume that the game is played in the open-loop sense; that is, for each $i\in [N]$, player $i$ chooses a control $\alpha^i$ from the set
\begin{align*}
\mathcal A^i[t_0,T]:=\left\{\begin{aligned}
&\alpha:\Omega\times[t_0,T]\to\R^d:\ \alpha\text{ is }\\
&\qquad(\mathcal F_t)_{t_0\leq t\leq T}\text{--progressively measurable, and }\Ex\left[\int_{t_0}^T|\alpha_t|^2dt\right]<\infty
\end{aligned}\right\}.
\end{align*}

We will refer to this game as the \textit{ finite-horizon $N$-player game} (or simply \textit{$N$-player game} when there is no potential confusion) and we are interested in equilibrium states of this game. In the following definition, for $\bm y\in (\R^{d})^N$ and $i\in [N]$, we use the notation $\bm y^{-i}=(y_1,\ldots,y_{i-1},y_{i+1},\ldots,y_N)$ for the vector $\bm y$ with its $i$ component missing. We also write $[\bm y^{-i};y^i]=\bm y$.

\begin{definition}[Nash equilibrium in the finite-horizon game]
\label{d:Nash_equilibrium_finite}
Fix $t_0<T$ and a square-integrable $\mathcal I$-measurable initial state vector $\bm\xi_0=(\xi_0^1,\ldots,\xi_0^N)$. We say that a strategy profile $\boldsymbol{\alpha} = (\alpha^{1}, \ldots, \alpha^{N})\in\prod_{j=1}^N\mathcal{A}^j[t_0,T]$ is an
{\rm
open-loop Nash equilibrium in the finite-horizon game} if for every $i \in [N]$ and every $\alpha \in \mathcal{A}^i[t_0,T]$,
$$J^i_T(t_0,\bm\xi_0;[\boldsymbol\alpha^{-i};\alpha])\geq J^i_T(t_0,\bm\xi_0;\boldsymbol\alpha).$$
\end{definition}

The characterization of open-loop equilibria of this game involves the Hamiltonians $H^i\colon\R^d\times\R^d\to \R$, where
\be\label{Hamiltonian}
H^i(x,p)=\sup_{\alpha\in \R^d}\{ -\alpha\cdot p -L^i(x,\alpha)\},\;\;\text{for }i\in [N].
\ee
We recall that the Pontryagin maximum principle implies that, under appropriate assumptions, any open-loop Nash equilibrium of the $N$-player game takes the form
\be\label{Topt}
\alpha^i_T(t)=-D_pH^i(X^{i,N}_t,Y^{i,N}_t),
\ee
where the tuple
$$(\bm X^N, \bm Y^N,\bm Z^N)=\left((X^{1,N},\ldots,X^{N,N}), (Y^{1,N},\ldots,Y^{N,N}),( Z^{1,1,N},Z^{1,2,N},\ldots, Z^{N,N-1,N},Z^{N,N,N})\right)$$
is a solution\footnote{We refer to Definition~\ref{def:N-FBSDE-solution} below for a precise definition of a solution of \eqref{N-FBSDE}.} to the following Pontryagin system of forward-backward stochastic differential equations (FBSDEs for short)
\be\label{N-FBSDE}
\begin{cases}
dX_t^{i,N}=-D_pH^i(X_t^{i,N},Y_t^{i,N})dt+\sqrt{2}\sigma^idW^i_t,\\[2pt]
dY_t^{i,N}=(D_xH^i(X_t^{i,N},Y_t^{i,N})-D_{x^i}F^i(\bm X_t))dt+\sqrt{2}\sum_{j=1}^NZ_t^{i,j,N}dW^j_t,\\[2pt]
X^{i,N}_{t_0}=\xi_0^i,\;\; Y^{i,N}_T=D_{x^i}G^i(\bm X_T),\qquad t\in[t_0,T].
\end{cases}
\ee
Under additional displacement monotonicity assumptions (see (A1)--(A2) from Section \ref{not+as}), we show in Section \ref{3} that the above system of FBSDEs has a unique solution in the sense of the following definition:

\begin{definition}\label{def:N-FBSDE-solution}
Let $T>t_0$. We say that the tuple $(\bm X^N,\bm Y^N,\bm Z^N)$ satisfies \eqref{N-FBSDE} if
\begin{enumerate}
    \item $X^{i,N},\; Y^{i,N}$ are continuous and adapted $\R^d$-valued processes for each $i\in[N]$;
    \item $Z^{i,j,N}$, for $i,j\in[N]$, are progressively measurable $\R^{d\times d}$-valued processes;
    \item $(\bm X^N,\bm Y^N,\bm Z^N)$ satisfies the condition
    $$\mathbb{E}\left[ \sup_{t_0\le t\le T}\left(|\bm X_t^N|^2+|\bm Y_t^N|^2\right)+\int_{t_0}^T|\bm Z_t^N|^2dt \right]<\infty;$$
    \item $(\bm X^N,\bm Y^N,\bm Z^N)$ satisfies the integral form of \eqref{N-FBSDE}.
\end{enumerate}
\end{definition}
\vspace{2mm}

Furthermore, the Pontryagin system of FBSDEs \eqref{N-FBSDE} can be decoupled by the Pontryagin system of partial differential equations (PDEs for short):
\be\label{N-PDE}
\begin{cases}
    -\partial_t v_T^i-\sum_{j=1}^N\text{tr}(\sigma^j(\sigma^j)^\top D^2_{x^jx^j}v_T^i)+D_xH^i(x^i,v_T^i)\\
    \hspace{2.5cm}+\sum_{j=1}^N D_{x^j}v_T^i D_pH^j(x^j, v_T^j)= D_{x^i}F^i(\bm x),&(t,\bm x)\in [t_0,T)\times(\R^d)^N,\\
    v_T^i(T,\bm x)=D_{x^i}G^i(\bm x), &\bm x\in (\R^d)^N.
\end{cases}
\ee

In the above system the unknown map $v_T=(v^1_T,\ldots,v_T^N)\colon [t_0,T]\times(\R^d)^N\to (\R^d)^N$ satisfies \eqref{N-PDE} in the classical sense\footnote{A classical solution is continuous on
$[t_0,T]\times(\R^d)^N$, is once continuously differentiable in
time and twice continuously differentiable in space on
$[t_0,T)\times(\R^d)^N$, satisfies the PDE there, and attains the
terminal condition continuously at $T$.
} and its association with \eqref{N-FBSDE} is through the following formulas
$$Y_t^{i,N}=v_T^i(t,\bm X_t^N)\quad\text{and}\quad Z_t^{i,j,N}=D_{x^j}v_T^i(t,\bm X_t^N)\sigma^j,\quad\text{for all } i,j\in[N],\;dt\otimes d\Pr\text{-a.e}.$$

\subsection{Notation}
We write $[N]=\{1,\ldots,N\}$ and identify $(\R^d)^N$ with
$\R^{Nd}$. For $\bm x=(x^1,\ldots,x^N)$, we use $|\bm x|^2=\sum_{i=1}^N|x^i|^2$ for the standard Euclidean norm. For comparisons across different population sizes, we use the
population-normalized inner product and norm
\[
\langle\bm x,\bm y\rangle_N
:=
\frac1N\bm x\cdot\bm y,
\qquad
|\bm x|_N
:=
\langle\bm x,\bm x\rangle_N^{1/2}
=
\frac{|\bm x|}{\sqrt N}.
\]
$I_d$ denotes the $d\times d$ identity matrix. For a matrix $A=(A_{k\ell})\in\R^{q\times m}$ we introduce its operator and Hilbert--Schmidt norms
\[\|A\|_{\mathrm{op}}:=\sup_{h\neq0}\frac{|Ah|}{|h|}=\sup_{|h|=1}|Ah|\quad\text{and}\quad \|A\|_{\mathrm{HS}}
:=
\left(
\sum_{k=1}^q\sum_{\ell=1}^m|A_{k\ell}|^2
\right)^{1/2}
=
\sqrt{\operatorname{tr}(AA^\top)}.\]
If $A$ $\in\R^{Nd\times Nd}$ is partitioned into $N\times N$ blocks $A^{i,j}\in\R^{d\times d}$, we use the notation $\|A\|_{\mathrm{HS},N}^2:=\frac1N\sum_{i,j=1}^N\|A^{i,j}\|_{\mathrm{HS}}^2=\frac1N\|A\|_{\mathrm{HS}}^2.$
In particular,
$\|D_{\bm x}v_T(t,\bm x)\|_{\mathrm{op}}$ denotes the operator norm
of the full $Nd\times Nd$ Jacobian of $v_T(t,\cdot)$.\\
For a map $w:\R^{Nd}\to\R^{Nd}$ with at most linear growth, define
\[
\|w\|_*
:=
\sup_{\bm x\in\R^{Nd}}
\frac{|w(\bm x)|}{1+|\bm x|},
\qquad
\|w\|_{*,N}
:=
\sup_{\bm x\in\R^{Nd}}
\frac{|w(\bm x)|_N}{1+|\bm x|_N}
=
\sup_{\bm x\in\R^{Nd}}
\frac{|w(\bm x)|}{\sqrt N+|\bm x|}.
\]
For every fixed $N$, $\|w\|_{*,N}\leq\|w\|_*\leq\sqrt N\,\|w\|_{*,N}.$ Thus the two norms are equivalent for a fixed game, whereas $\|\cdot\|_{*,N}$ is normalized for estimates that are uniform in the number of players.

When no confusion is possible, we suppress the superscript $N$ on $\bm X^N,\bm Y^N,\bm Z^N$ in proofs.

For subsets $A$ and $B$ of a Euclidean space, we write
$
A\Subset B
$
if $\overline A$ is compact and contained in $B$; equivalently, $A$ is
relatively compact in $B$. When $B$ is open, this notation in particular
implies that $A$ remains a positive distance from $\partial B$. By contrast,
$A\subset B$ denotes ordinary inclusion.

We denote by $\mathcal P_2(\R^d)$ the set of probability measures on
$\R^d$ with finite second moment, by $\Pi(m,\bar m)$ the set of
couplings of $m$ and $\bar m$, and by
\[
W_p(m,\bar m)
:=
\inf_{\pi\in\Pi(m,\bar m)}\left(
\int_{\R^d\times\R^d}|z-\bar z|^p\,\pi(dz,d\bar z)\right)^{1/p},\;\;\;\;p\in \{1,2\}
\]
their 1-Wasserstein and 2-Wasserstein distances, respectively.
\vspace{0.5mm}

For a sufficiently regular function
$\mathcal U:\R^d\times\mathcal P_2(\R^d)\to\R$,
the vector
$D_m\mathcal U(x,m,z)\in\R^d$ denotes its Lions derivative with
respect to the measure variable. Thus, if
$X,V\in L^2(\Omega;\R^d)$ and $m=\mathcal L(X)$, then
\[
\left.
\frac{d}{d\varepsilon}
\right|_{\varepsilon=0}
\mathcal U\bigl(x,\mathcal L(X+\varepsilon V)\bigr)
=
\Ex\left[
D_m\mathcal U(x,m,X)\cdot V
\right].
\]

\section{Assumptions and Main Results}\label{not+as}

\subsection{Assumptions} Throughout the paper, we will be using the following assumptions on the data $L^i(x,a)=L^i\colon \R^d\times\R^d\to \R$, $F^i(\bm x)=F^i\colon (\R^d)^N\to \R$ and $G^i(\bm x)=G^i\colon (\R^d)^N\to \R$. To lighten the notation, we denote $D\bm F(\bm x)=\bigl(D_{x^i}F^i(\bm x)\bigr)_{i=1}^N$ and $D\bm G(\bm x)=\bigl(D_{x^i}G^i(\bm x)\bigr)_{i=1}^N.$
\vspace{3mm}

\textbf{(A1).} For every $i\in [N]$, $L^i,\; F^i$ and $G^i$ are $C^2$ functions with bounded derivatives of order 2 and satisfy the following displacement monotonicity assumption
$$D^2_{(x,a)}L^i\ge C_L\begin{pmatrix}
I_d & 0 \\
0 & I_d
\end{pmatrix},\;\; \left( D^2_{x^jx^i}F^i\right)_{i,j=1,...,N}\ge 0\text{ and }\left( D^2_{x^jx^i}G^i\right)_{i,j=1,...,N}\ge 0,$$
where $C_L$ is a given positive constant. Here $D^2_{(x,a)}L^i$ is the full $2d\times2d$ Hessian. The inequalities for the coupling Jacobians mean nonnegativity of their quadratic forms, equivalently of their symmetric parts; symmetry of those Jacobians is not assumed (see also
\cite[p.~19, equations~(2.6)--(2.7)]
{cirant2025non}). Here, $\left( D^2_{x^jx^i}F^i\right)_{i,j=1,...,N}$ and $\left( D^2_{x^jx^i}G^i\right)_{i,j=1,...,N}$ are viewed as elements of $(\R^{d\times d})^{N\times N}$. Equivalently, for any $\overline{\bm x}, \bm x=(x^1,\ldots,x^N)\in (\R^d)^N$ and $x,a,\overline{x},\overline{a}\in \R^d$ we assume:
\begin{align}
   (D_xL^i(x,a)-D_xL^i(\overline{x},\overline{a}))\cdot (x-\overline{x})+(D_aL^i(x,a)-D_aL^i(\overline{x},\overline{a}))&\cdot (a-\overline{a})\nonumber\\
   &\ge C_L(|x-\overline{x}|^2+|a-\overline{a}|^2);\label{Ldispl}\\
   \sum_{j=1}^N(D_{x^j}F^j(\bm x)-D_{x^j}F^j(\overline{\bm x}))\cdot ( x^j-\overline{ x}^j)&\ge 0;\label{Fmon}\\
   \sum_{j=1}^N(D_{x^j}G^j(\bm x)-D_{x^j}G^j(\overline{\bm x}))\cdot ( x^j-\overline{ x}^j)&\ge 0.\label{Gmon}
\end{align}
\vspace{1.5mm}

\textbf{(A2).} The derivatives of $F^i,\; G^i$, $i\in[N]$, satisfy the following Lipschitz property for some constant $C_{\text{Lip}}\ge 0$
$$\sum_{i=1}^N|D_{x^i}F^i(\bm x)-D_{x^i}F^i(\overline{\bm x})|^2+\sum_{i=1}^N|D_{x^i}G^i(\bm x)-D_{x^i}G^i(\overline{\bm x})|^2\le C_{\text{Lip}}|\bm x- \overline{\bm x}|^2,\;\text{for any }\bm x,\overline{\bm x}\in (\R^d)^N.$$

\textbf{(A3).} The matrices $\sigma^j$ are constant and satisfy $\sigma^j(\sigma^j)^\top\geq \nu I_d,\; j=1,\ldots,N,$ for some $\nu>0$.
\vspace{4mm}

By adapting standard proofs (e.g. \cite[Remark 2.6]{meszaros2024mean}), we can show that \eqref{Ldispl} gives a displacement monotonicity condition for the Hamiltonians $H^i$ defined in \eqref{Hamiltonian}:
\be\label{Hdispl}
\begin{split}
(-D_xH^i(x,p)+D_xH^i(\overline{x},\overline{p}))\cdot (x-\overline{x})+&(D_pH^i(x,p)-D_pH^i(\overline{x},\overline{p}))\cdot (p-\overline{p})\\
&\ge C_L(|x-\overline{x}|^2+|D_pH^i(x,p)-D_pH^i(\overline{x},\overline{p})|^2),
\end{split}
\ee
for any $x,\overline{x},p,\overline{p}\in \R^d$. We now let $K_L:=\max_{i\in[N]}\sup_{(x,a)\in\R^{2d}}\|D^2_{(x,a)}L^i(x,a)\|_{\mathrm{op}}$, which is finite by \textup{(A1)}. Using \eqref{Hdispl} and the standard equality $D_xH^i(x,p)=-D_xL^i(x,-D_pH^i(x,p))$, we may also show that
\begin{equation}\label{eq:H-gradient-Lipschitz}
\left|
D_pH^i(x,p)-D_pH^i(\bar x,\bar p)
\right|
\leq
\frac{K_L}{C_L}|x-\bar x|
+\frac1{C_L}|p-\bar p|,
\end{equation}
and
\be\label{eq:Hx-grad}
\left|
D_xH^i(x,p)-D_xH^i(\bar x,\bar p)
\right|
\leq
K_L\left(1+\frac{K_L}{C_L}\right)|x-\bar x|
+\frac{K_L}{C_L}|p-\bar p|,
\ee
for any $x,\overline{x},p,\overline{p}\in \R^d$.

\vspace{3mm}

We note that the constants appearing in Assumptions (A1)-(A3) may depend on $N$. In order to obtain estimates that are independent of $N$ (uniformity in the number of players), we introduce the following assumption:
\vspace{2mm}

\textbf{(U).} For every $N$, Assumptions \textup{(A1)--(A3)} hold and $\nu$ from (A3) is independent of $N$. There exist
constants
$0<\overline C_L,\overline K_L,\overline L_F,\overline L_G,B_0$,
independent of $N$, such that the following bounds hold. For every $i\in[N]$ and $(x,a)\in\R^{2d}$,
\[
D^2_{(x,a)}L^{i}(x,a)\geq
\overline C_L I_{2d},
\qquad
\left\|D^2_{(x,a)}L^{i}(x,a)\right\|_{\mathrm{op}}
\leq\overline K_L.
\]
Moreover, for all $\bm x,\bar{\bm x}\in(\R^{d})^N$,
\[|D\bm F(\bm x)-D\bm F(\bar{\bm x})|\leq\overline L_F|\bm x-\bar{\bm x}|\;\;\text{ and }\;\;|D\bm G(\bm x)-D\bm G(\bar{\bm x})|\le\overline L_G|\bm x-\bar{\bm x}|.\]
Finally,
\[
\begin{aligned}
\sup_{N\geq1}\frac1N\Bigg[
\sum_{i=1}^N\left(
 \|\sigma^{i}\|_{\mathrm{HS}}^2
 +|D_aL^{i}(0,0)|^2
 +|D_xL^{i}(0,0)|^2
\right)
+|D\bm F(\bm0)|^2
+|D\bm G(\bm0)|^2
\Bigg]\leq B_0.
\end{aligned}
\]

\begin{remark}
Under \textup{(U)}, the Lipschitz constants in \eqref{eq:H-gradient-Lipschitz} and \eqref{eq:Hx-grad} are independent of $N$ and $i$.
\end{remark}

\subsubsection{Examples}
We have two main classes of examples of costs $L^i$, $F^i, G^i$ in mind:\\[1ex]
(1) (Displacement convex costs). For any $i\in [N]$, set
\[ L^i(x, a) = \frac{1}{2} a^\top R a + \frac{1}{2} x^\top Q_L x,\;\; F^i(\bm x) = \frac{\theta}{2N} \sum_{j=1}^N |x^i - x^j|^2,\;\; G^i(\bm x) = \frac{1}{2} |x^i - \xi^i|^2. \]
Here $R\geq C_LI_d$ and $Q_L\geq C_LI_d$, as required by joint strong convexity in \textup{(A1)}; $\theta>0$, and the targets $\xi^i\in\R^d$ are fixed deterministic parameters, distinct from the random initial states $\xi_0^i\in L^2$. This example satisfies (A1)-(A2).\\[1ex]
(2) (Mean field case). For any $i\in [N],$ let
\[
L^i=L,\qquad
F^i(\bm x)=\mathcal F(x^i,m^N_{\bm x}),\qquad
G^i(\bm x)=\mathcal G(x^i,m^N_{\bm x}),
\]
where $L\colon\R^d\times\R^d\to\R$,
$\mathcal F,\mathcal G\colon\R^d\times\mathcal P_2(\R^d)\to\R$
are sufficiently regular with bounded derivatives and $m^N_{\bm x}=\frac1N\sum_{j=1}^N\delta_{x^j}$ is the empirical measure. We assume that $L$ satisfies \eqref{Ldispl} and that for $\mathcal U\in\{\mathcal F,\mathcal G\}$, the Lions derivative $D_m\mathcal U(x,m,z)$ and the derivative $D_x\mathcal{U}(x,m)$ exist and are globally $\ell_{\mathcal U}$-Lipschitz with respect to all variables. We also assume that $\mathcal U\in \{\mathcal{F},\mathcal{G}\}$ is strongly
displacement monotone: there is $c_{\mathcal U}>0$ such that, for every
$m,n\in\mathcal P_2(\R^d)$ and every $X,Y\in L^2$ with $X\sim m,\; Y\sim n$
\[\mathbb{E}
\bigl[
(D_x\mathcal U(X,m)-D_x\mathcal U(Y,n))\cdot(X-Y)\bigr]
\geq
c_{\mathcal U}
\mathbb{E}[|X-Y|^2].
\]
We have
\[D_{x^i}F^i(\bm x)=D_x\mathcal F(x^i,m^N_{\bm x})+\frac1N D_m\mathcal F(x^i,m^N_{\bm x},x^i).\]
Using $\pi=\frac{1}{N}\sum_{i=1}^N\delta_{(x^i,\overline{x}^i)}$,  $W_2(m^N_{\bm x},m^N_{\bar{\bm x}})
\leq N^{-1/2}|\bm x-\bar{\bm x}|$ and Cauchy--Schwarz, we obtain
\[
\begin{aligned}
&\sum_{i=1}^N
\bigl(
D_{x^i}F^i(\bm x)
-D_{x^i}F^i(\bar{\bm x})
\bigr)\cdot(x^i-\bar x^i)\geq
\left(c_{\mathcal U}-\frac{3\ell_{\mathcal U}}N\right)
|\bm x-\bar{\bm x}|^2.
\end{aligned}
\]
Thus, \eqref{Fmon} is satisfied whenever $N\geq N_0:=
\left\lceil
\max\left\{
\frac{3\ell_{\mathcal F}}{c_{\mathcal F}},
\frac{3\ell_{\mathcal G}}{c_{\mathcal G}}
\right\}
\right\rceil$ and similarly for \eqref{Gmon}. Assumption (A1) is satisfied if $L^i,F^i,G^i$ have bounded second derivatives. Assumption (A2) follows from the Lipschitz continuity of $D_x\mathcal{U},D_m\mathcal{U}$, for $\mathcal{U}\in \{\mathcal{F},\mathcal{G}\}$.

We refer to the works \cite{meszaros2024mean,cirant2024long}, where displacement monotonicity assumptions for $\mathcal{F},\mathcal{G}$ appear in the context of MFGs, and to \cite{jackson2025quantitative,cirant2025non} in $N$-player games.

\subsection{Main results}  An important consequence of our displacement monotonicity assumptions (see Section~\ref{not+as}) is that an open-loop Nash equilibrium exists and the association between open-loop Nash equilibria and \eqref{Topt} is an equivalence (Proposition \ref{prop:N-wellposed}). Furthermore, we show that the Pontryagin system of PDEs \eqref{N-PDE} has a unique classical solution $v_T$, which is globally Lipschitz in space with the Lipschitz constant being independent of the time horizon $T$ (Propositions \ref{N-PDE-classical} and \ref{prop:uniform-Z}).

\subsubsection{Turnpike estimates}
Using this as our starting point, the purpose of this paper is to obtain turnpike estimates for the $N$-player game. We first obtain an estimate that illustrates how the solution of \eqref{N-PDE} changes if we change the time horizon $T$.

\begin{theorem}[Comparison of terminal costs and finite horizons]\label{thm:horizon-comparison}
Assume that the running costs $L^i,F^i$ and the diffusion matrices $\sigma^i$, $i\in [N]$, satisfy Assumptions \textup{(A1)-(A3)}. Also, let $\bm G,\widetilde{\bm G}\colon (\R^d)^N\to \R^N$ be terminal costs satisfying Assumptions \textup{(A1)-(A2)} and $T'\ge T\ge t_0$ be two terminal times. We consider $v_T^G\colon [t_0,T]\times(\R^d)^N\to (\R^d)^N$ and $v_{T'}^{\widetilde{G}}\colon [t_0,T]\times(\R^d)^N\to (\R^d)^N$ the solutions of \eqref{N-PDE} with terminal times $T,T'$ and terminal data $D\bm G,D\widetilde{\bm G}$, respectively. Then, there exist constants $c,C>0$ that are independent of $T,T'$, but might depend on $N$ such that
\begin{equation}\label{goal:estimate1}
\left\|
v_T^G(t,\cdot)-v_{T'}^{\widetilde G}(t,\cdot)
\right\|_*
=
\sup_{\bm x\in\R^{Nd}}
\frac{|v_T^G(t,\bm x)-v_{T'}^{\widetilde G}(t,\bm x)|}{1+|\bm x|}\leq Ce^{-c(T-t)},\;\;t\in [t_0,T].
\end{equation}
If, in addition, Assumption \textup{(U)} is in place, then there are constants $c_*,C_*$ independent of $N$ such that
\be\label{goal:unifNest}
\left\|
v_T^G(t,\cdot)-v_{T'}^{\widetilde G}(t,\cdot)
\right\|_{*,N}
=
\sup_{\bm x\in\R^{Nd}}
\frac{|v_T^G(t,\bm x)-v_{T'}^{\widetilde G}(t,\bm x)|}{\sqrt{N}+|\bm x|}\leq C_*e^{-c_*(T-t)},\;\;t\in [t_0,T].
\ee
\end{theorem}
\vspace{4mm}

Theorem~\ref{thm:horizon-comparison} is proved in Subsection~\ref{pf:horizon-comparison}. The above estimate \eqref{goal:estimate1} allows us to characterize the long-time behavior of $v_T$ through a solution $v_{\infty}:=(v_{\infty}^1,\ldots,v_{\infty}^N)\colon (\R^d)^N\to (\R^{d})^N$ of the stationary Pontryagin system of PDEs:
\be\label{N-PDEinf}
\begin{cases}
-\sum_{j=1}^N\operatorname{tr}\left(
\sigma^j(\sigma^j)^\top D^2_{x^jx^j}v_\infty^i\right)
+D_xH^i(x^i,v_\infty^i)\\[2ex]
\hspace{3cm}+\sum_{j=1}^ND_{x^j}v_\infty^iD_pH^j(x^j,v_\infty^j)
=D_{x^i}F^i(\bm x),\quad i\in[N].
\end{cases}
\ee

For any globally Lipschitz solution $w=(w_1,\ldots,w_N)\colon (\R^d)^N\to (\R^d)^N$ of \eqref{N-PDEinf}, we consider its induced forward diffusion $\bm X^w$ so that
\be\label{forward''}
dX_s^{i,w}=-D_pH^i(X_s^{i,w},w^i(\bm X_s^w))\,ds
+\sqrt{2}\sigma^i\,dW_s^i,\;\;s\ge t_0,\qquad \bm X_{t_0}^w=\bm \xi_0.
\ee

Our main turnpike result is stated as follows.

\begin{theorem}[The limiting Pontryagin field]\label{thm:stationary-limit}
Assume \textup{(A1)--(A3)}. Let $T>t_0$ be a terminal time and let $C,c$ be the fixed constants obtained in
Theorem~\ref{thm:horizon-comparison} for
$\widetilde{\bm G}=\bm G$. Let $v_T$ be the unique solution to \eqref{N-PDE}. There exists a globally Lipschitz classical solution
$v_\infty\colon(\R^d)^N\to(\R^d)^N$ of \eqref{N-PDEinf} with at most linear growth. It is unique among such globally Lipschitz solutions $w$ of \eqref{N-PDEinf} for which the forward diffusion \eqref{forward''} induced by $w$ satisfies
\be\label{condition}
\sup_{s\geq t_0}\Ex\left[|\bm X_s^w|^2\right]<\infty,\text{ whenever }\bm \xi_0=\bm x\in (\R^d)^N.
\ee
In addition, the following estimate holds
\begin{equation}\label{goal:estimate2}
\left\|v_T(t,\cdot)-v_\infty(\cdot)\right\|_*
=
\sup_{\bm x\in\R^{Nd}}
\frac{|v_T(t,\bm x)-v_\infty(\bm x)|}{1+|\bm x|}
\leq
C
e^{-c(T-t)},\qquad t\in[t_0,T).
\end{equation}
Furthermore, if Assumption \textup{(U)} is in place, for the constants $C_*,c_*$ from Theorem \ref{thm:horizon-comparison}, we have
\begin{equation}\label{goal:estimate2-uniform}
\left\|v_T(t,\cdot)-v_\infty(\cdot)\right\|_{*,N}
=
\sup_{\bm x\in\R^{Nd}}
\frac{
 |v_T(t,\bm x)-v_\infty(\bm x)|
}
{\sqrt N+|\bm x|}
\leq
C_*e^{-c_*(T-t)},\qquad t\in[t_0,T).
\end{equation}
\end{theorem}
\vspace{2mm}

Theorem~\ref{thm:stationary-limit} is proved in Subsection~\ref{pf:stationary-limit}. The next remark explains why the stationary PDE does not have a unique globally Lipschitz, linearly growing solution unless one imposes a stability condition (e.g. \eqref{condition}) on the state process generated by the solution.
\begin{remark}[Why \eqref{condition} is needed]
The restriction \eqref{condition} in the uniqueness statement is essential.
Take $t_0=0$, $N=d=1$, $L(x,a)=(a^2+x^2)/2$, $F=G=0$, and $\sigma>0$. Then, $H(x,p)=\frac{p^2-x^2}{2}$ and \eqref{N-PDE} becomes
$$-\partial_tv_T(t,x)-\sigma^2v_T''(t,x)+v_T'(t,x)v_T(t,x)=x,\qquad(t,x)\in [t_0,T]\times\R,$$
therefore $v_T(t,x)=x\tanh(T-t)\xrightarrow{T\to +\infty}x=v_{\infty}(x)$. Both $v_+(x)=x$ and $v_-(x)=-x$ are globally Lipschitz solutions of
the stationary Pontryagin equation
\[
-\sigma^2v''(x)-x+v'(x)v(x)=0.
\]
These functions induce different forward dynamics. The field $v_+(x)=x$ gives
\[dX_s=-X_s\,ds+\sqrt{2}\sigma\,dW_s,\qquad\Ex[|X_s|^2]=e^{-2s}|x|^2+\sigma^2(1-e^{-2s}),\]
and therefore has uniformly bounded second moments. In contrast,
$v_-(x)=-x$ gives
\[dX_s=X_s\,ds+\sqrt{2}\sigma\,dW_s,\qquad\Ex[|X_s|^2]=e^{2s}|x|^2+\sigma^2(e^{2s}-1),\]
whose second moments are unbounded. Thus $v_-$ lies outside the stability class, whereas the finite-horizon fields converge to the stable solution $v_\infty=v_+$.
\end{remark}
\vspace{2mm}

A corollary of the proof of Theorem \ref{thm:stationary-limit} is a turnpike estimate between a solution of \eqref{N-FBSDE} and the forward problem induced by $v_\infty$:
\be\label{forward'}
dX_s^{i,\infty}=-D_pH^i(X_s^{i,\infty},v_\infty^i(\bm X_s^\infty))\,ds
+\sqrt{2}\sigma^i\,dW_s^i,\;\;s\ge t_0,\qquad \bm X_{t_0}^\infty=\bm \xi_0,
\ee
where $\bm \xi_0$ is the common initial condition. We state this estimate in the following theorem.

\begin{theorem}[Turnpike estimates]\label{thm:turnpike}
    Suppose that Assumptions \textup{(A1)-(A3)} hold and $N$ is fixed. Then, for any $\mathcal I$-measurable initial state $\bm \xi_0\in L^2(\Omega;(\R^d)^N)$, \eqref{forward'} has a unique solution. Also, let $(\bm X^T,\bm Y^T,\bm Z^T)$ be the solution of \eqref{N-FBSDE} with initial state $\bm \xi_0$ at time $t_0$. Furthermore, we set $\bm Y^\infty_s=v_{\infty}(\bm X^\infty_s),\;s\ge t_0$ and the controls
    $$\alpha^i_T(s)=-D_pH^i(X^{i,T}_s,Y^{i,T}_s)\quad\text{and}\quad \alpha^i_\infty(s)=-D_pH^i(X^{i,\infty}_s,Y^{i,\infty}_s),\;\; i\in [N].$$
For every $R\in(t_0,T)$, there is a constant $C_{N,R}>0$,
independent of $T$ and depending on $R$ only through $R-t_0$, such that
\begin{equation}\label{eq:turnpike-fixed-window}
\begin{split}
&\Ex\left[
\sup_{s\in[t_0,R]}
\left\{
|\bm X_s^T-\bm X_s^\infty|^2
+|\bm Y_s^T-\bm Y_s^\infty|^2
+|\bm\alpha_T(s)-\bm\alpha_\infty(s)|^2
\right\}
\right]\\
&\hspace{9cm}\leq
C_{N,R}\left(1+\Ex[|\bm\xi_0|^2]\right)
e^{-2c(T-R)}.
\end{split}
\end{equation}
If \textup{(U)} holds, then
\begin{equation}\label{eq:turnpike-fixed-window-N}
\begin{split}
&\Ex\left[
\sup_{s\in[t_0,R]}
\left\{
|\bm X_s^{T}-\bm X_s^{\infty}|^2
+|\bm Y_s^{T}-\bm Y_s^{\infty}|^2
+|\bm\alpha_T(s)-\bm\alpha_\infty(s)|^2
\right\}
\right]\\
&\hspace{9cm}\leq
C_R^*\left(N+\Ex[|\bm\xi_0|^2]\right)
e^{-2c_*(T-R)},
\end{split}
\end{equation}
where $C_R^*$ is independent of $N$ and $T$, and depends on $R$ only through $R-t_0$.
Here $c$ and $c_*$ are the decay constants in \eqref{goal:estimate2} and \eqref{goal:estimate2-uniform}, respectively.
\end{theorem}

Theorem~\ref{thm:turnpike} is proved in Subsection~\ref{pf:turnpike}, using Step~1 of the proof of Theorem~\ref{thm:stationary-limit}.

So far, $\bm X^\infty$ and $\bm \alpha_\infty$ are obtained as limits of the equilibrium states $\bm X^T$ and strategies $\bm \alpha_T$, as $T\to +\infty$, respectively. The results below establish that $\bm \alpha_\infty$ is an open-loop equilibrium for an ergodic $N$-player game. We note that a solution of the stationary Pontryagin system of PDEs does not, by itself, specify an infinite-horizon cost criterion.

\subsubsection{Average-cost $N$-player game and its equilibria}

We introduce the long-run average running cost and state explicitly the exact strategy profiles in which the limiting control $\bm \alpha_\infty$ is an equilibrium.

Let $\bm x\in (\R^d)^N$. For $i\in [N]$ and $(\beta^i)_{s\ge t_0}$ an $(\mathcal{F}_t)_{t\ge t_0}$-progressively measurable process, we set
\begin{equation*}
X_s^{i, \beta^i}=x^i+\int_{t_0}^s\beta_u^i\,du+\sqrt2\,\sigma^i(W_s^i-W_{t_0}^i)\quad\text{and}\quad
\mathcal E^i(\beta^i)
:=\sup_{R\geq t_0+1}\frac{1}{R-t_0}\int_{t_0}^R
\Ex\bigl[|X_s^{i,\beta^i}|^2+|\beta_s^i|^2\bigr]ds\,.
\end{equation*}
We denote by $\bm X^{\bm \beta}$ the vector of $X^{i,\beta^i}$. To define the set of admissible strategies, we consider the set of $(\mathcal{F}_t)_{t\ge t_0}$-progressively measurable processes
\begin{equation*}
\mathcal{A}_\infty^i(\bm x)=\left\{\beta^i\colon[t_0,+\infty)\times\Omega\to \R^d:\;\int_{t_0}^{t}\Ex[|\beta_s^i|^2]ds<\infty,\;\forall t\ge t_0,\; \mathcal E^i(\beta^i)<\infty\right\},\;\;i\in [N]\,.
\end{equation*}
For a profile $\bm \beta \in \prod_{i=1}^N\mathcal{A}^i_\infty(\bm x)$, we write $\mathcal E(\bm\beta)=\sum_{i=1}^N\mathcal E^i(\beta^i)$ and we let
\begin{equation*}
A_T^i(\bm\beta)
=\frac{1}{T-t_0}\Ex\left[\int_{t_0}^T\left(L^i(X_s^{i,\beta^i},\beta_s^i)
+F^i(\bm X_s^{\bm\beta})\right)ds\right]\,.
\end{equation*}
Each player aims to minimize the functional
$$\overline J^i(\bm x;\bm\beta):=\limsup_{T\to\infty} A_T^i(\bm\beta),\;\;\;\;i\in[N],\text{ for }\bm\beta\in \mathcal{A}_\infty(\bm x):=\prod_{i=1}^N\mathcal{A}^i_\infty(\bm x).$$
We will refer to this game as the average-cost $N$-player game.

\begin{definition}[Strong average-cost equilibrium]\label{erg-strong}
A profile $\bm\beta\in\mathcal A_\infty(\bm x)$ is a \emph{strong
average-cost open-loop equilibrium} if, for every $i\in[N]$ and
every $\gamma\in\mathcal A^i_\infty(\bm x)$,
\begin{equation}\label{erg:strong}
\liminf_{T\to\infty}
\bigl(A_T^i([\bm\beta^{-i};\gamma])-A_T^i(\bm\beta)\bigr)\geq0.
\end{equation}
\end{definition}

\begin{theorem}[Average-cost open-loop equilibrium]\label{prop:average-Nash}
Under the hypotheses of Theorem~\ref{thm:turnpike} with $\bm \xi_0=\bm x\in (\R^d)^N$ fixed, for any $i\in [N]$, the control process constructed in Theorem \ref{thm:turnpike}, namely
$$\alpha_\infty^{i}(s):=-D_pH^i(X_s^{i,\infty},Y_s^{i,\infty}),$$belongs to $\mathcal A_\infty^i(\bm x)$.
In addition, for any $i\in [N]$ and any $\beta^i\in\mathcal A_\infty^i(\bm x)$,
$$\liminf_{T\to\infty}\left(A_T^i([\bm\alpha^{-i}_\infty;\beta^i])-A_T^i(\bm\alpha_\infty)\right)\geq0.$$
In particular, the limiting control profile $\bm\alpha_\infty$ is a strong average-cost Nash equilibrium, and
$$
\overline J^i(\bm x;\bm\alpha_\infty)
\leq\overline J^i(\bm x;[\bm\alpha^{-i}_\infty;\beta^i]).
$$
\end{theorem}
\vspace{2mm}

Theorem~\ref{prop:average-Nash} is proved in Subsection~\ref{pf:average-Nash}. The strong average-cost equilibrium is not necessarily unique (see Remark \ref{rem:erg-nonunique}); however, we show that for any $i\in [N]$, the value $\overline{J}^i(\bm x;\bm \beta)$ is the same for any equilibrium $\bm\beta\in \mathcal{A}_\infty(\bm x)$.

\begin{theorem}\label{thm:main-erg}
Suppose that (A1)-(A3) are in place. There exist constants $(\lambda^{i,N})_{i=1}^N$, independent of the initial configuration, such that, for every $\bm x\in(\R^d)^N$,
$$\lim_{T\to\infty}A_T^i(\bm\beta)=\lambda^{i,N},$$
for any $\bm\beta\in\mathcal{A}_\infty(\bm x)$ strong average-cost equilibrium in the sense of Definition \ref{erg-strong}. We refer to the constants $(\lambda^{i,N})_{i=1}^N$ as ergodic constants.
\end{theorem}

Theorem~\ref{thm:main-erg} follows from Proposition~\ref{thm:erg-value} (see Subsection~\ref{pf:main-erg}), which also gives the rate of convergence of $A_T^i(\bm\alpha_\infty)$ to $\lambda^{i,N}$.

\subsubsection{Long-time behavior of the $N$-player game} For $t_0\leq t\le T$, let $V_T^{i,N}(t,\bm x)$ be player $i$'s total cost
at the finite-horizon open-loop equilibrium started from $\bm x\in (\R^d)^N$ at
time $t$, including the terminal cost $G^i$; in particular, $V_T^{i,N}(T,\cdot)=G^i$. We also denote by $(\bm X_s^{T,\bm x})_{s\in[t_0,T]}$ the state of the $N$ players at an equilibrium, assuming that they start from $\bm x$ at time $t_0$.

We first characterize the long-time behavior of the distribution of the $N$ players in the finite-horizon game at an equilibrium. In particular, we show that $\mathcal{L}(\bm X_s^{T,\bm x})$ stays close to the invariant measure $\mu^N_\infty\in \mathcal{P}_2(\R^{Nd})$ corresponding to \eqref{forward'}.

\begin{theorem}\label{thm:mainLTMdistr}
In the setup of Theorem \ref{thm:turnpike}, there exists a measure $\mu^N_\infty\in \mathcal{P}_2(\R^{Nd})$ such that
$$W_2\!\left(\mathcal L(\bm X_s^{T,\bm x}),\mu_\infty^N\right)
\leq C(1+|\bm x|)
\left(e^{-\theta (s-t_0)}+e^{-c(T-s)}\right),
\quad s\in [t_0,T],$$
for some constants $\theta ,c,C>0$ which are independent of $T,t_0$, but might depend on $N$.

If \textup{(U)} also holds, there are constants
$C_*,c_*,\theta_*>0$, independent of $N$, $T$, and $t_0$, such that
\begin{equation}\label{eq:state-law-turnpike-uniform}
\frac{
W_2\!\left(\mathcal L(\bm X_s^{T,\bm x}),\mu_\infty^N\right)}
{\sqrt N+|\bm x|}
\leq
C_*\left(
e^{-\theta_*(s-t_0)}+e^{-c_*(T-s)}
\right),
\qquad s\in[t_0,T],
\end{equation}
where one may take
$\theta_*=\overline C_L/(2K_*)$.
\end{theorem}

Theorem~\ref{thm:mainLTMdistr} follows from Proposition~\ref{prop:erg-contraction}, which constructs $\mu^N_\infty$, and from Corollary~\ref{cor:erg-long-time-limit} (see Subsection~\ref{pf:mainLTMdistr}).

\begin{remark}
In fact, one may take $\theta= \frac{C_L}{2K}$, where $K$ is the horizon ($T$) independent Lipschitz bound of $v_T$ we will obtain in Proposition \ref{prop:uniform-Z}.
\end{remark}

Our last result is the long-time behavior of the value functions of the finite-horizon $N$-player game. We show that for any $i\in [N]$ and large $T$, $V^{i,N}_T$ is approximated by $(T-t)\lambda^{i,N}+w^{i,N}_G$, where $w^{i,N}_G\colon \R^{Nd}\to \R$ is an appropriate corrector.

\begin{theorem}\label{thm:mainLTB}
    In the setup of Theorem \ref{thm:main-erg}, for any $i\in [N]$, the limit
    $$\kappa^{i,N}_G:=\lim_{\tau\to +\infty}\left( \int V_{t_0+\tau}^{i,N}(t_0,\bm x)\,d\mu_\infty^N(\bm x)-  \tau \lambda^{i,N}\right),$$
    where $(\lambda^{i,N})_{i=1}^N$ are the ergodic constants from Theorem \ref{thm:main-erg}, exists. In addition, there is a unique classical solution $w^{i,N}_G\colon \R^{Nd}\to \R,\;\; i\in [N]$ with at most quadratic growth of
    \be\nonumber
    \begin{cases}
    \lambda^{i,N}-\sum_{j=1}^N\text{tr}\left(\sigma^j(\sigma^j)^\top D^2_{x^jx^j}w^{i,N}_G\right)+\sum_{j=1}^ND_pH^j(x^j,v^j_\infty(\bm x))\cdot D_{x^j}w^{i,N}_G\\
    \hspace{7cm}=L^i(x^i,-D_pH^i(x^i,v^i_\infty(\bm x)))+F^i(\bm x),\;\;\;\bm x\in (\R^d)^N,\\
    \int w^{i,N}_G(\bm y)d\mu^N_\infty(\bm y)=\kappa^{i,N}_G.
    \end{cases}
    \ee
    Furthermore, there exist constants $C_E,c_E>0$, which are independent of $t,T,\bm x$ and $t_0$, but might depend on $N$ such that
    $$\left|V^{i,N}_T(t,\bm x)-(T-t)\lambda^{i,N}-w^{i,N}_G(\bm x)   \right|\le C_E(1+|\bm x|^2)e^{-c_E(T-t)},\;\;\text{for all }(t,\bm x)\in [t_0,T]\times \R^{Nd}.$$
    \end{theorem}
\vspace{2mm}

Theorem~\ref{thm:mainLTB} follows from Proposition~\ref{prop:erg-poisson} and Theorem~\ref{thm:expansion-general} (see Subsection~\ref{pf:mainLTB}), with $w^{i,N}_G=u^{i,N}+\kappa^{i,N}_G$, where $u^{i,N}$ is the corrector of Proposition~\ref{prop:erg-poisson}.
\vspace{2mm}

\textbf{Convention on the initial time.}
The data $L^i$, $F^i$, $G^i$ and $\sigma^i$ do not depend on time. Hence, if $(\bm X,\bm Y,\bm Z)$ solves~\eqref{N-FBSDE} on $[t_0,T]$, then $(\bm X_{t_0+s},\bm Y_{t_0+s},\bm Z_{t_0+s})_{s\in[0,T-t_0]}$ solves the same system on $[0,T-t_0]$, with the Brownian motions $W^j_{t_0+s}-W^j_{t_0}$ and the filtration $(\mathcal F_{t_0+s})_{s\geq0}$. The same shift maps the admissible controls, the costs and the equilibria of the game on $[t_0,T]$, and of the average-cost game on $[t_0,+\infty)$, to those of the games on $[0,T-t_0]$ and on $[0,+\infty)$. The statements of this section involve the times only through the differences $T-t$, $s-t_0$ and $R-t_0$, so it suffices to prove them for $t_0=0$. From now on, in Sections~\ref{3}--\ref{7} and in the appendices, we assume that $t_0=0$.

\section{Finite-horizon system of FBSDEs: well-posedness and estimates}\label{3}

In this section, we first establish well-posedness of the finite-horizon Pontryagin
FBSDE and characterize the resulting finite-horizon open-loop Nash equilibrium. We then construct its Markovian decoupling field by establishing globally Lipschitz and classical solutions of \eqref{N-PDE} under \textup{(A1)--(A3)}.

\begin{proposition}[Well-posedness and equilibrium characterization]\label{prop:N-wellposed}
Assume \textup{(A1)--(A2)}. For every $T>0$ and every
square-integrable $\mathcal I$-measurable $\bm\xi_0$, system
\eqref{N-FBSDE} has a unique solution in the sense of
Definition~\ref{def:N-FBSDE-solution}. Moreover, a profile
$\bm\alpha\in\prod_{i=1}^N\mathcal A^i[0,T]$
is an open-loop Nash equilibrium if and only if it is induced by this solution through $\alpha_t^i=-D_pH^i(X_t^{i,N},Y_t^{i,N})$. Consequently, the open-loop equilibrium exists and is unique, up to equality $dt\otimes d\Pr$-almost everywhere.
\end{proposition}

\begin{proof}
Fix the initial condition $\bm \xi_0$ and write
\[
\mathcal H:=\prod_{i=1}^N\mathcal A^i[0,T],
\qquad
\langle\bm a,\bm b\rangle_{\mathcal H}
=\Ex\left[\int_{0}^T\bm a_t\cdot\bm b_t\,dt\right].
\]
For $\bm a\in\mathcal H$, define
$$X_t^{\bm a,i}=\xi_0^i+\int_{0}^ta_s^i ds+\sqrt{2}\sigma^iW_t^i,$$
and set
$$\begin{aligned}
Y_t^{\bm a,i}
=\Ex\left[D_{x^i}G^i(\bm X_T^{\bm a})
+\int_t^T\left(D_xL^i(X_s^{\bm a,i},a_s^i)
+D_{x^i}F^i(\bm X_s^{\bm a})\right)ds
\,\middle|\,\mathcal F_t\right],\;\; \mathcal{T}(\bm a)_t^i
=D_aL^i(X_t^{\bm a,i},a_t^i)+Y_t^{\bm a,i}.
\end{aligned}
$$
Bounded second derivatives and \textup{(A2)} imply that $\mathcal T:\mathcal H\to\mathcal H$ is globally Lipschitz, with a finite constant $K_{N,T}$. To check its monotonicity, compare $\bm a,\bm b\in\mathcal H$ and write $\Delta$ for their differences. Since $\Delta\bm X_t=\int_{0}^t\Delta\bm a_sds$, conditional expectation and Fubini's theorem give the following monotonicity

\begin{align}
\langle\mathcal T(\bm a)-\mathcal T(\bm b),\Delta&\bm a\rangle_{\mathcal H}
=\Ex\biggl[\int_{0}^T\sum_{i=1}^N
\left(\Delta D_aL_t^i\cdot\Delta a_t^i
+\Delta D_xL_t^i\cdot\Delta X_t^i\right)dt\biggr]\nonumber\\
&\;\;+\Ex\biggl[\int_{0}^T
\bigl(D\bm F(\bm X_t^{\bm a})-D\bm F(\bm X_t^{\bm b})\bigr)
\cdot\Delta\bm X_tdt\biggr]+\Ex\left[
\bigl(D\bm G(\bm X_T^{\bm a})-D\bm G(\bm X_T^{\bm b})\bigr)
\cdot\Delta\bm X_T\right]\nonumber\\
\ge &C_L\Ex\biggl[\int_{0}^T
\left(|\Delta\bm X_t|^2+|\Delta\bm a_t|^2\right)dt\biggr].\label{monotonicity}
\end{align}
Here the derivatives of $L^i$ are evaluated at $(X_t^{\bm a,i},a_t^i)$ and $(X_t^{\bm b,i},b_t^i)$, respectively.
We write $\Delta\mathcal T=\mathcal T(\bm a)-\mathcal T(\bm b)$.
For every $\varepsilon>0$, the preceding monotonicity inequality
and the Lipschitz continuity of $\mathcal T$ give
\begin{align*}
\|\Delta\bm a-\varepsilon\Delta\mathcal T\|_{\mathcal H}^2
={}&\|\Delta\bm a\|_{\mathcal H}^2
-2\varepsilon
 \langle\Delta\mathcal T,\Delta\bm a\rangle_{\mathcal H}
+\varepsilon^2\|\Delta\mathcal T\|_{\mathcal H}^2\\
\le{}&
\|\Delta\bm a\|_{\mathcal H}^2
-2\varepsilon C_L
\left(
\|\Delta\bm a\|_{\mathcal H}^2
+\Ex\biggl[\int_{0}^T|\Delta\bm X_t|^2\,dt\biggr]
\right)
+\varepsilon^2K_{N,T}^2
 \|\Delta\bm a\|_{\mathcal H}^2\\
={}&
\left(1-2\varepsilon C_L+\varepsilon^2K_{N,T}^2\right)
\|\Delta\bm a\|_{\mathcal H}^2
-2\varepsilon C_L
 \Ex\biggl[\int_{0}^T|\Delta\bm X_t|^2\,dt\biggr]\\
\le{}&
\left(1-2\varepsilon C_L+\varepsilon^2K_{N,T}^2\right)
\|\Delta\bm a\|_{\mathcal H}^2.
\end{align*}

Choosing $0<\varepsilon<\frac{2C_L}{K_{N,T}^2}$ and setting $\sqrt{1-2\varepsilon C_L+\varepsilon^2K_{N,T}^2}=:\rho_\varepsilon\in [0,1)$, this implies that the map $\mathcal T_\varepsilon:\mathcal H\longrightarrow\mathcal H,$ with
$$\mathcal T_\varepsilon(\bm a)
:=\bm a-\varepsilon\mathcal T(\bm a)$$
is a contraction in the Hilbert space $\mathcal{H}$. Thus, the Banach fixed-point theorem gives
a unique $\bm\alpha\in\mathcal H$ satisfying $\mathcal T_\varepsilon(\bm\alpha)=\bm\alpha.$ Equivalently, there is a unique $\bm \alpha\in \mathcal{H}$ such that $\mathcal{T}(\bm \alpha)=0$.
\vspace{2mm}

For $i\in [N]$, since
\[Y_t^{\bm\alpha,i}+\int_0^t D_xL^i(X_r^{\bm\alpha,i},\alpha_r^i)+D_{x^i}F^i(\bm X_r^{\bm\alpha})\,dr=\Ex\left[D_{x^i}G^i(\bm X_T^{\bm\alpha})+\int_0^T D_xL^i(X_r^{\bm\alpha,i},\alpha_r^i) +D_{x^i}F^i(\bm X_r^{\bm\alpha})\,dr\,\middle|\,\mathcal F_t\right]\]
is a martingale with respect to the filtration $(\mathcal F_t)_{t\ge0}$, martingale representation implies that there exists a progressively measurable square integrable process $\bm Z$ such that
$$dY_t^{\bm\alpha,i}
=-\left(D_xL^i(X_t^{\bm\alpha,i},\alpha_t^i)
+D_{x^i}F^i(\bm X_t^{\bm\alpha})\right)dt
+\sqrt{2}\sum_{j=1}^NZ_t^{i,j,N}dW_t^j,\;\;i\in [N].$$

We next identify the control $\bm\alpha$ using the Hamiltonian convention \eqref{Hamiltonian}. For fixed $(x,p)\in \R^d\times \R^d$, the function $u\mapsto L^i(x,u)+p\cdot u$ is strongly convex and coercive. Indeed,
\[
L^i(x,u)+p\cdot u
\ge
L^i(x,0)+\bigl(D_aL^i(x,0)+p\bigr)\cdot u
+\frac{C_L}{2}|u|^2.
\]
Consequently it has a unique minimizer, denoted by
$\widehat u^i(x,p)$, characterized by the necessary and sufficient
first-order condition
\begin{equation}\label{eq:legendre-control-foc}
D_aL^i(x,\widehat u^i(x,p))+p=0.
\end{equation}
Since $D_{aa}^2L^i\ge C_LI_d$, the implicit function theorem
shows that $\widehat u^i$ is continuously differentiable, hence differentiating
\[
H^i(x,p)
=-p\cdot\widehat u^i(x,p)
-L^i(x,\widehat u^i(x,p))
\]
and using \eqref{eq:legendre-control-foc}, we obtain
\begin{align*}
D_pH^i(x,p)
&=-\widehat u^i(x,p)
-\bigl(D_p\widehat u^i(x,p)\bigr)^\top
 \bigl(p+D_aL^i(x,\widehat u^i(x,p))\bigr)=-\widehat u^i(x,p),\\
D_xH^i(x,p)
&=-D_xL^i(x,\widehat u^i(x,p))
-\bigl(D_x\widehat u^i(x,p)\bigr)^\top
 \bigl(p+D_aL^i(x,\widehat u^i(x,p))\bigr)=-D_xL^i(x,\widehat u^i(x,p)).
\end{align*}
By the definition of $\mathcal T$, its zero satisfies
\[0=\mathcal T(\bm\alpha)_t^i=D_aL^i(X_t^{\bm\alpha,i},\alpha_t^i)+Y_t^{\bm\alpha,i},\qquad dt\otimes d\mathbb P\text{-almost everywhere}.\]
Thus $\alpha_t^i$ satisfies
\eqref{eq:legendre-control-foc} with
$x=X_t^{\bm\alpha,i}$ and $p=Y_t^{\bm\alpha,i}$.
Uniqueness of the minimizer gives
\[
\alpha_t^i
=\widehat u^i(X_t^{\bm\alpha,i},Y_t^{\bm\alpha,i})
=-D_pH^i(X_t^{\bm\alpha,i},Y_t^{\bm\alpha,i}).
\]
The second Hamiltonian identity also gives
\[
-D_xL^i(X_t^{\bm\alpha,i},\alpha_t^i)
=D_xH^i(X_t^{\bm\alpha,i},Y_t^{\bm\alpha,i}).
\]
Substitution into the equations satisfied by $X^{\bm \alpha,i}$ (forward) and $Y^{\bm \alpha,i}$ (backward) results in the system \eqref{N-FBSDE}. Square integrability (Definition \ref{def:N-FBSDE-solution}(3)) follows directly from the linear growth of $D_xL^i$, $D\bm F$, and $D\bm G$.

\vspace{2mm}

We now show that any solution $(\tilde{\bm X},\tilde{\bm Y},\tilde{\bm Z})$ of \eqref{N-FBSDE} is equal to $(\bm X,\bm Y,\bm Z)$. Indeed, uniqueness of the zero of $\mathcal{T}$ implies that $-D_pH^i(\tilde{X}_t^i,\tilde{Y}_t^i)=\alpha^i_t$, where $\bm \alpha$ is the zero of $\mathcal{T}$. This in turn determines $\tilde{\bm X}$. The conditional representation of $\tilde{\bm Y}$ then uniquely determines $\tilde{\bm Y}$, and uniqueness of martingale representation determines $\tilde{\bm Z}$ up to $dt\otimes d\mathbb P$ null sets.
\vspace{2mm}

It remains to identify the equilibrium. Fix a player $i$ and
the opponents' open-loop control processes. Their states remain
unchanged under a unilateral deviation because each player's
forward dynamics involve only that player's control and Brownian motion. Testing monotonicity of $D\bm F$ and $D\bm G$ on configurations differing only in coordinate $i$ shows that $x^i\mapsto D_{x^i}F^i(\bm x)$ and $x^i\mapsto D_{x^i}G^i(\bm x)$ are monotone when the other coordinates are fixed. Hence $F^i$ and $G^i$ are convex in $x^i$. Together with joint convexity of $L^i$ and the affine
dependence of $X^i$ on its control, this makes player $i$'s
cost a convex differentiable functional on $\mathcal{A}^i[0,T]$. For every $\beta^i\in\mathcal{A}^i[0,T]$, its derivative is
\[D_{\alpha^i}J_T^i
(0,\bm\xi_0;\bm\alpha)[\beta^i]=\Ex\left[\int_{0}^T\mathcal T(\bm\alpha)_t^i\cdot\beta_t^i\,dt\right].
\]
Suppose that $\alpha^i$ is optimal against $\bm\alpha^{-i}$. For every $\beta^i\in\mathcal{A}^i[0,T]$, the controls
$\alpha^i+r\beta^i\in \mathcal{A}^i[0,T]$, for all $r\in\R$. Therefore the differentiable function
\[r\longmapsto J_T^i(0,\bm\xi_0;\alpha^i+r\beta^i,\bm\alpha^{-i})\]
has a minimum at $r=0$, and
\[
\Ex\left[\int_{0}^T
\mathcal T(\bm\alpha)_t^i\cdot\beta_t^i\,dt\right]=0,
\quad\text{for every }\beta^i\in\mathcal{A}^i[0,T].
\]
Since $\mathcal T:\mathcal H\to\mathcal H$, we may choose $\beta^i=\mathcal T(\bm\alpha)^i$. This yields
\[0=\Ex\left[\int_{0}^T|\mathcal T(\bm\alpha)_t^i|^2\,dt\right].
\]
Consequently, $\mathcal T(\bm\alpha)_t^i=0$ $dt\otimes d\mathbb P\text{-almost everywhere}$.

Conversely, suppose $\mathcal T(\bm\alpha)^i=0$. For any competing control $\gamma^i\in\mathcal{A}^i[0,T]$, convexity gives
\begin{align*}
J_T^i(0,\bm\xi_0;\gamma^i,\bm\alpha^{-i})
-J_T^i(0,\bm\xi_0;\bm\alpha)\ge
D_{\alpha^i}J_T^i
(0,\bm\xi_0;\bm\alpha)[\gamma^i-\alpha^i]=
\Ex\biggl[\int_{0}^T
\mathcal T(\bm\alpha)_t^i\cdot
(\gamma_t^i-\alpha_t^i)\,dt\biggr]
=0.
\end{align*}
Thus $\alpha^i$ is optimal against $\bm\alpha^{-i}$.
Applying these two implications to every player proves that
the Nash equilibrium condition in
Definition~\ref{d:Nash_equilibrium_finite} is equivalent to
$\mathcal T(\bm\alpha)=0$ in $\mathcal H$.
The existence and uniqueness of that zero establish existence and uniqueness of the open-loop Nash equilibrium.
\end{proof}

Our next step is to establish the well-posedness of the Pontryagin system of PDEs \eqref{N-PDE}, as well as its connection with \eqref{N-FBSDE}. We start with the existence and we note that, at this stage, our Lipschitz and linear growth estimates depend on the horizon $T>0$ and the number of players $N$.

\begin{lemma}[Markovian decoupling field]\label{lem:N-decoupling}
Under the hypotheses of Proposition~\ref{prop:N-wellposed}, there is a deterministic decoupling field $v_T:[0,T]\times\R^{Nd}\to\R^{Nd}$ satisfying
$$\bm Y_s=v_T(s,\bm X_s)
\quad\text{and}\quad v_T(T,\bm x)=D\bm G(\bm x),$$
for every deterministic $s\in[0,T]$, almost surely. Furthermore, there exists a constant $C_{N,T}<\infty$ that might depend on $N,T$ such that
$$|v_T(s,\bm x)-v_T(s,\bar{\bm x})|
\le C_{N,T}|\bm x-\bar{\bm x}|,
\qquad
|v_T(s,\bm x)|\le C_{N,T}(1+|\bm x|),
$$
for every $s\in[0,T]$ and $\bm x,\bar{\bm x}\in\R^{Nd}$.
\end{lemma}

\begin{proof}

Fix $(s,\bm x)\in[0,T]\times\R^{Nd}$. Let
$(\bm X^{s,\bm x},\bm Y^{s,\bm x},\bm Z^{s,\bm x})$ be the solution
of \eqref{N-FBSDE} on $[s,T]$, starting from
$\bm X_s^{s,\bm x}=\bm x$, in the usual augmentation
of the future-increment filtration
\[
\mathcal G_r^s
:=
\sigma(W_u^j-W_s^j:s\leq u\leq r,\ j\in[N]),
\qquad r\in[s,T].
\]
Proposition \ref{prop:N-wellposed}, applied with initial time
$s$ and trivial initial sigma-algebra, gives existence and
uniqueness in this filtration.
Since $\mathcal G_s^s$ is trivial up to null sets,
$\bm Y_s^{s,\bm x}$ is deterministic. Moreover, this solution is adapted to the original common
filtration and satisfies the same FBSDE there.
Uniqueness in Proposition~\ref{prop:N-wellposed} identifies it
with the solution constructed in that filtration.
We therefore define the deterministic field
\[
v_T(s,\bm x):=\bm Y_s^{s,\bm x}.
\]

We next show the Lipschitz and growth estimates for $v_T$, which are uniform over $s\in[0,T]$. Denote the control operator in the preceding proof (Proposition \ref{prop:N-wellposed}) by $\mathcal T_{s,\bm x}$. For a fixed control, changing $\bm x$ to $\bar{\bm x}$ changes the forward state by the constant $\bm x-\bar{\bm x}$. Therefore,
$$
\begin{gathered}
\|\mathcal T_{s,\bm x}(\bm a)
-\mathcal T_{s,\bar{\bm x}}(\bm a)\|_{\mathcal H_s}
\le C_{N,T}\sqrt{T-s}\,|\bm x-\bar{\bm x}|,\quad \text{where }\mathcal H_s:=\prod_{i=1}^N\mathcal{A}^i[s,T].
\end{gathered}
$$
Let $\bm\alpha^{s,\bm x}$ and
$\bm\alpha^{s,\bar{\bm x}}$ denote the zeros of $\mathcal T_{s,\bm x}, \mathcal T_{s,\bar{\bm x}},$ respectively. We set
$\Delta\bm\alpha
=\bm\alpha^{s,\bm x}-\bm\alpha^{s,\bar{\bm x}}$.
Then the monotonicity \eqref{monotonicity} and the above inequality yield
\[
\begin{aligned}
C_L\|\Delta\bm\alpha\|_{\mathcal H_s}^2
&\le
\left\langle
\mathcal T_{s,\bm x}(\bm\alpha^{s,\bm x})
-\mathcal T_{s,\bm x}(\bm\alpha^{s,\bar{\bm x}}),
\Delta\bm\alpha
\right\rangle_{\mathcal H_s}
\le C_{N,T}\sqrt{T-s}\,|\bm x-\bar{\bm x}|\,
\|\Delta\bm\alpha\|_{\mathcal H_s}.
\end{aligned}
\]
Canceling the final factor, unless it is zero, gives
$$\|\bm\alpha^{s,\bm x}-\bm\alpha^{s,\bar{\bm x}}\|_{\mathcal H_s}\le C_{N,T}\sqrt{T-s}\,|\bm x-\bar{\bm x}|.$$
The forward equation, the representation of $\bm Y$ by conditional expectation, and standard martingale estimates then imply
\be\label{lipest}
\Ex\left[\sup_{s\le r\le T}
\left(|\Delta\bm X_r|^2+|\Delta\bm Y_r|^2\right)
+\int_s^T|\Delta\bm Z_r|^2dr\right]
\le C_{N,T}|\bm x-\bar{\bm x}|^2.
\ee
Likewise, since $\mathcal T_{s,\bm x}(\bm\alpha^{s,\bm x})=0$, the monotonicity \eqref{monotonicity} applied to $\bm\alpha^{s,\bm x}$ and $0$ gives $C_L\|\bm\alpha^{s,\bm x}\|_{\mathcal H_s}^2\le-\langle\mathcal T_{s,\bm x}(0),\bm\alpha^{s,\bm x}\rangle_{\mathcal H_s}$, hence
$\|\bm\alpha^{s,\bm x}\|_{\mathcal H_s}
\le C_L^{-1}\|\mathcal T_{s,\bm x}(0)\|_{\mathcal H_s}$.
 Here $\mathcal T_{s,\bm x}(0)$ is computed along the uncontrolled state $X^{0,j}_r=x^j+\sqrt2\sigma^j(W^j_r-W^j_s)$, $j\in[N]$, $r\in[s,T]$. The linear growth of $D_aL^i$, $D_xL^i$, $D\bm F$ and $D\bm G$, together with $\Ex[\sup_{s\le r\le T}|\bm X^0_r|^2]\le C(|\bm x|^2+T-s)$, gives $\|\mathcal T_{s,\bm x}(0)\|_{\mathcal H_s}\le C_{N,T}\sqrt{T-s}\,(1+|\bm x|)$. Applying the same forward and backward estimates gives
\be\label{lingrowth}
\Ex\left[\sup_{s\le r\le T}
\left(|\bm X_r^{s,\bm x}|^2+|\bm Y_r^{s,\bm x}|^2\right)
+\int_s^T|\bm Z_r^{s,\bm x}|^2dr\right]
\le C_{N,T}(1+|\bm x|^2).
\ee
The claimed spatial Lipschitz estimate for $v_T$ follows from \eqref{lipest}. The linear growth estimate follows from \eqref{lingrowth}.

We finally show that if $(\bm X,\bm Y,\bm Z)$ is the solution of \eqref{N-FBSDE} on $[0,T]$ starting from $\bm X_0=\bm \xi_0$, then we have $v_T(s,\bm X_s)=\bm Y_s$, almost surely.
\vspace{1mm}

Let $\bm\xi\in L^2(\Omega,\mathcal F_s;\R^{Nd})$, and let
$(\bm X^{s,\bm\xi},\bm Y^{s,\bm\xi},\bm Z^{s,\bm\xi})$ be the solution of
\eqref{N-FBSDE} on $[s,T]$ with $\bm X^{s,\bm\xi}_s=\bm\xi$, given by
Proposition~\ref{prop:N-wellposed} with initial sigma-algebra $\mathcal F_s$.
It suffices to show that $\bm Y^{s,\bm\xi}_s=v_T(s,\bm\xi)$ almost surely. Indeed, for the solution of \eqref{N-FBSDE} on $[0,T]$ starting from $\bm\xi_0$,
its restriction to $[s,T]$ is, by uniqueness, the solution starting from $\bm X_s$ (at time $s$);
hence $\bm Y_s=v_T(s,\bm X_s)$ almost surely.
\vspace{1mm}

If $\bm\xi=\sum_{k=1}^m\bm x_k\ind_{A_k}$, for some $\bm x_1,...,\bm x_m\in\R^{Nd}$, where $\{A_k\}\subset\mathcal F_s$
is a finite partition of $\Omega$, then
$\sum_{k=1}^m(\bm X^{s,\bm x_k},\bm Y^{s,\bm x_k},\bm Z^{s,\bm x_k})\ind_{A_k}$
solves \eqref{N-FBSDE} on $[s,T]$ starting from $\bm\xi$, because the $\ind_{A_k}$ are
$\mathcal F_s$-measurable and sum to one. Uniqueness and the definition of $v_T$ give
$\bm Y^{s,\bm\xi}_s=\sum_{k=1}^m\ind_{A_k}v_T(s,\bm x_k)=v_T(s,\bm\xi)$.

In general, let $\bm\xi^n$ be $\bm\xi$ rounded to the grid
$2^{-n}\Z^{Nd}$ on $\{|\bm\xi|\leq n\}$, and $\bm\xi^n=0$ elsewhere.
Then $\bm\xi^n$ takes finitely many values and $\bm\xi^n\to\bm\xi$ in $L^2$, as $n\to +\infty$.
The argument above, with $\bm x-\bar{\bm x}$ replaced by
$\bm\xi-\bm\xi^n$, gives
$\Ex[|\bm Y^{s,\bm\xi}_s-\bm Y^{s,\bm\xi^n}_s|^2]\leq C_{N,T}\Ex[|\bm\xi-\bm\xi^n|^2]$.
Together with the Lipschitz bound for $v_T(s,\cdot)$, this proves the claim, if we let $n\to +\infty$ in $v_T(s,\bm \xi^n)=\bm Y_s^{s,\bm \xi^n}$.
\end{proof}

In the next Proposition, we show that \eqref{N-PDE} has a unique solution. The constants appearing in the estimates might once again depend on $N,T$.

\begin{proposition}[Classical solvability of the Pontryagin system]\label{N-PDE-classical}
Assume \textup{(A1)--(A3)}. For every fixed $N$ and finite horizon $[0,T]$, the PDE system \eqref{N-PDE}
has a unique classical solution
$$v_T\in C^{1,2}([0,T)\times\R^{Nd};\R^{Nd})
\cap C([0,T]\times\R^{Nd};\R^{Nd})$$
with at most linear growth and bounded spatial derivative; that is
\[\sup_{t\in[0,T]}
\|v_T(t,\cdot)\|_*+\sup_{t\in[0,T],\,\bm x\in\R^{Nd}}
\|D_{\bm x}v_T(t,\bm x)\|_{\mathrm{op}}<\infty.\]
This solution is the decoupling field of \eqref{N-FBSDE}. In particular,
$$
Y_s^{i,N}=v_T^i(s,\bm X_s),
\qquad
Z_s^{i,j,N}=D_{x^j}v_T^i(s,\bm X_s)\sigma^j,
$$
where the first identity holds almost surely for every $s$ and the second holds $ds\otimes d\Pr$-almost everywhere.

\end{proposition}

\begin{proof}
The proof is divided into five steps. We first derive a linear Markov
representation and establish the regularity of its coefficients. We
then obtain interior classical regularity by localization, identify
the nonlinear Pontryagin system,
 identify the martingale integrand, and finally prove
uniqueness in the stated class.

\medskip
\noindent\textit{Step 1: Linear Markov representation and coefficient
regularity.}
Write
\[
g^i(\bm x):=D_{x^i}G^i(\bm x),
\qquad
\Sigma:=\operatorname{diag}(\sigma^1,\ldots,\sigma^N),
\]
and let $v_T(t,\bm x)=\bm Y_t^{t,\bm x}$ be the decoupling
field from Lemma~\ref{lem:N-decoupling}. Thus
\[
|v_T(t,\bm x)-v_T(t,\bar{\bm x})|
\leq C_{N,T}|\bm x-\bar{\bm x}|\quad\text{and}\quad
|v_T(t,\bm x)|\leq C_{N,T}(1+|\bm x|).
\]
Define
\[
\begin{aligned}
\bar b^j(t,\bm x):=-D_pH^j(x^j,v_T^j(t,\bm x))\quad\text{and}\quad
\bar c^i(t,\bm x)=D_xH^i(x^i,v_T^i(t,\bm x))
  -D_{x^i}F^i(\bm x).
\end{aligned}
\]
The estimates in \eqref{eq:H-gradient-Lipschitz}, \eqref{eq:Hx-grad} together with \textup{(A2)} and the preceding spatial Lipschitz estimate, show that
$\bar{\bm b}$ and $\bar{\bm c}$ are globally Lipschitz in $\bm x$, uniformly on $[0,T]$, and have at most linear growth. Consequently, for every $(t,\bm x)$, the equation
\[
d\bm X_r^{t,\bm x}
=\bar{\bm b}(r,\bm X_r^{t,\bm x})\,dr
+\sqrt{2}\Sigma\,d\bm W_r,
\qquad
\bm X_t^{t,\bm x}=\bm x,
\]
has a unique strong solution. It is the forward component of the
FBSDE \eqref{N-FBSDE} started from $(t,\bm x)$. The flow property and the backward
equation therefore give, for $t\leq s\leq T$,
\begin{equation}\label{eq:linear-Markov-representation}
v_T(t,\bm x)
=
\Ex\left[
v_T(s,\bm X_s^{t,\bm x})
-\int_t^s\bar{\bm c}(r,\bm X_r^{t,\bm x})\,dr
\right].
\end{equation}
We first establish the time regularity needed below. If
$|\bm x|\leq R$ and $0\leq s-t\leq1$, standard SDE estimates give
\[
\Ex\left[
\sup_{t\leq r\leq s}
|\bm X_r^{t,\bm x}-\bm x|^2
\right]
\leq C_{N,T,R}(s-t).
\]
Using \eqref{eq:linear-Markov-representation}, the spatial Lipschitz
bound for $v_T$, and the linear growth of $\bar{\bm c}$, we obtain
\[
\begin{aligned}
|v_T(t,\bm x)-v_T(s,\bm x)|
\leq
C_{N,T}\Ex[|\bm X_s^{t,\bm x}-\bm x|]
+\Ex\left[\int_t^s
|\bar{\bm c}(r,\bm X_r^{t,\bm x})|\,dr\right]\leq C_{N,T,R}|s-t|^{1/2}.
\end{aligned}
\]
Thus $v_T$ is jointly continuous, locally $1/2$-H\"older continuous
in time, and globally Lipschitz in space. The boundedness of the second derivatives of $H^i$ and $F^i$, together with the preceding time and space estimates, implies that, on every compact space--time cylinder $K$, there is a constant $C_K$
such that
\[
\begin{aligned}
|\bar{\bm b}(t,\bm x)-\bar{\bm b}(s,\bm y)|
 +|\bar{\bm c}(t,\bm x)-\bar{\bm c}(s,\bm y)|\leq
C_K\bigl(|t-s|^{1/2}+|\bm x-\bm y|\bigr),
\qquad
(t,\bm x),(s,\bm y)\in K.
\end{aligned}
\]
In particular, for any fixed $\gamma\in(0,1)$, this stronger estimate implies\footnote{We work with the exponent $\gamma\in(0,1)$ in the Schauder
estimates. The exponent $1/2$ above describes the stronger time
regularity obtained directly from the stochastic representation,
while the corresponding spatial regularity is Lipschitz continuity. Hence, for every $\gamma\in(0,1)$, we have  \(C_{\mathrm{loc}}^{1/2,1}\subset C_{\mathrm{loc}}^{\gamma/2,\gamma}.\)
We retain $\gamma$ because the standard parabolic Schauder estimates are formulated for H\"older exponents strictly between zero and one.} $\bar{\bm b},\bar{\bm c}\in C_{\mathrm{loc}}^{\gamma/2,\gamma}([0,T]\times\R^{dN}).$

\medskip
\noindent\textit{Step 2: Interior classical regularity by localization.} For the
solution started from a deterministic point, Lemma~\ref{lem:N-decoupling}
gives $\bm Y_r^{t,\bm x}=v_T(r,\bm X_r^{t,\bm x})$ almost surely, for every deterministic $r$.
Moreover,
\[
\bm Y_r^{t,\bm x}
-\int_t^r\bar{\bm c}(q,\bm X_q^{t,\bm x})\,dq
\]
is a square-integrable martingale. Optional sampling consequently
implies that, for every bounded stopping time $\rho\in[t,T]$,
\begin{equation}\label{eq:stopped-Markov-representation}
v_T(t,\bm x)
=
\Ex\left[
v_T(\rho,\bm X_\rho^{t,\bm x})
-\int_t^\rho
\bar{\bm c}(r,\bm X_r^{t,\bm x})\,dr
\right].
\end{equation}
Let $D\Subset\R^{dN}$
be a bounded smooth domain, let
$0\leq a<b<T$, and set
\[
Q=(a,b)\times D,
\qquad
\partial_p^+Q
=
\bigl(\{b\}\times\overline D\bigr)
\cup\bigl([a,b]\times\partial D\bigr).
\]
Thus $\partial_p^+Q$ is the backward parabolic boundary. For
$(t,\bm x)\in Q$, define
\[
\tau=\inf\{r\geq t:\bm X_r^{t,\bm x}\notin D\}\wedge b.
\]
Fix $i\in[N]$ and one scalar coordinate $k\in[d]$ of $v_T^i$.
Consider the scalar Cauchy--Dirichlet problem
\begin{equation}\label{eq:localized-linear-problem}
\begin{cases}
\displaystyle
\partial_tu+
\sum_{j=1}^N
\operatorname{tr}\!\left(
\sigma^j(\sigma^j)^\top D^2_{x^jx^j}u
\right)
+\sum_{j=1}^ND_{x^j}u\,\bar b^j
=\bar c^{i,k},
&\text{in }Q,\\[4pt]
u=v_T^{i,k},
&\text{on }\partial_p^+Q.
\end{cases}
\end{equation}
Here and below, equations for the vector $v_T^i$ are understood
coordinatewise. Choose smooth compatible boundary data $h_m$ on
$\partial_p^+Q$ such that $h_m\xrightarrow{m\to+ \infty} v_T^{i,k},
\;\text{uniformly on }\partial_p^+Q.$ Standard scalar Cauchy--Dirichlet theory gives a classical solution
$u_m$ of the equation in \eqref{eq:localized-linear-problem} with
boundary value $h_m$. Since the forcing term is the same for every
$m$, the backward parabolic maximum principle gives
\[
\|u_m-u_\ell\|_{C(\overline Q)}
\leq
\|h_m-h_\ell\|_{C(\partial_p^+Q)}.
\]
Hence $(u_m)_m$ converges uniformly on $\overline Q$ to a continuous
function $u$. In addition, for every $Q'\Subset Q$, the interior Schauder estimates give
\[
\|u_m\|_{C^{1+\gamma/2,\,2+\gamma}(Q')}
\leq
C_{Q'}\left(
\|u_m\|_{C(\overline Q)}
+\|\bar c^{i,k}\|_{C^{\gamma/2,\gamma}(\overline Q)}
\right).
\]
The right-hand side is bounded uniformly in $m$. Therefore, using the Arzel\`a--Ascoli theorem and passing to the limit gives $u\in C_{\mathrm{loc}}^{1+\gamma/2,\,2+\gamma}(Q)$ and $u$ satisfies the equation in
\eqref{eq:localized-linear-problem} in the classical sense.

We now identify the constructed solution $u$ with $v_T^{i,k}$ in
$Q$. Choose an increasing sequence $(D_\ell)_{\ell\geq1}$ of bounded open
sets with smooth boundaries such that
\[
D_\ell\Subset D_{\ell+1}\Subset D,
\qquad
\bigcup_{\ell\geq1}D_\ell=D.
\]
For every
sufficiently large $\ell$ such that
$\bm x\in D_\ell$ and $b-\ell^{-1}>t$, define
\[
\tau_\ell
=
\inf\{r\geq t:\bm X_r^{t,\bm x}\notin D_\ell\}
\wedge(b-\ell^{-1}).
\]
It\^o's formula applied to $u_m(r,\bm X_r^{t,\bm x})$ up to
$\tau_\ell$ gives
\[
u_m(t,\bm x)
=
\Ex\left[
u_m(\tau_\ell,\bm X_{\tau_\ell}^{t,\bm x})
-\int_t^{\tau_\ell}
\bar c^{i,k}(r,\bm X_r^{t,\bm x})\,dr
\right].
\]
Letting $\ell\to\infty$ and then using the boundary condition yields
\[
u_m(t,\bm x)
=
\Ex\left[
h_m(\tau,\bm X_\tau^{t,\bm x})
-\int_t^\tau
\bar c^{i,k}(r,\bm X_r^{t,\bm x})\,dr
\right].
\]
Since all quantities are bounded on the localized cylinder, dominated
convergence permits $m\to\infty$, and therefore
\[
u(t,\bm x)
=
\Ex\left[
v_T^{i,k}(\tau,\bm X_\tau^{t,\bm x})
-\int_t^\tau
\bar c^{i,k}(r,\bm X_r^{t,\bm x})\,dr
\right].
\]
By \eqref{eq:stopped-Markov-representation}, the right-hand side is
$v_T^{i,k}(t,\bm x)$. Thus $u=v_T^{i,k}$ in $Q$. Since $D,a,b,i$, and $k$ were arbitrary, this proves
\[
v_T\in
C_{\mathrm{loc}}^{1+\gamma/2,\,2+\gamma}
((0,T)\times\R^{dN};\R^{dN}),
\]
and hence $v_T\in C^{1,2}((0,T)\times\R^{dN};\R^{dN})$. At $t=0$, the same local argument on one-sided cylinders gives the stated regularity on $[0,T)$.

\medskip
\noindent\textit{Step 3: Identification of the nonlinear Pontryagin
system.}
The equation obtained above is
\begin{equation}\label{eq:frozen-linear-v}
\partial_tv_T^i+
\sum_{j=1}^N
\operatorname{tr}\!\left(
\sigma^j(\sigma^j)^\top
D^2_{x^jx^j}v_T^i
\right)
+\sum_{j=1}^ND_{x^j}v_T^i\,\bar b^j
=\bar c^i.
\end{equation}
Substituting the definitions of $\bar b^j$ and $\bar c^i$ into
\eqref{eq:frozen-linear-v}, and rearranging, gives
\[
\begin{aligned}
-\partial_tv_T^i
-\sum_{j=1}^N
\operatorname{tr}\!\left(
\sigma^j(\sigma^j)^\top
D^2_{x^jx^j}v_T^i
\right)
+D_xH^i(x^i,v_T^i)
+\sum_{j=1}^N
D_{x^j}v_T^iD_pH^j(x^j,v_T^j)
=D_{x^i}F^i(\bm x),
\end{aligned}
\]
which is precisely \eqref{N-PDE}. Since $v_T$ is differentiable and
has spatial Lipschitz constant $C_{N,T}$,
\[
\sup_{t\in [0,T],\,\bm x\in\R^{Nd}}
\|D_{\bm x}v_T(t,\bm x)\|_{\mathrm{op}}
\leq C_{N,T}.
\]

\medskip
\noindent\textit{Step 4: Identification of the martingale integrand.}
We now verify the relation between $v_T$ and $\bm Z$ asserted in the
proposition. Let
$(\bm X,\bm Y,\bm Z)$ be the FBSDE solution. Fix $0\leq s<r<T$ and define
\[\tau_R=\inf\{q\geq s:|\bm X_q|\geq R\}\wedge r.\]
Applying It\^o's formula to
$v_T^i(q,\bm X_q)$ on $[s,\tau_R]$ and then letting $R\to\infty$
gives, by \eqref{eq:frozen-linear-v},
\[\begin{aligned}v_T^i(r,\bm X_r)-v_T^i(s,\bm X_s)=\int_s^r\bar c^i(q,\bm X_q)\,dq+\sqrt{2}\sum_{j=1}^N\int_s^rD_{x^j}v_T^i(q,\bm X_q)\sigma^j\,dW_q^j.\end{aligned}\]
The passage to the limit follows from the boundedness of
$D_{\bm x}v_T$, the linear growth of $\bar{\bm c}$, and the
square-integrability of $\bm X$.
Since $\bm Y_q=v_T(q,\bm X_q)$, subtracting the two
semimartingale decompositions yields
\[
\sqrt{2}\sum_{j=1}^N\int_s^r
\left(
Z_q^{i,j,N}
-D_{x^j}v_T^i(q,\bm X_q)\sigma^j
\right)dW_q^j=0.
\]
It\^o's isometry therefore gives
\[
2\Ex\left[\int_s^r\sum_{j=1}^N
\left\|
Z_q^{i,j,N}
-D_{x^j}v_T^i(q,\bm X_q)\sigma^j
\right\|_{\mathrm{HS}}^2dq\right]=0,
\]
which implies $Z_q^{i,j,N}=D_{x^j}v_T^i(q,\bm X_q)\sigma^j$, $dq\otimes d\Pr\text{-almost everywhere},\;i,j\in[N]$.

\medskip
\noindent\textit{Step 5: Uniqueness.}
Let $w$ be another classical solution of
\eqref{N-PDE} so that
\[
\sup_{t\in[0,T],\,\bm x\in\R^{Nd}}
\frac{|w(t,\bm x)|}{1+|\bm x|}
+
\sup_{t\in[0,T),\,\bm x\in\R^{Nd}}
\|D_{\bm x}w(t,\bm x)\|_{\mathrm{op}}
<\infty.
\]
For deterministic $(s,\bm x)$, consider
\[
dX_r^j
=
-D_pH^j(X_r^j,w^j(r,\bm X_r))\,dr
+\sqrt{2}\sigma^j\,dW_r^j,
\qquad
\bm X_s=\bm x.
\]
The bounded second derivatives of $H^j$ and boundedness of
$D_{\bm x}w$ make this drift globally Lipschitz in $\bm x$, uniformly
in $r$, and its growth is at most linear. Hence the SDE has a unique
strong solution satisfying
\[
\Ex\bigl[\sup_{s\leq r\leq T}|\bm X_r|^2\bigr]<\infty.
\]
For $r<T$, set
\[
Y_r^i=w^i(r,\bm X_r),
\qquad
Z_r^{i,j}=D_{x^j}w^i(r,\bm X_r)\sigma^j.
\]
Fix $\varepsilon>0$ and define the stopping time
\[
\rho_R
=
\inf\{r\geq s:|\bm X_r|\geq R\}\wedge(T-\varepsilon).
\]
It\^o's formula on $[s,\rho_R]$ and the PDE for $w$ give, after
letting $R\to\infty$,
\begin{equation}\label{eq:ito-w-terminal}
\begin{aligned}
w^i(T-\varepsilon,\bm X_{T-\varepsilon})
&=w^i(s,\bm x)
+\int_s^{T-\varepsilon}
\left(
D_xH^i(X_r^i,w^i(r,\bm X_r))
-D_{x^i}F^i(\bm X_r)
\right)dr\\
&\hspace{5cm}+
\sqrt{2}\sum_{j=1}^N
\int_s^{T-\varepsilon}
D_{x^j}w^i(r,\bm X_r)\sigma^j\,dW_r^j.
\end{aligned}
\end{equation}
The passage $R\to\infty$ follows from the at-most-linear growth of
$w$ and of the $dr$-integrand in the first line of
\eqref{eq:ito-w-terminal}, the moment estimate for $\bm X$, and the boundedness of
$D_{\bm x}w$.
We now let $\varepsilon\downarrow0$. First,
\[
\begin{aligned}
\left\|
w(T-\varepsilon,\bm X_{T-\varepsilon})
-D\bm G(\bm X_T)
\right\|_{L^2}\leq
\|D_{\bm x}w\|_\infty
\|\bm X_{T-\varepsilon}-\bm X_T\|_{L^2}
+
\left\|
w(T-\varepsilon,\bm X_T)-D\bm G(\bm X_T)
\right\|_{L^2}.
\end{aligned}
\]
The first term tends to zero by the $L^2$-continuity of $\bm X$.
The second tends to zero by continuity of $w$ at $T$, its terminal
condition, and dominated convergence, since
\[
|w(T-\varepsilon,\bm X_T)|^2+|D\bm G(\bm X_T)|^2
\leq C(1+|\bm X_T|^2).
\]
Thus, \(w(T-\varepsilon,\bm X_{T-\varepsilon})
\xrightarrow{\e\to 0}D\bm G(\bm X_T),\;\;\text{in }L^2.\) If
\[c_w^i(r,\bm x)
=
D_xH^i(x^i,w^i(r,\bm x))-D_{x^i}F^i(\bm x),\;\;i\in [N],\]
then $|c_w(r,\bm x)|\leq C(1+|\bm x|)$, and hence
\[
\Ex\left[\left|
\int_{T-\varepsilon}^T
c_w(r,\bm X_r)\,dr
\right|^2\right]
\leq
\varepsilon
\Ex\left[\int_{T-\varepsilon}^T
C(1+|\bm X_r|^2)\,dr\right]
\xrightarrow{\e\to 0}0.
\]
Likewise, boundedness of $D_{\bm x}w$ and It\^o's isometry give
\[
\Ex\left[\left|
\sum_{j=1}^N\int_{T-\varepsilon}^T
D_{x^j}w^i(r,\bm X_r)\sigma^j\,dW_r^j
\right|^2\right]
\leq C\varepsilon
\xrightarrow{\e\to 0}0.
\]
Therefore \eqref{eq:ito-w-terminal} passes to $T$. The processes $\bm Y_r=w(r,\bm X_r),$ $Z_r^{i,j}=D_{x^j}w^i(r,\bm X_r)\sigma^j$ belong to the square-integrable class of
Definition~\ref{def:N-FBSDE-solution} and solve \eqref{N-FBSDE}.
Uniqueness in Proposition~\ref{prop:N-wellposed} shows that
\[w(s,\bm x)=\bm Y_s^{s,\bm x}=v_T(s,\bm x).\]
Since $(s,\bm x)$ was arbitrary, $w=v_T$. The proof is complete.
\end{proof}

\begin{remark}\label{rem:classical-constants}
Proposition~\ref{N-PDE-classical} is a fixed-$N$, fixed-horizon
regularity result. The constants in the local Schauder estimates used
in its proof may depend on $N$, on $T$ and the localized
space--time cylinders, and on the fixed diffusion matrices, including
the inverse ellipticity bound $\nu^{-1}$. These constants are distinct from the horizon-independent global spatial Lipschitz and $Z$ bounds we will establish in Proposition~\ref{prop:uniform-Z} below.
\end{remark}

\section{Turnpike estimates and limiting Pontryagin systems}\label{5}

In this section, we will show the main turnpike estimates for \eqref{N-PDE}. The following Proposition strengthens the Lipschitz estimate of Lemma \ref{lem:N-decoupling}. In particular, we show that the Lipschitz constant is independent of the horizon $T$.
\vspace{2mm}

We consider the constants
$$K_L=\max_{1\leq i\leq N}\sup_{(x,a)\in\R^{2d}}
\|D^2_{(x,a)}L^i(x,a)\|_{\op},\;\;
L_F=\operatorname{Lip}(D\bm F),\;\;
L_G=\operatorname{Lip}(D\bm G),\;\; B=K_L+L_F.$$
By assumption (A2), we know that $L_F,L_G\leq\sqrt{C_{\rm Lip}}$. We set
\be\label{lipconstant}
K:=\max\left\{
\frac{2K_L^2+B^2}{C_L},\;
2L_G+\frac{(L_G+B)^2}{C_L}
\right\},
\ee
which is independent of $T$, but might depend on the number of players $N$. When the constants $K_L,L_F,L_G,B$ are bounded independently of $N$, that is, if Assumption (U) holds, we set
\be\label{lipconstantU}
\overline B:=\overline K_L+\overline L_F,
\qquad
K_*:=
\max\left\{
\frac{2\overline K_L^2+\overline B^2}{\overline C_L},
\;
2\overline L_G+
\frac{(\overline L_G+\overline B)^2}{\overline C_L}
\right\}.
\ee

The following lemma contains the computation which allows us to compare solutions of the Pontryagin equations with the same Brownian motions.

\begin{lemma}\label{lem:pairing}
Assume \textup{(A1)--(A2)} and let $0\leq a<b<\infty$. For $k\in\{1,2\}$, let $\bm X^k,\bm Y^k$ be continuous adapted $(\R^d)^N$-valued processes and let $\bm Z^k=(Z^{k,i,j})_{i,j\in[N]}$ be a progressively measurable process on $[a,b]$ such that
\begin{equation*}
\Ex\Bigl[\sup_{a\leq s\leq b}|\bm X^k_s|^2\Bigr]+\sup_{a\leq s\leq b}\Ex\bigl[|\bm Y^k_s|^2\bigr]+\Ex\left[\int_a^b\|\bm Z^k_s\|_{\mathrm{HS}}^2\,ds\right]<\infty\,,
\end{equation*}
and which satisfy, with the same Brownian motions, for $i\in[N]$,
\begin{equation*}
\begin{cases}
dX^{k,i}_s=\bigl(-D_pH^i(X^{k,i}_s,Y^{k,i}_s)+\phi^{k,i}_s\bigr)\,ds+\sqrt2\,\sigma^i\,dW^i_s\,,\\[2pt]
dY^{k,i}_s=\bigl(D_xH^i(X^{k,i}_s,Y^{k,i}_s)-D_{x^i}F^i(\bm X^k_s)+\psi^{k,i}_s\bigr)\,ds+\sqrt2\sum_{j=1}^NZ^{k,i,j}_s\,dW^j_s\,,
\end{cases}
\end{equation*}
where $\bm\phi^k=(\phi^{k,i})_{i=1}^N$ and $\bm\psi^k=(\psi^{k,i})_{i=1}^N$ are progressively measurable and $\Ex[\int_a^b\bigl(|\bm\phi^k_s|^2+|\bm\psi^k_s|^2\bigr)ds]<\infty$. Let $\alpha^{k,i}_s=-D_pH^i(X^{k,i}_s,Y^{k,i}_s)$ and $\bm f^k_s=\bigl(D_xL^i(X^{k,i}_s,\alpha^{k,i}_s)+D_{x^i}F^i(\bm X^k_s)\bigr)_{i=1}^N$, write $\Delta$ for the difference between the quantities with $k=1$ and with $k=2$, and set
\begin{equation*}
q(s)=\Ex[\Delta\bm X_s\cdot\Delta\bm Y_s]\qquad\text{and}\qquad e(s)=\Ex\bigl[|\Delta\bm X_s|^2+|\Delta\bm\alpha_s|^2\bigr],\qquad s\in[a,b]\,.
\end{equation*}
Then the following hold.
\begin{enumerate}[label=(\roman*)]
\item The function $q$ is absolutely continuous on $[a,b]$ and, for almost every $s\in[a,b]$,
\begin{equation*}
q'(s)\leq-C_Le(s)+\Ex\bigl[\Delta\bm\phi_s\cdot\Delta\bm Y_s+\Delta\bm X_s\cdot\Delta\bm\psi_s\bigr]\,.
\end{equation*}
\item For every $s\in[a,b]$, $\Ex[|\Delta\bm Y_s|^2]\leq K_L^2e(s)$, $\Ex[|\Delta\bm f_s|^2]\leq B^2e(s)$ and $|q(s)|\leq\frac{1+K_L^2}{2}e(s)$.
\item If $\bm\phi^k=\bm\psi^k=0$ for $k\in\{1,2\}$ and $b-a\geq1$, then
\begin{equation*}
\Ex\bigl[|\Delta\bm Y_a|^2\bigr]\leq(2K_L^2+B^2)\int_a^{a+1}e(r)\,dr\,.
\end{equation*}
\end{enumerate}
\end{lemma}

\begin{proof}
By the first-order condition~\eqref{eq:legendre-control-foc} and the identities for $D_pH^i$ and $D_xH^i$ that follow it in the proof of Proposition~\ref{prop:N-wellposed},
\begin{equation}\label{eq:pairing-legendre}
Y^{k,i}_s=-D_aL^i(X^{k,i}_s,\alpha^{k,i}_s)\qquad\text{and}\qquad D_xH^i(X^{k,i}_s,Y^{k,i}_s)=-D_xL^i(X^{k,i}_s,\alpha^{k,i}_s)\,.
\end{equation}
Hence $d\bm Y^k_s=(-\bm f^k_s+\bm\psi^k_s)\,ds+\sqrt2\,\bm Z^k_s\,d\bm W_s$. The two forward equations have the same diffusion term, so $d\Delta\bm X_s=(\Delta\bm\alpha_s+\Delta\bm\phi_s)\,ds$, and It\^o's formula gives
\begin{equation*}
d(\Delta\bm X_s\cdot\Delta\bm Y_s)=\bigl(\Delta\bm\alpha_s\cdot\Delta\bm Y_s-\Delta\bm X_s\cdot\Delta\bm f_s+\Delta\bm\phi_s\cdot\Delta\bm Y_s+\Delta\bm X_s\cdot\Delta\bm\psi_s\bigr)\,ds+dM_s\,,
\end{equation*}
where $M_s=\sqrt2\int_a^s\Delta\bm X_r\cdot\Delta\bm Z_r\,d\bm W_r$. By the Cauchy--Schwarz inequality and the integrability assumption,
\begin{equation*}
\Ex\left[\left(\int_a^b|\Delta\bm X_r|^2\|\Delta\bm Z_r\|_{\mathrm{HS}}^2\,dr\right)^{1/2}\right]\leq\left(\Ex\Bigl[\sup_{a\leq r\leq b}|\Delta\bm X_r|^2\Bigr]\right)^{1/2}\left(\Ex\biggl[\int_a^b\|\Delta\bm Z_r\|_{\mathrm{HS}}^2\,dr\biggr]\right)^{1/2}<\infty\,,
\end{equation*}
so $M$ is a martingale by the Burkholder--Davis--Gundy inequality. The $ds$-integrand is integrable on $[a,b]\times\Omega$, by the linear growth of $D_pH^i$, $D_xL^i$ and $D\bm F$ and the integrability assumptions. Taking expectations shows that $q$ is absolutely continuous, and that, for almost every $s$, $q'(s)$ is the expectation of the $ds$-integrand. By~\eqref{eq:pairing-legendre} and~\eqref{Ldispl},
\begin{equation*}
\Delta\alpha^i_s\cdot\Delta Y^i_s-\Delta X^i_s\cdot\Delta D_xL^i_s=-\bigl(\Delta D_aL^i_s\cdot\Delta\alpha^i_s+\Delta D_xL^i_s\cdot\Delta X^i_s\bigr)\leq-C_L\bigl(|\Delta X^i_s|^2+|\Delta\alpha^i_s|^2\bigr)\,,
\end{equation*}
where $\Delta D_aL^i_s$ and $\Delta D_xL^i_s$ are the differences of $D_aL^i$ and $D_xL^i$ evaluated at $(X^{1,i}_s,\alpha^{1,i}_s)$ and at $(X^{2,i}_s,\alpha^{2,i}_s)$. Moreover, $\Delta\bm X_s\cdot\bigl(D\bm F(\bm X^1_s)-D\bm F(\bm X^2_s)\bigr)\geq0$ by~\eqref{Fmon}. Summing over $i$ proves~(i).

By~\eqref{eq:pairing-legendre} and the bound $K_L$ on the Hessians of $L^i$, we have $|\Delta\bm Y_s|\leq K_L\bigl(|\Delta\bm X_s|^2+|\Delta\bm\alpha_s|^2\bigr)^{1/2}$, and the Lipschitz bound for $D\bm F$ gives $|\Delta\bm f_s|\leq(K_L+L_F)\bigl(|\Delta\bm X_s|^2+|\Delta\bm\alpha_s|^2\bigr)^{1/2}$. Taking expectations, and using $|q(s)|\leq\frac12\Ex\bigl[|\Delta\bm X_s|^2+|\Delta\bm Y_s|^2\bigr]$, proves~(ii).

For~(iii), let $u\in[a,a+1]$. Since $\int\Delta\bm Z\,d\bm W$ is a martingale, the backward equations give $\Delta\bm Y_a=\Ex\bigl[\Delta\bm Y_u+\int_a^u\Delta\bm f_r\,dr\,\big|\,\mathcal F_a\bigr]$. The conditional Jensen inequality, the Cauchy--Schwarz inequality and~(ii) give
\begin{equation*}
\Ex\bigl[|\Delta\bm Y_a|^2\bigr]\leq2K_L^2e(u)+2B^2(u-a)\int_a^ue(r)\,dr\,.
\end{equation*}
We integrate over $u\in[a,a+1]$. By Fubini's theorem and $\int_r^{a+1}(u-a)\,du\leq\frac12$ for $r\in[a,a+1]$, this proves~(iii).
\end{proof}

Lemma~\ref{lem:pairing} is used in two ways. Without perturbations, it compares two equilibria of finite-horizon games, or two solutions of the stationary system, and the applications differ in the information available at the endpoints of the time interval: the monotonicity of $D\bm G$ at the terminal time (Proposition~\ref{prop:uniform-Z} and Lemma~\ref{lem:W-lip-general}), equal initial states together with uniform second moments (Theorems~\ref{thm:horizon-comparison} and~\ref{thm:stationary-limit}, and Proposition~\ref{prop:erg-stable-Pontryagin}), or bounds on $q$ given by the monotonicity of $v_\infty$ (Lemma~\ref{lem:erg-strong-monotonicity} and Proposition~\ref{prop:erg-contraction}). With perturbations, it compares the equilibrium of the finite-horizon game, or an arbitrary admissible profile, with the stationary dynamics; the drift error then enters through the last term in~(i) (Lemma~\ref{lem:erg-pointwise-turnpike} and Theorem~\ref{thm:erg-Pontryagin-residual}).

\begin{proposition}\label{prop:uniform-Z}
Assume \textup{(A1)--(A3)}. For any $s\in [0, T]$ and any $\bm x\in\R^{Nd}$, let $(\bm X_t^{\bm x},\bm Y_t^{\bm x},\bm Z_t^{\bm x})_{t\ge s}$ be the solution to \eqref{N-FBSDE} starting from $\bm X_s=\bm x$.  Then, the decoupling field $v_T(s,\bm x)=\bm Y_s^{s,\bm x}$ satisfies
\be\label{lipbound}
|v_T(s,\bm x)-v_T(s,\bar{\bm x})|
\leq K|\bm x-\bar{\bm x}|,\;\;\text{for any }s\in [0,T],\;\bm x,\bar{\bm x}\in \R^{Nd},
\ee
where $K$ is the constant from \eqref{lipconstant}. Furthermore, for any $\bm x\in\R^{Nd}$ the matrices $\bm Z^{\bm x}$ satisfy
\be\label{ZL2}
\sum_{j=1}^N\|Z_t^{i,j,N}\|_{\mathrm{HS}}^2\leq dK^2S^2,\;\;\text{for any }i\in [N],\quad dt\otimes d\Pr\text{-almost everywhere on }[0,T),
\ee
where $S=\max_{j\in [N]}\|\sigma^j\|_{\op}$ and $K$ is the constant from \eqref{lipconstant}. Under \textup{(U)}, $v_T$ satisfies
\be\label{lipboundU}
|v_T(s,\bm x)-v_T(s,\bar{\bm x})|_N
\leq K_*|\bm x-\bar{\bm x}|_N,
\ee
where $K_*$ is the constant from \eqref{lipconstantU} and is independent of $N$, $s$, and $T$. Moreover,
\[
\|\bm Z_t^N\|_{\mathrm{HS},N}^2
\leq
K_*^2\frac1N\sum_{j=1}^N
\|\sigma^{j,N}\|_{\mathrm{HS}}^2
\leq K_*^2B_0, \text{ for }dt\otimes d\Pr\text{-almost every }t<T.
\]
\end{proposition}

\begin{proof}
Fix $s<T$ and couple the solutions starting from $\bm x$ and
$\bar{\bm x}$ (at time $t=s$) using the same Brownian motions. Write $\Delta$ for their
difference, suppress the superscript $N$ on the states, and
use the notation of Lemma~\ref{lem:pairing} on the interval $[s,T]$, with $\bm\phi^k=\bm\psi^k=0$; the integrability conditions of the lemma hold by Definition~\ref{def:N-FBSDE-solution}. In particular, $e(t)=\Ex\bigl[|\Delta\bm X_t|^2+|\Delta\bm\alpha_t|^2\bigr]$ and $q(t)=\Ex[\Delta\bm X_t\cdot\Delta\bm Y_t]$. By Lemma~\ref{lem:pairing}(i), $q'\leq-C_Le$ almost everywhere on $[s,T]$.
Since the terminal pairing is nonnegative by the monotonicity of
$D\bm G$, integration yields
\be\label{eq:vT-pairing}
C_L\int_s^T e(t)\,dt
\leq(\bm x-\bar{\bm x})\cdot
\bigl(v_T(s,\bm x)-v_T(s,\bar{\bm x})\bigr).
\ee
Let $r:=|\bm x-\bar{\bm x}|$ and
$u:=|v_T(s,\bm x)-v_T(s,\bar{\bm x})|$. We separate the proof of the Lipschitz bound \eqref{lipbound} into two cases.

If $T-s\geq1$,
 Lemma~\ref{lem:pairing}(iii) and~\eqref{eq:vT-pairing} give
$$u^2\leq(2K_L^2+B^2)\int_s^{s+1}e(t)\,dt\leq\frac{2K_L^2+B^2}{C_L}\,ru.$$
This proves \eqref{lipbound} in this case, including $u=0$.

If $h=T-s\leq1$, the terminal condition, the backward equation, the bound for $\Delta\bm f$ in Lemma~\ref{lem:pairing}(ii), and
$\Delta\bm X_T=\bm x-\bar{\bm x}+\int_s^T\Delta\bm\alpha_t\,dt$
give, by Cauchy--Schwarz,
$$
u\leq L_Gr+(L_G+B)\left(h\int_s^T e(t)\,dt\right)^{1/2}
\leq L_Gr+\frac{L_G+B}{\sqrt{C_L}}\sqrt{ru}.
$$
Young's inequality then gives
$$u\leq\left(2L_G+\frac{(L_G+B)^2}{C_L}\right)r.$$
This implies \eqref{lipbound} in this case as well.
\vspace{2mm}

We finish the proof by showing the $L^2$-bound \eqref{ZL2}. For notational convenience we write $(\bm X_t,\bm Y_t,\bm Z_t)$ for $(\bm X^{\bm x}_t,\bm Y^{\bm x}_t,\bm Z^{\bm x}_t)$. Also, let $\Sigma=\text{diag}(\sigma^1,\ldots,\sigma^N)$. By It\^o's formula, it is easy to derive that $\bm Z_t=D_{\bm x}v_T(t,\bm X_t)\Sigma$. Thus, using \eqref{lipbound}, we get
$$\sum_{j=1}^N\|Z_t^{i,j,N}\|_{\mathrm{HS}}^2
=\|D_{\bm x}v_T^i(t,\bm X_t)\Sigma\|_{\mathrm{HS}}^2
\leq S^2\|D_{\bm x}v_T^i(t,\bm X_t)\|_{\mathrm{HS}}^2
\leq dK^2S^2.$$

Under \textup{(U)}, the same estimates hold with $K_*$ in place of $K$ and \eqref{lipboundU} follows. The Lipschitz bound and \textup{(U)} give
\[
\|\bm Z_t\|_{\mathrm{HS},N}^2
=\frac1N\|D_{\bm x}v_T(t,\bm X_t)\Sigma\|_{\mathrm{HS}}^2
\leq K_*^2\frac1N\sum_{j=1}^N\|\sigma^{j}\|_{\mathrm{HS}}^2
\leq K_*^2B_0.\qedhere
\]
\end{proof}
\vspace{2mm}

We note that Proposition \ref{prop:uniform-Z} is an improvement of the stability estimate \eqref{lipest} that gave us the Lipschitz continuity (with constant depending on the terminal time $T$) in Lemma \ref{lem:N-decoupling}. In the following Proposition, we establish an analogous improvement of the $L^2$-bound \eqref{lingrowth} and hence obtain a uniform $L^2$-bound of the solution of the FBSDE system \eqref{N-FBSDE}.

\begin{proposition}[Uniform second moments]\label{prop:uniform-moments}
Assume \textup{(A1)--(A3)} and fix $N$ and a terminal time $T>0$. There exists a constant $C_N<\infty$, independent of $T$ and of $\bm\xi_0$, depending only on the structural constants, the values of the first derivatives at the origin, and the diffusion matrices, such that every solution of \eqref{N-FBSDE} with initial condition $\bm\xi_0$ satisfies
\be\label{unifL2}
\sup_{0\leq t\leq T}\Ex\left[|\bm X_t|^2+|\bm Y_t|^2+|\bm\alpha_t|^2\right]\leq C_N\left(1+\Ex[|\bm\xi_0|^2]\right),
\ee
where $\alpha_t^i=-D_pH^i(X_t^i,Y_t^i)$, $i\in [N]$. In particular, if $\bm\xi_0=\bm x\in (\R^d)^N$, the solution $v_T$ of \eqref{N-PDE} satisfies
\be\label{Tuniflingrowth}
\|v_T(t,\cdot)\|_*\le \sqrt{C_N},\qquad t\in[0,T]\,.
\ee
If Assumption \textup{(U)} is also in place, then there
exists a constant $C_*>0$ independent of $T,N$ such that
\be\label{eq:uniform-moments-N}
\sup_{0\leq t\leq T}
\Ex\left[
 |\bm X_t|^2
 +|\bm Y_t|^2
 +|\bm\alpha_t|^2
\right]
\leq
C_*\left(N+\Ex[|\bm\xi_0|^2]\right).
\ee
In particular, if $v_T$ is the solution of \eqref{N-PDE} established in Proposition \ref{N-PDE-classical}, then
\begin{equation}\label{eq:uniform-v-weighted}
\sup_{N\geq1}
\sup_{T>0}
\sup_{0\leq t<T}
\|v_T(t,\cdot)\|_{*,N}
\leq
\sqrt{C_*}.
\end{equation}
\end{proposition}

\begin{proof}
  We first assume $T\geq1$; the case $0<T\leq1$ is treated at the end of the proof, in both settings. We write $M_0:=1+\Ex[|\bm\xi_0|^2]$. Throughout the proof, constants will depend on the stated structural data, but not on the terminal time $T$.
\vspace{2mm}

Let $\Sigma=\operatorname{diag}(\sigma^1,\ldots,\sigma^N)$ and let $\bm R$ be the Ornstein--Uhlenbeck process driven by the same Brownian motions as $\bm X$:
$$
d\bm R_t=-\bm R_t\,dt+\sqrt{2}\Sigma\,d\bm W_t,\;\;t\ge 0,
\qquad \bm R_0=\bm\xi_0.
$$
Its explicit representation yields the uniform $L^2$-bound
$$
\sup_{t\geq0}\Ex[|\bm R_t|^2]
\leq\Ex[|\bm\xi_0|^2]+\|\Sigma\|_{\mathrm{HS}}^2
\leq C M_0.
$$
We set $\bm U_t=\bm X_t-\bm R_t$ and $\bm A_t=\bm\alpha_t+\bm R_t$. Then $\bm U_0=0$ and $d\bm U_t=\bm A_t\,dt.$ In particular, $\bm U$ has finite variation. We introduce
$$e(t):=\Ex\left[|\bm U_t|^2+|\bm A_t|^2\right]\quad\text{and}\quad P(t):=\Ex[\bm U_t\cdot\bm Y_t],$$
and we use the vector notation
$$
\bm\ell_a(\bm x,\bm a)
=\bigl(D_aL^i(x^i,a^i)\bigr)_{i=1}^N\quad\text{and}\quad
\bm\ell_x(\bm x,\bm a)
=\bigl(D_xL^i(x^i,a^i)\bigr)_{i=1}^N.
$$
By duality (recall the proof of Proposition \ref{prop:N-wellposed}) we have $\bm Y_t=-\bm\ell_a(\bm X_t,\bm\alpha_t)$ and
$$d\bm Y_t=-\bm f_t\,dt+\sqrt{2}\bm Z_t\,d\bm W_t,\qquad\text{where }\;\;\bm f_t=\bm\ell_x(\bm X_t,\bm\alpha_t)+D\bm F(\bm X_t).$$
The linear growth of $\bm\ell_a$, $\bm\ell_x$, and $D\bm F$ shows that
$$
\sup_{t\geq0}\Ex\left[
|\bm\ell_a(\bm R_t,-\bm R_t)|^2+
|\bm\ell_x(\bm R_t,-\bm R_t)+D\bm F(\bm R_t)|^2
\right]\leq C M_0.
$$
The function $P$ is absolutely continuous, so by the product rule and using \textup{(A1)} we have
\begin{align*}
P'(t)
={}&-\Ex\left[
\bm A_t\cdot\bigl(\bm\ell_a(\bm X_t,\bm\alpha_t)
-\bm\ell_a(\bm R_t,-\bm R_t)\bigr)\right]-\Ex\left[
\bm U_t\cdot\bigl(\bm\ell_x(\bm X_t,\bm\alpha_t)
-\bm\ell_x(\bm R_t,-\bm R_t)\bigr)\right]\\
&-\Ex\left[\bm U_t\cdot\bigl(D\bm F(\bm X_t)-D\bm F(\bm R_t)\bigr)\right]-\Ex\left[\bm A_t\cdot\bm\ell_a(\bm R_t,-\bm R_t)
+\bm U_t\cdot\bigl(\bm\ell_x(\bm R_t,-\bm R_t)+D\bm F(\bm R_t)\bigr)\right]\\
\leq{}&-\frac{C_L}{2}e(t)+C_0 M_0.
\end{align*}
Taking expectation in It\^o's formula and using the same Lipschitz and growth bounds, we get that there are constants $C_P\geq1$ and $C_1>0$ such that
$$
|P(t)|\leq C_Pe(t)+C_1M_0\quad\text{and}\quad
\Ex\left[|\bm Y_t|^2+|\bm f_t|^2\right]\leq C\bigl(e(t)+M_0\bigr).$$
Since $\bm U_T=\bm U_q+\int_q^T\bm A_r\,dr$, integration over $q\in[T-1,T]$ yields
$$\Ex[|\bm U_T|^2]\leq2\int_{T-1}^T e(r)\,dr.$$
Thus, at the terminal time, monotonicity of $D\bm G$ gives
\begin{align}
P(T)&=\Ex\left[\bm U_T\cdot\bigl(D\bm G(\bm X_T)-D\bm G(\bm R_T)\bigr)+\bm U_T\cdot D\bm G(\bm R_T)\right]\nonumber\\
&\geq-\varepsilon\Ex[|\bm U_T|^2]-C_\varepsilon M_0\ge -2\e\int_{T-1}^Te(t)dt-C_\e M_0,\label{termforpT}
\end{align}
for some constant $C_\varepsilon>0$, for every $\varepsilon>0$.
Fix any $t\in[0,T]$ and put $\eta_t(r)=e^{-\beta|r-t|}$, where
$$\beta=\min\left\{1,\frac{C_L}{8C_P}\right\}>0\quad\text{and}\quad \varepsilon=\frac{C_L}{16e^\beta}.$$
Integration by parts, using $P(0)=0$, $|\eta_t'|\leq\beta\eta_t$, and $\eta_t(T)\leq e^\beta\eta_t(r)$ for $r\in[T-1,T]$, gives
\begin{align*}
\frac{C_L}{2}\int_0^T\eta_t(r)e(r)\,dr
&\leq-\eta_t(T)P(T)+\int_0^T\eta_t'(r)P(r)\,dr
+C_0M_0\int_0^T\eta_t(r)\,dr\\
&\leq\left(2\varepsilon e^\beta+\beta C_P\right)
\int_0^T\eta_t(r)e(r)\,dr+C M_0.
\end{align*}
Here we used $\eta_t(T)\leq1$ and $\int_0^T\eta_t(r)\,dr\leq2/\beta$ to bound the remaining terms. By the choices of $\beta$ and $\varepsilon$, the coefficient in parentheses is at most $C_L/4$. Consequently,
$$\int_0^T \eta_t(r)e(r)\,dr\leq C M_0,\qquad t\in[0,T].$$
It remains to recover a pointwise estimate for $e(t)$ from this integral estimate. Let
\[I_t=\begin{cases}
[t,t+1],&t\leq T-1,\\
[T-1,T],&t>T-1.
\end{cases}\]
Then $I_t$ has length one and contains $t$. For every $q\in I_t$, we have $\bm U_t=\bm U_q+\int_q^t\bm A_r\,dr,$ hence
\[
\Ex[|\bm U_t|^2]
\leq
2\Ex[|\bm U_q|^2]
+2|t-q|
\int_{q\wedge t}^{q\vee t}\Ex[|\bm A_r|^2]\,dr.
\]
Integrating with respect to $q\in I_t$ and using
$|t-q|\leq1$ gives
\[
\begin{aligned}
\Ex[|\bm U_t|^2]
\leq
2\int_{I_t}\Ex[|\bm U_q|^2]\,dq
+2\int_{I_t}\Ex[|\bm A_r|^2]\,dr=2\int_{I_t}e(r)\,dr\leq
2e^\beta\int_0^Te^{-\beta|r-t|}e(r)\,dr
\leq CM_0.
\end{aligned}
\]
Together with the uniform second-moment bound for $\bm R$, this yields
\[\sup_{0\leq t\leq T}\Ex[|\bm X_t|^2]\leq CM_0.\]

We next estimate $\bm Y$. If $t\leq T-1$, the backward equation and conditional Jensen's inequality give, for $q\in[t,t+1]$,
$$
\Ex[|\bm Y_t|^2]
\leq2\Ex[|\bm Y_q|^2]+2(q-t)\int_t^q\Ex[|\bm f_r|^2]\,dr.
$$
Integrating over $q\in I_t$ and using the preceding weighted estimate yields $\Ex[|\bm Y_t|^2]\leq C M_0$. If $t\in[T-1,T]$, the terminal representation instead gives
$$
\Ex[|\bm Y_t|^2]
\leq2\Ex[|D\bm G(\bm X_T)|^2]
+2(T-t)\int_t^T\Ex[|\bm f_r|^2]\,dr
\leq C M_0.
$$
For the fixed-$N$ assertion, the identity
$\bm Y_t=-\bm\ell_a(\bm X_t,\bm\alpha_t)$, strong convexity, and the
Hessian bound give $C_L|\bm\alpha_t|\leq|\bm Y_t|+K_L|\bm X_t|+|\bm\ell_a(\bm0,\bm0)|.$ Therefore,
\[\Ex[|\bm\alpha_t|^2]\leq C_N\left(1+\Ex[|\bm X_t|^2]+\Ex[|\bm Y_t|^2]\right),\]
which proves the control part of the fixed-$N$ moment estimate. This finishes the proof of \eqref{unifL2}.
\vspace{2mm}

For the uniform in $N$ bound \eqref{eq:uniform-moments-N} under Assumption (U), we set $M_{0,N}:=N+\Ex[|\bm\xi_0|^2]$. Assumption \textup{(U)} gives
\[\sup_{t\geq0}\Ex[|\bm R_t|^2]\leq
\Ex[|\bm\xi_0|^2]+\sum_{i=1}^N\|\sigma^{i}\|_{\mathrm{HS}}^2\leq C M_{0,N},\]
where $C$ is independent of $N$. The origin and Lipschitz bounds in
\textup{(U)} also give
\[
\begin{aligned}
\sup_{t\geq0}\Ex\Big[
|\bm\ell_a(\bm R_t,-\bm R_t)|^2+
|\bm\ell_x(\bm R_t,-\bm R_t)
+D\bm F(\bm R_t)|^2
\Big]
\leq C M_{0,N}\;\;\text{and}\;\;\Ex[|D\bm G(\bm R_T)|^2]\leq C M_{0,N}.
\end{aligned}
\]
Consequently, for $T\geq1$, the preceding proof applies with
$M_{0,N}$ in place of $M_0$, and every constant is independent of
$N$.

For $0<T\leq1$, in both settings (with $M_0$ in the fixed-$N$ case and $M_{0,N}$ under \textup{(U)}), one has
\[
\bm U_T=\int_0^T\bm A_r\,dr,\;\;\text{so}\;\;\;
\Ex[|\bm U_T|^2]
\leq\int_0^T e(r)\,dr.
\]
Using this inequality in the terminal estimate for $P(T)$ \eqref{termforpT} and
repeating the preceding argument gives
\[
\int_0^T e^{-\beta|r-t|}e(r)\,dr
\leq C M_{0,N}\quad(\text{resp. }CM_0).
\]
With a similar argument as before (now with $I_t:=[0,T]$), this integral bound gives the bound for $\Ex[|\bm U_t|^2]$ and the terminal representation of $\bm Y_t$ gives the corresponding bound for $\Ex[|\bm Y_t|^2]$. Finally, the identity
$\bm Y_t=-\bm\ell_a(\bm X_t,\bm\alpha_t)$ and uniform strong
convexity give
\[C_L|\bm\alpha_t|\leq|\bm Y_t|+K_L|\bm X_t|+|\bm\ell_a(\bm0,\bm0)|.\]
Therefore,
\(|\bm\alpha_t|^2\leq C\left(N+|\bm X_t|^2+|\bm Y_t|^2\right),\) where $C$ is independent of $N$. This proves
\eqref{eq:uniform-moments-N}; in the fixed-$N$ case, the bound on $\bm\alpha_t$ obtained above for $T\geq1$ applies verbatim and completes the proof of \eqref{unifL2}.
\vspace{2mm}

For \eqref{eq:uniform-v-weighted}, we rely on \eqref{eq:uniform-moments-N} with the
initial time equal to $t$ and with deterministic initial state $\bm X_t=\bm x$. Since $\bm Y_t=v_T(t,\bm x)$, the $L^2$-bound gives
\[
|v_T(t,\bm x)|_N^2
\leq C_*\left(1+|\bm x|_N^2\right).
\]
Dividing by $(1+|\bm x|_N)^2$ and taking the supremum over $t\in [0,T]$ proves
\eqref{eq:uniform-v-weighted}. This completes the proof.
\end{proof}

\subsection{Proofs of the turnpike estimates} We are now ready to prove the turnpike estimates.

\subsubsection*{Proof of Theorem~\ref{thm:horizon-comparison}}\phantomsection\label{pf:horizon-comparison}
Fix $(t,\bm x)\in [0,T]\times(\R^d)^N$ and consider the solutions of \eqref{N-FBSDE} with terminal data
$(T,D\bm G)$ and $(T',D\widetilde{\bm G})$ driven by the same Brownian
motions, starting both at $\bm x$ at time $t$. On their common
interval $[t,T]$, denote the differences by
$\Delta\bm X$, $\Delta\bm Y$, and $\Delta\bm\alpha$, where
$\alpha_r^i=-D_pH^i(X_r^i,Y_r^i)$ for each solution. We set
$$e(r)=\Ex\left[|\Delta\bm X_r|^2+|\Delta\bm\alpha_r|^2\right],\;\;\Phi(r)=\Ex\left[|\Delta\bm X_r|^2+|\Delta\bm Y_r|^2\right],\;\;
q(r)=\Ex[\Delta\bm X_r\cdot\Delta\bm Y_r].$$
These are the quantities of Lemma~\ref{lem:pairing} on $[t,T]$, with $\bm\phi^k=\bm\psi^k=0$; the integrability conditions of the lemma hold by Definition~\ref{def:N-FBSDE-solution}. By Lemma~\ref{lem:pairing}(i)--(ii), $q'\leq-C_Le$ almost everywhere on $[t,T]$ and $|q|\leq\frac{1+K_L^2}{2}e$. Since $\Delta\bm X_t=0$, we have $q(t)=0$, hence $q\leq0$ and, for almost every $r\in[t,T]$,
\begin{equation}\label{ineq1}
-q'(r)\geq C_Le(r)\geq-2c\,q(r)\,,\qquad c=\frac{C_L}{1+K_L^2}\,.
\end{equation}
Integration of this differential inequality yields
$$-q(r)\leq -e^{-2c(T-r)}q(T),\qquad t\leq r\leq T.$$
Proposition~\ref{prop:uniform-moments}, applied separately to the
original solutions with terminal data $(T,D\bm G)$ and
$(T',D\widetilde{\bm G})$, gives
$$-q(T)\leq\frac12\Phi(T)\leq2M(1+|\bm x|^2),$$
for some constant $M>0$, where the last inequality follows from \eqref{unifL2}.
Suppose first that $T-t\geq1$.
 Lemma~\ref{lem:pairing}(iii), and the integral of~\eqref{ineq1} over $[t,t+1]$, give
$$|\Delta\bm Y_t|^2
\leq(2K_L^2+B^2)\int_t^{t+1}e(r)\,dr
\leq-\frac{2K_L^2+B^2}{C_L}q(t+1),$$
because $q(t)=0$. Combining this with the preceding bounds proves
$$
|\Delta\bm Y_t|^2
\leq\frac{2M(2K_L^2+B^2)}{C_L}
e^{-2c(T-t-1)}(1+|\bm x|^2).
$$
If $0\leq T-t\leq1$, the same uniform moment estimate at the initial
time instead gives
$$
|\Delta\bm Y_t|^2
\leq4M(1+|\bm x|^2)
\leq4Me^{2c}e^{-2c(T-t)}(1+|\bm x|^2).
$$
Taking square roots in these two cases and using
$\sqrt{1+|\bm x|^2}\leq1+|\bm x|$ proves \eqref{goal:estimate1}, since $\Delta\bm Y_t=v_T^{G}(t,\bm x)-v_{T'}^{\widetilde{G}}(t,\bm x)$.
\vspace{2mm}

Suppose that \textup{(U)} holds for both sets of terminal costs,
with common constants independent of $N$. Let $M_*$ be the maximum
of the two constants furnished by
\eqref{eq:uniform-moments-N}, and set
$c_*:=\overline C_L/(1+\overline K_L^2)$.
Recall $\overline B=\overline K_L+\overline L_F$.
The same pairing argument, with the uniform structural bounds,
gives
\[
-q(r)\leq-e^{-2c_*(T-r)}q(T),
\qquad
-q(T)\leq2M_*(N+|\bm x|^2).
\]
Consequently,
\[
|\Delta\bm Y_t|^2
\leq
\begin{cases}
\displaystyle
\frac{2M_*(2\overline K_L^2+\overline B^2)}
     {\overline C_L}
e^{-2c_*(T-t-1)}(N+|\bm x|^2),
& T-t\geq1,\\[1.2ex]
\displaystyle
4M_*e^{2c_*}e^{-2c_*(T-t)}(N+|\bm x|^2),
& 0\leq T-t\leq1.
\end{cases}
\]
Taking square roots and using
$\sqrt{N+|\bm x|^2}\leq\sqrt N+|\bm x|$ proves
\eqref{goal:unifNest} with a constant $C_*$ independent of $N$.

\hfill\qed

\begin{remark}
When $T'>T$, the longer-horizon solution restricted to $[t,T]$ has
terminal value $\bm Y_T=v_{T'}^{\widetilde G}(T,\bm X_T)$.
The proof uses this solution only as the restriction of the original
system with horizon $T'$. The endpoint pairing is controlled
directly by Proposition~\ref{prop:uniform-moments} for the two
original solutions. Taking $T'=T$ compares two terminal costs on
the same horizon, while taking $\widetilde{\bm G}=\bm G$ compares
two horizons with the same terminal cost.
\end{remark}

We now turn to the proof of Theorem \ref{thm:stationary-limit}.

\subsubsection*{Proof of Theorem~\ref{thm:stationary-limit}}\phantomsection\label{pf:stationary-limit}
Translating the system \eqref{N-PDE} on $[t,T]$ by $r\mapsto r-t$ and invoking the uniqueness yields
\[v_T(t,\bm x)=v_{T-t}(0,\bm x)=:u(T-t,\bm x),\]
where $u\colon [0,+\infty)\times (\R^d)^N\to (\R^d)^N$. In particular, $u(h,\bm x):=v_h(0,\bm x)$, for any $h\ge 0$ and $\bm x\in (\R^d)^N$.
By Theorem~\ref{thm:horizon-comparison}, $u(h,\cdot)$ is Cauchy as
$h\to\infty$, in the weighted norm $\|\cdot\|_*$. Let $v_\infty\colon (\R^d)^N\to (\R^d)^N$ be its limit. Sending $T'\to +\infty$ in \eqref{goal:estimate1} gives \eqref{goal:estimate2}. The proof of \eqref{goal:estimate2-uniform} follows from the same argument if we send $T'\to +\infty$ in \eqref{goal:unifNest}. Setting $T=T'\to +\infty$ in \eqref{goal:estimate1}, also shows that $v_\infty$ is independent of the choice of terminal cost $G$.
\vspace{2mm}

We now show that $v_\infty$ is globally Lipschitz and has at most linear growth. Indeed, the Lipschitz bound and the linear growth follow from Proposition \ref{prop:uniform-Z} (relation \eqref{lipbound}) and Proposition \ref{prop:uniform-moments} (relation \eqref{Tuniflingrowth}) by sending $T\to +\infty$, respectively. If (U) holds, by considering the limit as $T\to +\infty$ in \eqref{lipboundU} and \eqref{eq:uniform-v-weighted}, then we can get the Lipschitz bound and the linear growth with constants independent of $N$.
\vspace{2mm}

To show that $v_\infty$ satisfies \eqref{N-PDEinf} in the classical sense, we set, for every $i\in [N]$,
\be\label{coef}
b_\infty^i(\bm x):=-D_pH^i(x^i,v_\infty^i(\bm x))\quad\text{and}\quad c_\infty^i(\bm x):=D_xH^i(x^i,v_\infty^i(\bm x))-D_{x^i}F^i(\bm x),\;\;\;\; \bm x\in (\R^d)^N.
\ee
We observe that both $b_\infty,c_\infty$ are globally Lipschitz and have at most linear growth. This implies that
\be\label{forward}
dX_s^{i,\infty}=-D_pH^i(X_s^{i,\infty},v_\infty^i(\bm X_s^\infty))\,ds
+\sqrt{2}\sigma^i\,dW_s^i,\;\;i\in [N],\qquad \bm X_{0}^\infty=\bm x\in (\R^d)^N,
\ee
is well posed.
We couple $\bm X^\infty$ synchronously with the forward equation of \eqref{N-FBSDE} starting from $\bm\xi_0=\bm x\in (\R^d)^N$ at time $0$ and horizon $T$. Let $\bm Y_s^T=v_T(s,\bm X_s^T)$ and write $\bm Y_s^\infty=v_\infty(\bm X_s^\infty)$. The proof is then divided into two steps.
\vspace{2.5mm}

\textit{Step 1.} (Turnpike estimates and uniform $L^2$-integrability)\\[0.5ex]
Let $R\in [0,T)$. In light of the global Lipschitz estimate \eqref{eq:H-gradient-Lipschitz} and Proposition~\ref{prop:uniform-Z}, the difference of the drifts of $X^{i,\infty}$ and $X^{i}$, $i\in [N]$, is bounded by a constant times $|\bm X_s^T-\bm X_s^\infty|+e^{-c(T-R)}(1+|\bm X_s^\infty|)$ on $[0,R]$. Standard arguments involving the Lipschitz continuity of $v_T,v_\infty$ together with Gronwall's inequality then imply
\[\Ex\left[\sup_{0\leq s\leq R}\left\{|\bm X_s^T-\bm X_s^\infty|^2+|\bm Y_s^T-\bm Y_s^\infty|^2\right\}\right]\leq C_{N,R}(1+|\bm x|^2)e^{-2c(T-R)},\]
for some constant $C_{N,R}$ depending only on $N,R$.
Under \textup{(U)}, the
same argument, gives
\[\Ex\left[\sup_{0\leq s\leq R}\left\{|\bm X_s^{T}-\bm X_s^{\infty}|^2+|\bm Y_s^{T}-\bm Y_s^{\infty}|^2\right\}\right]\leq C_R^*(N+|\bm x|^2)e^{-2c_*(T-R)},\]
where $C_R^*$ is independent of $N$ and $T$.
\vspace{2mm}

Recall $\alpha^{i}_\infty(s):=-D_pH^i(X_s^{i,\infty},Y_s^{i,\infty}),\;\;i\in [N]$.
Using the first estimate and passing to the limit $T\to +\infty$ at each fixed $s$ in the first estimate of Proposition~\ref{prop:uniform-moments} proves
\[\sup_{s\geq 0}\Ex\left[|\bm X_s^\infty|^2+|\bm Y_s^\infty|^2+|\bm \alpha_\infty(s)|^2\right]\leq C_N(1+|\bm x|^2).\]
The displacement monotonicity of $L^i,\;i\in[N],$ from Section \ref{not+as} and the first-order condition $D_aL^i(x,-D_pH^i(x,p))=-p$ yield
\[C_L|\bm\alpha_\infty(s)|
\leq |\bm Y_s^\infty|+K_L|\bm X_s^\infty|
+\left(\sum_{i=1}^N|D_aL^i(0,0)|^2\right)^{1/2}.\]
Consequently, under \textup{(U)}, \eqref{eq:uniform-moments-N}
and its origin bound give
\[\sup_{s\geq0}\Ex\left[|\bm X_s^{\infty}|^2+|\bm Y_s^{\infty}|^2+|\bm\alpha_\infty(s)|^2\right]\leq C_*(N+|\bm x|^2),\]
with $C_*$ independent of $N$ and $s$.
\vspace{2mm}

\textit{Step 2.} ($v_\infty$ satisfies \eqref{N-PDEinf})\\[0.5ex]
Define
\[
\bm c_T(s,\bm x)
:=
\bigl(
D_xH^i(x^i,v_T^i(s,\bm x))-D_{x^i}F^i(\bm x)
\bigr)_{i=1}^N.
\]
Taking expectations in the finite-horizon backward equation from \eqref{N-FBSDE} gives
\be\label{eq1}
v_T(0,\bm x)=\Ex\left[
v_T(r\wedge\tau,\bm X_{r\wedge\tau}^T)-\int_{0}^{r\wedge\tau}\bm c_T(s,\bm X_s^T)\,ds\right],
\ee
for any stopping time $\tau$ and $r>0$. Since $v_T$ is globally Lipschitz and has linear growth both with constants independent of $T$, the above turnpike estimates imply after sending $T\to +\infty$ in \eqref{eq1}
\be\label{eq2}
v_\infty(\bm x)=\Ex\left[v_\infty(\bm X_{r\wedge\tau}^\infty)-\int_{0}^{r\wedge\tau}\bm c_\infty(\bm X_s^\infty)\,ds\right].
\ee
Now let $B\subset (\R^d)^N$ be a bounded ball and recall the coefficients $\bm b_\infty$ and $\bm c_\infty$ from \eqref{coef} associated with $v_\infty$. We solve the scalar linear Dirichlet problem
$$\begin{cases}\sum_{j=1}^N\operatorname{tr}\left(\sigma^j(\sigma^j)^\top
D^2_{x^jx^j}u_B^i\right)
+\sum_{j=1}^ND_{x^j}u_B^i\, b_\infty^j=c_\infty^i,&\hbox{in }B,\\[1ex]
u_B^i=v_\infty^i,& \hbox{on }\partial B.
\end{cases}$$
The coefficients are locally H\"older continuous and the diffusion is uniformly elliptic. Classical elliptic regularity theory gives a continuous solution on $\overline B$ that is $C^{2,\eta}_{\mathrm{loc}}(B)$, for any $\eta\in(0,1)$.

Suppose that $\bm x\in B$ and let $\tau_B$ be the first exit time of $\bm X^\infty$, defined in \eqref{forward}, from $B$.
Uniform ellipticity and boundedness of the drift on $\overline B$ give $\Ex_{\bm x}[\tau_B]<\infty$ \footnote{Indeed, let $M:=\sup_{\overline B}|b_\infty|$, fix a unit vector $\bm n=(n^1,\ldots,n^N)\in(\R^d)^N$, and set $\psi(\bm y):=e^{\gamma\,\bm n\cdot\bm y}$ with $\gamma:=(M+1)/\nu$. By \textup{(A3)}, the generator of~\eqref{forward} applied to $\psi$ equals
$$\Bigl(\gamma^2\sum_{j=1}^N|(\sigma^j)^\top n^j|^2+\gamma\,b_\infty\cdot\bm n\Bigr)\psi\geq\gamma(\gamma\nu-M)\psi=\gamma\psi,\quad\text{on }\overline B.$$
Taking expectations in It\^o's formula for $\psi(\bm X^\infty)$ up to $t\wedge\tau_B$ therefore gives $\gamma\,\Ex_{\bm x}[t\wedge\tau_B]\min_{\overline B}\psi\leq\max_{\overline B}\psi$ for every $t>0$, and letting $t\to\infty$ monotone convergence proves the claim.}. Taking expectations in It\^o's formula yields
$$u_B(\bm x)=\mathbb{E}\left[v_\infty(\bm X^\infty_{\tau_B})-\int_{0}^{\tau_B}\bm c_\infty(\bm X^\infty_s)ds \right].$$
Combining this with \eqref{eq2} gives $u_B=v_\infty$ in $B$. The arbitrariness of $B$ gives that $v_\infty$ satisfies \eqref{N-PDEinf} in the classical sense, which finishes the proof of Step 2.
\vspace{3mm}

We finally show the uniqueness statement. Suppose that $w\colon (\R^d)^N\to (\R^d)^N$ is a globally Lipschitz solution of \eqref{N-PDEinf} satisfying \eqref{condition}. We couple its induced forward diffusion (recall \eqref{forward''}) with the induced forward diffusion of $v_\infty$ (given by \eqref{forward}) with the same Brownian motions and initial states $\bm\xi_0=\bm x\in (\R^d)^N$. We set
\begin{align*}
&\Delta\bm X_s= \bm X^\infty_s-\bm X^w_s,\qquad \Delta \bm Y_s=\bm Y^\infty_s-\bm Y^w_s=v_\infty(\bm X^\infty_s)-w(\bm X_s^w),\\
&\hspace{1cm}\Delta\bm\alpha=\left(-D_pH^i(X_s^{i,\infty},Y_s^{i,\infty})+D_pH^i(X_s^{i,w},Y_s^{i,w})\right)_{i=1}^N.
\end{align*}
We also let
\[
q(s):=\Ex[\Delta\bm X_s\cdot\Delta\bm Y_s]\quad\text{and}\quad
e(s):=\Ex\left[|\Delta\bm X_s|^2+|\Delta\bm\alpha_s|^2\right],\;\;s\ge 0.
\]
 By It\^o's formula and~\eqref{N-PDEinf}, the state $\bm X^\infty$ with $\bm Y=v_\infty(\bm X^\infty)$ and $Z^{i,j}=D_{x^j}v^i_\infty(\bm X^\infty)\sigma^j$, and the state $\bm X^w$ with $\bm Y=w(\bm X^w)$ and $Z^{i,j}=D_{x^j}w^i(\bm X^w)\sigma^j$, satisfy the equations of Lemma~\ref{lem:pairing} with $\bm\phi^k=\bm\psi^k=0$ on every interval $[0,b]$; the integrability conditions hold because $v_\infty$ and $w$ are Lipschitz. Lemma~\ref{lem:pairing}(i)--(ii) gives
\[
q(0)=0,\qquad q'(s)\leq-C_Le(s),\qquad
|q(s)|\leq\frac{1+K_L^2}{2}e(s),
\]
for almost every $s\ge 0$. Therefore $q$ is nonpositive and satisfies
$q'(s)\leq2C_Lq(s)/(1+K_L^2)$, almost everywhere. Assumption \eqref{condition} and the linear growth of both fields $w,v_\infty$ make $q$ uniformly
bounded on $[0,\infty)$.

Suppose that there exists $s_0>0$ such that $q(s_0)<0$. Letting $u=\frac{2C_L}{1+K_L^2}$, we have that $q$ satisfies
$$q(s)\le e^{u(s-s_0)}q(s_0),\;\text{ for any }s\ge s_0.$$
This means that $q$ is unbounded and we get a contradiction. Thus, $q(s)=0$, for all $s\ge 0$ and this yields $e(s)=0$, or
$$\bm X_s^\infty=\bm X^w_s\quad\text{and}\quad D_pH^i(X_s^{i,\infty},Y_s^{i,\infty})=D_pH^i(X_s^{i,w},Y_s^{i,w}),\;\;i\in [N],$$
for any $s\ge 0$. Setting $s=0$, we get from the injectivity\footnote{This is true because $D_aL^i(x,-D_pH^i(x,p))=-p$, for any $x,p\in \R^d$ and $i\in [N]$.} of $p\mapsto D_pH^i(x,p)$ that $v_\infty(\bm x)=w(\bm x)$. The proof is complete.\qed
\vspace{3mm}

We are now ready to prove Theorem \ref{thm:turnpike}.

\subsubsection*{Proof of Theorem~\ref{thm:turnpike}}\phantomsection\label{pf:turnpike}
Couple the horizon-$T$ forward process and $\bm X^\infty$
synchronously, starting from the same square-integrable
$\mathcal I$-measurable initial configuration $\bm\xi_0$ at $0$.
The coupling argument in the preceding proof of Theorem \ref{thm:stationary-limit} (Step 1), conditioned on $\mathcal I$, gives
\[
\Ex\left[\sup_{0\leq s\leq R}
\left\{|\bm X_s^T-\bm X_s^\infty|^2
+|\bm Y_s^T-\bm Y_s^\infty|^2\right\}\right]
\leq C_{N,R}(1+\Ex[|\bm\xi_0|^2])e^{-2c(T-R)}.
\]
Indeed, \eqref{goal:estimate2}, the uniform spatial Lipschitz bound,
and \eqref{eq:H-gradient-Lipschitz} bound the drift difference;
Gronwall's inequality on $[0,R]$ then bounds the state difference.
The identity $\bm Y_s^T=v_T(s,\bm X_s^T)$ gives the second estimate. Under \textup{(U)}, the same argument uses
\eqref{goal:estimate2-uniform} and the constant $K_*$ to give
\be\label{turn}
\Ex\left[\sup_{0\leq s\leq R}
\left\{|\bm X_s^{T}-\bm X_s^{\infty}|^2
+|\bm Y_s^{T}-\bm Y_s^{\infty}|^2\right\}\right]
\leq C_R^*(N+\Ex[|\bm\xi_0|^2])e^{-2c_*(T-R)},
\ee
with $C_R^*$ independent of $N$ and $T$.
The control estimate for fixed $N$ follows directly from
\eqref{eq:H-gradient-Lipschitz}:
\[
|\bm\alpha_T(s)-\bm\alpha_\infty(s)|
\leq\frac{K_L}{C_L}|\bm X_s^T-\bm X_s^\infty|
+\frac1{C_L}|\bm Y_s^T-\bm Y_s^\infty|.
\]
Under \textup{(U)}, the constants $K_L/C_L$ and $1/C_L$ in the preceding bound are at most $\overline K_L/\overline C_L$ and $1/\overline C_L$. Squaring, taking the supremum over $s\in[0,R]$ and expectations, and combining the control bound with the first display of this proof gives \eqref{eq:turnpike-fixed-window}; combining it with \eqref{turn} gives \eqref{eq:turnpike-fixed-window-N}, after increasing the constants.\qed

\subsection{Monotonicity property of the Pontryagin fields} As a corollary of our proof methods for the above Propositions and Theorem \ref{thm:stationary-limit}, we establish the following monotonicity property for the Pontryagin fields $v_T,v_\infty$.

\begin{lemma}[Monotonicity property]
\label{lem:erg-strong-monotonicity}
Under Assumptions \textup{(A1)-(A3)}, the finite-horizon field $v_T$ satisfies the monotonicity property
\begin{equation}\label{eq:erg-finite-monotonicity}
\bigl(v_T(s,\bm x)-v_T(s,\bm y)\bigr)\cdot(\bm x-\bm y)
\geq0,\;\;\text{for all }\bm x,\bm y\in (\R^d)^N\text{ and }s\in [0,T].
\end{equation}
Moreover, the stationary field $v_\infty$ constructed in Theorem \ref{thm:stationary-limit} satisfies
\be\label{eq:erg-strong-monotonicity}
\bigl(v_\infty(\bm x)-v_\infty(\bm y)\bigr)\cdot(\bm x-\bm y)
\geq C_L|\bm x-\bm y|^2,\;\;\text{for all }\bm x,\bm y\in(\R^{d})^N.
\ee
In particular, $C_L\le K$, where $K$ is the Lipschitz constant of $v_\infty$ obtained in Proposition \ref{prop:uniform-Z}.
Under \textup{(U)}, one may take $C_L=\overline C_L$ in \eqref{eq:erg-strong-monotonicity}, so it holds with a constant independent of $N$.
\end{lemma}

\begin{proof}
The first statement is immediate. Inequality \eqref{eq:vT-pairing} from the proof of
Proposition~\ref{prop:uniform-Z} gives
\[
\bigl(v_T(s,\bm x)-v_T(s,\bm y)\bigr)\cdot(\bm x-\bm y)
\geq C_L\int_s^T\Ex\bigl[|\Delta\bm X_t|^2+|\Delta\bm\alpha_t|^2\bigr]dt
\geq0,
\]
where $\Delta$ refers there to the synchronously coupled finite-horizon
solutions started from $\bm x$ and $\bm y$ at time $s$.

For \eqref{eq:erg-strong-monotonicity}, we couple the solutions $\bm X^{\bm x}:=\bm X^{\infty,\bm x},\;\bm X^{\bm y}:=\bm X^{\infty, \bm y}$ of the forward problem \eqref{forward'} started at time $0$ from $\bm x\in (\R^d)^N$ and
$\bm y\in (\R^d)^N$, respectively. We denote by $\Delta \bm X$ the difference of the states and by $\Delta \bm \alpha$ the difference of their drifts $b_\infty(\bm X_r^{\bm x})-b_\infty(\bm X_r^{\bm y})$, where we used the notation from \eqref{coef}. We also set for $r\ge 0$
\[
f(r):=\Ex[|\Delta\bm X_r|^2],\qquad
g(r):=\Ex[|\Delta\bm\alpha_r|^2],\;\;\text{and}\;\;\;
q(r):=\Ex\bigl[\Delta\bm X_r\cdot
\bigl(v_\infty(\bm X_r^{\bm x})
-v_\infty(\bm X_r^{\bm y})\bigr)\bigr].
\]

By \eqref{goal:estimate2},
$v_T(0,\bm z)\to v_\infty(\bm z)$ as $T\to\infty$, for every
$\bm z\in(\R^d)^N$. Taking $s=0$ in
\eqref{eq:erg-finite-monotonicity} and letting $T\to\infty$ yields
\begin{equation}\label{eq:vinf-monotone}
\bigl(v_\infty(\bm x)-v_\infty(\bm y)\bigr)\cdot(\bm x-\bm y)\geq0,
\qquad \bm x,\bm y\in(\R^d)^N.
\end{equation}
Since $b_\infty$ is globally Lipschitz, \eqref{forward'} has a unique
strong solution for every square-integrable $\mathcal I$-measurable
initial state $\bm\xi$, and standard SDE estimates give, for every
$R>0$,
\begin{equation}\label{eq:canonical-finite-moments}
\Ex\Bigl[\sup_{0\leq r\leq R}|\bm X^{\infty,\bm\xi}_r|^2\Bigr]
\leq C_{N,R}\bigl(1+\Ex[|\bm\xi|^2]\bigr).
\end{equation}
By It\^o's formula and~\eqref{N-PDEinf}, the states $\bm X^{\bm x}$ and $\bm X^{\bm y}$, with $\bm Y=v_\infty(\bm X)$ and $Z^{i,j}=D_{x^j}v^i_\infty(\bm X)\sigma^j$ for $\bm X\in\{\bm X^{\bm x},\bm X^{\bm y}\}$, satisfy the equations of Lemma~\ref{lem:pairing} with $\bm\phi^k=\bm\psi^k=0$ on every interval $[0,R]$; the integrability conditions hold by~\eqref{eq:canonical-finite-moments}, the linear growth of $v_\infty$ and the boundedness of $D_{\bm x}v_\infty$. The process $\Delta\bm\alpha$ defined above is the one of Lemma~\ref{lem:pairing}, so the function $e$ there equals $f+g$, and Lemma~\ref{lem:pairing}(i) gives
\begin{equation}\label{eq:erg-pairing-energy}
q'(r)\leq-C_L\bigl(f(r)+g(r)\bigr),
\quad\text{for almost every }r\geq 0.
\end{equation}
Applying \eqref{eq:vinf-monotone} pointwise at
$(\bm X^{\bm x}_r,\bm X^{\bm y}_r)$ shows $q(r)\geq0$, for every
$r\geq 0$. Integrating \eqref{eq:erg-pairing-energy} over $[0,R]$
therefore gives
\[
C_L\int_{0}^R\bigl(f(r)+g(r)\bigr)dr\leq q(0)-q(R)\leq q(0),
\]
and, letting $R\to\infty$,
\begin{equation}\label{eq:energy-integrable}
C_L\int_{0}^\infty\bigl(f(r)+g(r)\bigr)dr
\leq(\bm x-\bm y)\cdot\bigl(v_\infty(\bm x)-v_\infty(\bm y)\bigr)<\infty.
\end{equation}
On the other hand, since $\Delta\bm X$ is absolutely continuous with
derivative $\Delta\bm\alpha$,
\[
f(r)=f(0)+2\int_{0}^r\Ex\bigl[\Delta\bm X_u\cdot\Delta\bm\alpha_u\bigr]du\ge f(0)-\int_{0}^r2\sqrt{f(u)g(u)}\;du\ge |\bm x-\bm y|^2-\int_{0}^r(f(u)+g(u))du,
\]
for any $r\ge 0$. By \eqref{eq:energy-integrable}, $f\in L^1(0,\infty)$, so there is a
sequence $r_n\to\infty$ with $f(r_n)\to0$. Evaluating the last
inequality along $r_n$ gives
$|\bm x-\bm y|^2\leq\int_{0}^\infty(f(r)+g(r))\,dr$, and combining this
with \eqref{eq:energy-integrable} proves
\eqref{eq:erg-strong-monotonicity}.
\vspace{2mm}

Under \textup{(U)}, \textup{(A1)} holds with $C_L=\overline C_L$ for
every $N$, and the argument above is unchanged.
\end{proof}

\section{Average-cost game interpretation of the limiting dynamics}\label{6}

\subsection{The limiting control is an average-cost equilibrium}
We start by showing that the process $\bm X^\infty$ from \eqref{forward'} is associated with an open-loop Nash equilibrium of the average-cost $N$-player game.

\subsubsection*{Proof of Theorem \ref{prop:average-Nash}}\phantomsection\label{pf:average-Nash}
Recall the estimates from Theorems \ref{thm:stationary-limit}, \ref{thm:turnpike}. The linear growth of $D_pH^i$ and the uniform moment bounds imply square integrability of $\bm\alpha_\infty$ and bounded second moments for $X^{i,\infty}_t,\;t\ge 0$, for all $i\in [N]$, hence $\mathcal{E}(\bm \alpha_\infty)<\infty$. Thus the profile $\bm \alpha_\infty\in \mathcal{A}_\infty(\bm x)$, and its expected average costs $A^i_T(\bm \alpha_\infty)$, for $T>0,\;i\in [N]$, are finite and uniformly bounded in absolute value.

Fix a player $i\in [N]$ and a unilateral deviation $\beta^i\in \mathcal{A}^i_{\infty}(\bm x)$. We write $\Delta X_s=X_s^{i,\beta^i}-X_s^{i,\infty}$ and $\Delta\alpha_s=\beta_s^i-\alpha_\infty^{i}(s)$. Since the opponents' open-loop controls are fixed, their states coincide in the two profiles. Let $R>0$ and set
$$f_s=D_xL^i(X_s^{i,\infty},\alpha_\infty^{i}(s))
+D_{x^i}F^i(\bm X_s^\infty).$$
The identity $D_aL^i(X_s^{i,\infty},\alpha_\infty^{i}(s))=-Y_s^{i,\infty},$ together with joint convexity of $L^i$ and own-coordinate convexity of $F^i$, implies
\begin{align}
R\left(A_R^i([\bm\alpha^{-i}_\infty;\beta^i])-A_R^i(\bm\alpha_\infty)\right)&\geq\int_{0}^R
\Ex\left[f_s\cdot\Delta X_s-Y_s^{i,\infty}\cdot\Delta\alpha_s\right]ds+\frac{C_L}{2}\int_{0}^R\Ex\left[|\Delta X_s|^2+|\Delta\alpha_s|^2\right]ds\nonumber\\
&=-\Ex[Y_R^{i,\infty}\cdot\Delta X_R]
+\frac{C_L}{2}\int_{0}^R\Ex\left[|\Delta X_s|^2+|\Delta\alpha_s|^2\right]ds\nonumber\\
&\ge -\Ex[Y_R^{i,\infty}\cdot\Delta X_R].\label{ineq101}
\end{align}
The equality follows from It\^o's formula on $Y^{i,\infty}_s\cdot \Delta X_s$, using $\Delta X_{0}=0$, $d\Delta X_s=\Delta\alpha_sds$. We also observe that
$$\frac{|\Ex[Y_R^{i,\infty}\cdot\Delta X_R]|}{R}
\leq\left(\sup_{s\geq{0}}\Ex[|Y_s^{i,\infty}|^2]\right)^{\frac{1}{2}}
\frac{(\Ex[|\Delta X_R|^2])^{\frac{1}{2}}}{R}\le \left(C(1+|\bm x|^2)\right)^{\frac{1}{2}}\frac{(\Ex[|\Delta X_R|^2])^{\frac{1}{2}}}{R}\xrightarrow{R\to +\infty}0,$$
where we used the linear growth of $v_\infty$ and the uniform square integrability of $\bm X^\infty$ to bound the first factor, and Lemma \ref{rem:erg-B-admissible} for the limit as $R\to +\infty$.
We divide the inequality \eqref{ineq101} by $R>0$ and take the lower limit as $R\to +\infty$.
The result follows.\qed
\vspace{3mm}

The next proposition provides a way to characterize strong average-cost equilibria.

\begin{proposition}\label{thm:erg-uniqueness}
Suppose that \textup{(A1)-(A3)} hold and let $\bm x\in (\R^d)^N$.
A profile $\bm\beta\in\mathcal{A}_\infty(\bm x)$ is a strong
average-cost equilibrium if and only if
\begin{equation}\label{erg:cesaro}
\lim_{R\to\infty}\frac1R\Ex\left[\int_0^R
\left(|\bm X_s^{\bm\beta}-\bm X_s^\infty|^2
+|\bm\beta_s-\bm\alpha_\infty(s)|^2\right)ds\right]=0.
\end{equation}
Every such profile satisfies
\begin{equation}\label{erg:same-cost}
\lim_{R\to\infty}
\bigl(A_R^i(\bm\beta)-A_R^i(\bm\alpha_\infty)\bigr)=0,\;\;\text{for any }i\in [N].
\end{equation}
\end{proposition}

\begin{proof}
Suppose that $\bm\beta$ is a strong average-cost equilibrium. Set
$\Delta\bm X=\bm X^{\bm\beta}-\bm X^\infty$ and
$\Delta\bm a=\bm\beta-\bm\alpha_\infty$. Apply
Lemma~\ref{lem:erg-first-order} at each of the two profiles with the
other profile's control as the deviation. Subtracting the resulting
zero limits, and summing over $i$, yields
\begin{align*}
0=\lim_{R\to\infty}\frac1R\int_0^R\sum_{i=1}^N\Ex
\bigl[
&\bigl(D_xL^i(X_s^{i,\bm\beta},\beta_s^i)
-D_xL^i(X_s^{i,\infty},\alpha_\infty^{i}(s))\bigr)
\cdot\Delta X_s^i\\
&\;\;\;+\bigl(D_aL^i(X_s^{i,\bm\beta},\beta_s^i)
-D_aL^i(X_s^{i,\infty},\alpha_\infty^{i}(s))\bigr)
\cdot\Delta a_s^i\\
&\hspace{3cm}+\bigl(D_{x^i}F^i(\bm X_s^{\bm\beta})
-D_{x^i}F^i(\bm X_s^\infty)\bigr)\cdot\Delta X_s^i
\bigr]ds.
\end{align*}
By \eqref{Ldispl} and \eqref{Fmon}, the integrand is at least
$C_L(|\Delta\bm X_s|^2+|\Delta\bm a_s|^2)$, proving
\eqref{erg:cesaro}.

Conversely, suppose \eqref{erg:cesaro} holds. Bounded Hessians give
the weighted Lipschitz bounds
\begin{equation}\label{erg:growth}
\begin{split}
|L^i(z)-L^i(z')|\leq C(1+|z|+|z'|)|z-z'|\;\;\text{and}\;\;
|F^i(\bm y)-F^i(\bm y')|
\leq C(1+|\bm y|+|\bm y'|)|\bm y-\bm y'|.
\end{split}
\end{equation}
Cauchy--Schwarz and the mean-energy bounds therefore prove
\eqref{erg:same-cost}. Fix $\gamma\in\mathcal{A}_\infty^i(\bm x)$ and
write the cost difference as
\begin{align*}
A_R^i([\bm\beta^{-i};\gamma])-A_R^i(\bm\beta)
=\bigl(A_R^i([\bm\beta^{-i};\gamma])
-&A_R^i([\bm\alpha^{-i}_\infty;\gamma])\bigr)\\
&+\bigl(A_R^i([\bm\alpha^{-i}_\infty;\gamma])
-A_R^i(\bm\alpha_\infty)\bigr) +\bigl(A_R^i(\bm\alpha_\infty)-A_R^i(\bm\beta)\bigr).
\end{align*}
The last term tends to zero. In the first term, player $i$ has the
same state and control in both costs, so only $F^i$ changes. Its
absolute value is bounded by
\[
C\bigl(1+\mathcal E(\bm\beta)+\mathcal E(\bm\alpha_\infty)
+\mathcal E^i(\gamma)\bigr)^{1/2}
\left(\frac1R\int_0^R\Ex[|\Delta\bm X_s|^2]ds\right)^{1/2},
\]
and also tends to zero. The middle term has nonnegative lower limit
by Theorem \ref{prop:average-Nash}. This proves \eqref{erg:strong} and completes the proof.
\end{proof}

We note that $\bm\alpha_\infty$ is not the unique strong average-cost equilibrium.

\begin{remark}[Non-uniqueness of strong average-cost equilibria]\label{rem:erg-nonunique}
Choose vectors $\zeta^1,\ldots,\zeta^N\in\R^d$, not all zero, and
define the progressively measurable profile
\[
\beta_s^j=\alpha_\infty^{j}(s)
+\zeta^j\bigl(\ind_{[0,1]}(s)-\ind_{(1,2]}(s)\bigr),\;\;j\in [N].
\]
The controls here are the canonical open-loop processes plus a pulse,
not the canonical rule evaluated at the perturbed states. Its state
difference is
\[
X_s^{j,\beta^j}-X_s^{j,\infty}
=\zeta^j\int_0^s
\bigl(\ind_{[0,1]}(u)-\ind_{(1,2]}(u)\bigr)du.
\]
Thus both profiles and states agree for $s\geq2$, and
$\bm\beta\in\mathcal{A}_\infty(\bm x)$ satisfies \eqref{erg:cesaro}, hence by Proposition \ref{thm:erg-uniqueness}, it is a distinct strong equilibrium with the same long-run costs.
\end{remark}

\subsection{Contraction and invariant law}
All constants in this section may depend on $N$.
The canonical state process $(\bm X^{\infty, \bm x}_s)_{s\ge 0}$ from \eqref{forward'} satisfies, by the proof of Theorem \ref{thm:stationary-limit},
\begin{equation}\label{eq:erg-canonical-moments}
\sup_{s\geq0}\Ex [|\bm X_s^{\infty,\bm x}|^2]
\leq C_N(1+|\bm x|^2),\;\;\bm x\in (\R^d)^N.
\end{equation}
The same estimate holds with $\Ex[|\bm\xi|^2]$ in place of
$|\bm x|^2$ for a square-integrable initial state $\bm\xi$
independent of the future Brownian increments.
\vspace{2mm}

Let $P_s$ be the Markov semigroup of the process $\bm X^{\infty}$ from \eqref{forward'}: $P_s\varphi(\bm x):=\Ex\bigl[\varphi(\bm X_s^{\infty,\bm x})\bigr]$. If the process $\bm X^\infty$ starts from distribution $\nu\in \mathcal{P}_2(\R^{Nd})$, then $\nu P_s$, $s\ge 0$, is the distribution of $\bm X^\infty_s$.

\begin{proposition}[Contraction and invariant distribution]\label{prop:erg-contraction}
Assume \textup{(A1)-(A3)} and let $\kappa:=\frac{C_L}{K}$. For square-integrable initial states $\bm\xi,\bm\eta$, coupled
independently of the future Brownian increments, synchronous
solutions of \eqref{forward'} satisfy
\begin{equation}\label{eq:erg-synchronous-contraction}
\Ex\left[|\bm X_s^{\infty,\bm\xi}
-\bm X_s^{\infty,\bm\eta}|^2\right]
\leq\frac{1}{\kappa}e^{-\kappa s}
\Ex[|\bm\xi-\bm\eta|^2],
\quad s\geq0.
\end{equation}
In particular, there is a unique $\mu_\infty^N\in\mathcal P_2(\R^{Nd})$ such that
\begin{equation}\label{eq:erg-Wasserstein-contraction}
W_2(\nu P_s,\mu_\infty^N)
\leq\frac{1}{\sqrt{\kappa}}e^{-\frac{\kappa}{2} s}
W_2(\nu,\mu_\infty^N),
\qquad\nu\in\mathcal P_2(\R^{Nd}).
\end{equation}
The measure $\mu_\infty^N$ has full support and is called the invariant measure, because $\mu_\infty^NP_s=\mu_\infty^N$, $s\ge 0$.
\end{proposition}

\begin{proof}
Let $\bm X^{\bm\xi}:=\bm X^{\infty,\bm\xi}$ and $\bm X^{\bm\eta}:=\bm X^{\infty,\bm\eta}$ solve~\eqref{forward'} from time $0$ with the same Brownian motions. Define $\Delta\bm X$, $f$, $g$ and $q$ as in the proof of Lemma~\ref{lem:erg-strong-monotonicity}, with $(\bm\xi,\bm\eta)\in L^2\times L^2$ in place of $(\bm x,\bm y)$. By~\eqref{eq:canonical-finite-moments}, which holds for square-integrable initial states, Lemma~\ref{lem:pairing} applies as in the proof of~\eqref{eq:erg-pairing-energy}, so $q$ is absolutely continuous and
\[
q'(s)\leq-C_L\bigl(f(s)+g(s)\bigr)\leq-C_Lf(s),\quad\text{for almost every }s\geq0.
\]
Applying~\eqref{eq:erg-strong-monotonicity} and the $K$-Lipschitz bound for $v_\infty$ at the points $\bm X_s^{\bm\xi}$ and $\bm X_s^{\bm\eta}$ gives $C_Lf(s)\leq q(s)\leq Kf(s)$, for all $s\geq0$. Therefore $q'(s)\leq-\kappa q(s)$ for almost every $s\geq0$, so $q(s)\leq e^{-\kappa s}q(0)$ and
\[
f(s)\leq\frac{q(s)}{C_L}\leq\frac{e^{-\kappa s}}{C_L}\,q(0)\leq\frac{K}{C_L}e^{-\kappa s}f(0)=\frac1\kappa e^{-\kappa s}\Ex[|\bm\xi-\bm\eta|^2].
\]
This is~\eqref{eq:erg-synchronous-contraction}.

For the second statement we compare two laws $\nu,\nu'\in\mathcal P_2(\R^{Nd})$. By~\eqref{eq:erg-canonical-moments}, $\nu P_s$ and $\nu'P_s$ belong to $\mathcal P_2(\R^{Nd})$. Let $(\bm \xi,\bm \eta)\sim\pi\in\Pi(\nu,\nu')$ be an optimal coupling for $W_2(\nu,\nu')$. We have
\begin{equation}\label{eq:erg-W2-two-laws}
W_2^2(\nu P_s,\nu'P_s)\leq\Ex\bigl[|\bm X_s^{\infty,\bm \xi}-\bm X_s^{\infty,\bm \eta}|^2\bigr]\leq \frac1\kappa e^{-\kappa s}\mathbb{E}[|\bm \xi -\bm \eta|^2]=\frac1\kappa e^{-\kappa s}W_2^2(\nu,\nu'),\quad s\geq0.
\end{equation}
Choose $s_0>0$ such that $\kappa^{-1/2}e^{-\kappa s_0/2}<1$. By~\eqref{eq:erg-W2-two-laws}, the map $\nu\mapsto\nu P_{s_0}$ is a strict contraction of the complete metric space $(\mathcal P_2(\R^{Nd}),W_2)$, so it has a unique fixed point $\mu_\infty^N$. For every $s\geq0$, the semigroup property gives $(\mu_\infty^NP_s)P_{s_0}=(\mu_\infty^NP_{s_0})P_s=\mu_\infty^NP_s$. Thus $\mu_\infty^NP_s$ is also a fixed point, and $\mu_\infty^NP_s=\mu_\infty^N$. Taking $\nu'=\mu_\infty^N$ in~\eqref{eq:erg-W2-two-laws} proves~\eqref{eq:erg-Wasserstein-contraction}. If $\mu'\in\mathcal P_2(\R^{Nd})$ also satisfies~\eqref{eq:erg-Wasserstein-contraction}, then the choice $\nu=\mu_\infty^N$ there gives $W_2(\mu_\infty^N,\mu')\leq\kappa^{-1/2}e^{-\kappa s/2}W_2(\mu_\infty^N,\mu')$ for every $s\geq0$, hence $\mu'=\mu_\infty^N$.

It remains to show that $\mu_\infty^N$ has full support. Fix $s>0$, $\bm x\in\R^{Nd}$ and a nonempty open set $O\subset\R^{Nd}$. The drift $b_\infty$ in~\eqref{coef} has at most linear growth and $\Sigma=\operatorname{diag}(\sigma^1,\ldots,\sigma^N)$ is invertible. Girsanov's theorem therefore applies on $[0,s]$ and shows that the law of $\bm X_s^{\infty,\bm x}$ is equivalent to the law of $\bm x+\sqrt2\Sigma\bm W_s$, which is a nondegenerate Gaussian law. Hence $\Pr(\bm X_s^{\infty,\bm x}\in O)>0$, for all $\bm x\in \R^{Nd}$, and the invariance of $\mu_\infty^N$ gives
\[
\mu_\infty^N(O)=\int\Pr\bigl(\bm X_s^{\infty,\bm x}\in O\bigr)\,\mu_\infty^N(d\bm x)>0.\hspace{3cm}\qedhere
\]
\end{proof}

\subsection{Ergodic constants and ergodic values} To construct the ergodic constants we recall the notation that we introduced in the previous sections.
\begin{align}
b^i(\bm x)=-D_pH^i(&x^i,v_\infty^i(\bm x)),\qquad
b_T^i(s,\bm x)=-D_pH^i(x^i,v_T^i(s,\bm x)),
\label{eq:erg-canonical-drifts}\\
&c^i(\bm x)=D_xH^i(x^i,v_\infty^i(\bm x))-D_{x^i}F^i(\bm x).
\label{eq:erg-canonical-backward-drift}
\end{align}
Write $b=(b^i)_{i=1}^N$, $\bm c=(c^i)_{i=1}^N$, and
$\Sigma=\operatorname{diag}(\sigma^1,\ldots,\sigma^N)$.
The process $\bm X^\infty$ with
\begin{equation}\label{eq:erg-canonical-profile}
d\bm X_s^\infty=b(\bm X_s^\infty)\,ds
 +\sqrt{2}\Sigma\,d\bm W_s,
\qquad \bm X_0^\infty=\bm x,
\qquad \bm\alpha_\infty(s)=b(\bm X_s^\infty)
\end{equation}
is referred to as the canonical state. Its generator is denoted by
\begin{equation}\label{eq:erg-generator}
\mathcal L_\infty\varphi(\bm x)
=\sum_{j=1}^N\operatorname{tr}\!\left(
 \sigma^j(\sigma^j)^\top D_{x^jx^j}^2\varphi(\bm x)\right)
 +\sum_{j=1}^N b^j(\bm x)\cdot D_{x^j}\varphi(\bm x).
\end{equation}
For $i\in[N]$, we define the running cost
\begin{equation}\label{eq:erg-running-cost}
r^{i,N}(\bm x)
:=L^i(x^i,b^i(\bm x))+F^i(\bm x).
\end{equation}
Bounded Hessians of $L^i,F^i$ and the global Lipschitz property of
$b$ \footnote{The gradients of $L^i$ and $F^i$ have at most linear growth, and $|b(\bm x)|\leq C_N(1+|\bm x|)$.} imply
\begin{equation}\label{eq:erg-running-cost-bounds}
\begin{split}
|r^{i,N}(\bm x)|\leq C_{N,i}(1+|\bm x|^2)\;\;\;\text{and}\;\;\;
|r^{i,N}(\bm x)-r^{i,N}(\bm y)|\leq C_{N,i}(1+|\bm x|+|\bm y|)|\bm x-\bm y|.
\end{split}
\end{equation}

\subsubsection*{Proof of Theorem \ref{thm:main-erg}}\phantomsection\label{pf:main-erg}
Theorem \ref{thm:main-erg} follows directly from Proposition~\ref{thm:erg-value} below, where we introduce the ergodic constants and show that they are the values of the average-cost game at equilibrium.

\begin{proposition}\label{thm:erg-value}
Assume \textup{(A1)-(A3)} and consider $\mu_\infty^N$ the invariant measure constructed in Proposition \ref{prop:erg-contraction}. For each player $i\in [N]$, the constant $\lambda^{i,N}:=\int r^{i,N}\,d\mu_\infty^N$ is finite and does not depend on the initial state. For the canonical state
started from $\bm x$,
\begin{equation}\label{eq:erg-value-rate}
\left|\Ex\left[r^{i,N}(\bm X_s^{\infty,\bm x})\right]
-\lambda^{i,N}\right|
\leq C_{N,i}(1+|\bm x|^2)e^{-\frac{\kappa}{2} s},
\qquad
\left|A_T^i(\bm\alpha_\infty)-\lambda^{i,N}\right|
\leq\frac{C_{N,i}(1+|\bm x|^2)}{T}.
\end{equation}
Every strong average-cost equilibrium
$\bm\beta\in\mathcal{A}_\infty(\bm x)$ has the same value for player $i$:
\begin{equation}\label{eq:erg-all-equilibrium-values}
\lim_{T\to\infty}A_T^i(\bm\beta)=\lambda^{i,N}.
\end{equation}
\end{proposition}

\begin{proof}
Fix $i\in[N]$. The quadratic bound in~\eqref{eq:erg-running-cost-bounds} and $\mu_\infty^N\in\mathcal P_2(\R^{Nd})$ show that $\lambda^{i,N}$ is finite. It does not depend on the initial state, because $\mu_\infty^N$ does not.

We first prove an estimate for the semigroup that will be used again in Section~\ref{7}. Let $\varphi\colon\R^{Nd}\to\R$ and $C_\varphi\geq0$ satisfy
\[
|\varphi(\bm y)-\varphi(\bm z)|\leq C_\varphi(1+|\bm y|+|\bm z|)|\bm y-\bm z|,\qquad\bm y,\bm z\in\R^{Nd}.
\]
Then there is a constant $C$, independent of $\varphi$, $\bm x$ and $s$, such that
\begin{equation}\label{eq:erg-semigroup-cost-rate}
\left|P_s\varphi(\bm x)-\int\varphi\,d\mu_\infty^N\right|\leq C\,C_\varphi(1+|\bm x|^2)e^{-\frac{\kappa}{2}s},\qquad\bm x\in\R^{Nd},\ s\geq0.
\end{equation}
To see this, let $\pi_s\in\Pi(\delta_{\bm x}P_s,\mu_\infty^N)$ be an optimal coupling for $W_2$. The Cauchy--Schwarz inequality gives
\begin{align*}
\left|P_s\varphi(\bm x)-\int\varphi\,d\mu_\infty^N\right|
&\leq\int|\varphi(\bm y)-\varphi(\bm z)|\,\pi_s(d\bm y,d\bm z)\\
&\leq C_\varphi\left(\int(1+|\bm y|+|\bm z|)^2\,\pi_s(d\bm y,d\bm z)\right)^{1/2}W_2(\delta_{\bm x}P_s,\mu_\infty^N).
\end{align*}
By~\eqref{eq:erg-canonical-moments} and $\mu_\infty^N\in\mathcal P_2(\R^{Nd})$, the first factor is at most $C(1+|\bm x|)$. By~\eqref{eq:erg-Wasserstein-contraction} with $\nu=\delta_{\bm x}$,
\[
W_2(\delta_{\bm x}P_s,\mu_\infty^N)\leq\frac{1}{\sqrt{\kappa}}e^{-\frac{\kappa}{2}s}\left(\int|\bm x-\bm z|^2\,\mu_\infty^N(d\bm z)\right)^{1/2}\leq C(1+|\bm x|)e^{-\frac{\kappa}{2}s}.
\]
Multiplying the two bounds and using $(1+|\bm x|)^2\leq2(1+|\bm x|^2)$ proves~\eqref{eq:erg-semigroup-cost-rate}.

By~\eqref{eq:erg-running-cost-bounds}, we may take $\varphi=r^{i,N}$ in~\eqref{eq:erg-semigroup-cost-rate}. Since $P_sr^{i,N}(\bm x)=\Ex[r^{i,N}(\bm X_s^{\infty,\bm x})]$ and $\int r^{i,N}\,d\mu_\infty^N=\lambda^{i,N}$, this is the first estimate in~\eqref{eq:erg-value-rate}. To get the second estimate in \eqref{eq:erg-value-rate}, we integrate the first for $s\in [0,T]$ and we divide by $T$ (and rename the constant $C_{N,i}$).
Finally, let $\bm\beta\in\mathcal{A}_\infty(\bm x)$ be a strong average-cost equilibrium. By~\eqref{erg:same-cost} in Proposition~\ref{thm:erg-uniqueness}, $A_T^i(\bm\beta)-A_T^i(\bm\alpha_\infty)\to0$ as $T\to\infty$, and the estimate just proved gives~\eqref{eq:erg-all-equilibrium-values}.
\end{proof}

\section{Long-time behavior of the finite-horizon $N$-player game}\label{7}

We divide this last section of the paper into two parts. We first construct the correctors and then we prove the long-time behavior results.

\subsection{Construction of the correctors}

\begin{proposition}[Poisson correctors]
\label{prop:erg-poisson}
Assume \textup{(A1)-(A3)}. For each $i\in[N]$, there exists a function $u^{i,N}\colon (\R^d)^N\to \R$ which is in $C^2$, satisfies
\begin{equation}\label{eq:poisson-general}
\lambda^{i,N}-\mathcal L_\infty u^{i,N}
=r^{i,N}\quad\text{on }\R^{Nd}
\end{equation}
in the classical sense, and
\begin{equation}\label{eq:erg-corrector-bounds}
\begin{split}
\int &u^{i,N}\,d\mu_\infty^N=0,\qquad|u^{i,N}(\bm x)|\leq C_{N,i}(1+|\bm x|^2),\\
|u^{i,N}(\bm x)-u^{i,N}(\bm y)|
&\leq C_{N,i}(1+|\bm x|+|\bm y|)|\bm x-\bm y|,\quad\text{and}\quad|Du^{i,N}(\bm x)|\leq C_{N,i}(1+|\bm x|),
\end{split}
\end{equation}
for some constant $C_{N,i}\ge 0$. Moreover, if $w\in C^2(\R^{Nd})$ has at most quadratic growth and
$\widetilde\lambda-\mathcal L_\infty w=r^{i,N}$, then
$\widetilde\lambda=\lambda^{i,N}$ and $w-u^{i,N}$ is constant.
\end{proposition}

\begin{proof}
Fix $i\in[N]$ and define
$$u^{i,N}(\bm x):=\int_0^\infty\Bigl(P_sr^{i,N}(\bm x)-\lambda^{i,N}\Bigr)ds,\qquad\bm x\in\R^{Nd}.$$
\textit{Step 1.} (Bounds and normalization).\\
By the first estimate in~\eqref{eq:erg-value-rate}, the integral converges absolutely and $|u^{i,N}(\bm x)|\leq2\kappa^{-1}C_{N,i}(1+|\bm x|^2)$. Since $\mu_\infty^N\in\mathcal P_2(\R^{Nd})$, Fubini's theorem and the invariance of $\mu_\infty^N$ give
\[
\int u^{i,N}\,d\mu_\infty^N=\int_0^\infty\left(\int P_sr^{i,N}\,d\mu_\infty^N-\lambda^{i,N}\right)ds=0.
\]
Let $\bm x,\bm y\in\R^{Nd}$ and couple the canonical states $\bm X^{\infty,\bm x}$ and $\bm X^{\infty,\bm y}$ synchronously. The Lipschitz bound in~\eqref{eq:erg-running-cost-bounds}, the Cauchy--Schwarz inequality, \eqref{eq:erg-canonical-moments} and~\eqref{eq:erg-synchronous-contraction} give
\[
\bigl|P_sr^{i,N}(\bm x)-P_sr^{i,N}(\bm y)\bigr|\leq C_{N,i}(1+|\bm x|+|\bm y|)e^{-\frac{\kappa}{2}s}|\bm x-\bm y|,\qquad s\geq0.
\]
Integrating over $s\in(0,\infty)$ proves the Lipschitz bound in~\eqref{eq:erg-corrector-bounds}. In particular, $u^{i,N}$ is continuous.
\vspace{2mm}

\textit{Step 2.} (The Poisson equation \eqref{eq:poisson-general})\\
Let $t\geq0$. The semigroup property and Fubini's theorem, which applies by the first estimate in~\eqref{eq:erg-value-rate} and~\eqref{eq:erg-canonical-moments}, give
\[
P_tu^{i,N}(\bm x)=\int_0^\infty\Bigl(P_{t+s}r^{i,N}(\bm x)-\lambda^{i,N}\Bigr)ds=u^{i,N}(\bm x)-\int_0^t\Bigl(P_sr^{i,N}(\bm x)-\lambda^{i,N}\Bigr)ds.
\]
By the Markov property of $\bm X^{\infty,\bm x}$, this identity means that the process
\[
u^{i,N}(\bm X_t^{\infty,\bm x})+\int_0^t\bigl(r^{i,N}(\bm X_s^{\infty,\bm x})-\lambda^{i,N}\bigr)ds,\quad t\geq0,
\]
is a martingale. Optional sampling then gives, for every bounded stopping time $\rho$,
\begin{equation}\label{eq:erg-corrector-stopped}
u^{i,N}(\bm x)=\Ex\left[u^{i,N}(\bm X_\rho^{\infty,\bm x})+\int_0^\rho\bigl(r^{i,N}(\bm X_s^{\infty,\bm x})-\lambda^{i,N}\bigr)ds\right].
\end{equation}
Let $B\subset\R^{Nd}$ be a bounded ball. The operator $\mathcal L_\infty$ is uniformly elliptic by \textup{(A3)}, the functions $b$ and $r^{i,N}$ are locally Lipschitz, and $u^{i,N}$ is continuous. As in Step~2 of the proof of Theorem~\ref{thm:stationary-limit}, the Dirichlet problem
$$\mathcal L_\infty w=\lambda^{i,N}-r^{i,N}\ \text{ in }B,\qquad w=u^{i,N}\ \text{ on }\partial B,$$
therefore has a solution $w\in C(\overline B)\cap C^2(B)$. Let $\bm x\in B$ and let $\tau_B$ be the first exit time of $\bm X^{\infty,\bm x}$ from $B$. As in Step~2 of the proof of Theorem~\ref{thm:stationary-limit}, $\Ex[\tau_B]<\infty$. Taking expectations in It\^o's formula for $w(\bm X^{\infty,\bm x})$ up to $\tau_B$ yields
$$w(\bm x)=\Ex\left[u^{i,N}(\bm X_{\tau_B}^{\infty,\bm x})+\int_0^{\tau_B}\bigl(r^{i,N}(\bm X_s^{\infty,\bm x})-\lambda^{i,N}\bigr)ds\right].$$
On the other hand, taking $\rho=t\wedge\tau_B$ in~\eqref{eq:erg-corrector-stopped} and letting $t\to\infty$ \footnote{Dominated convergence applies because $u^{i,N}$ and $r^{i,N}$ are bounded on $\overline B$ and $\Ex[\tau_B]<\infty$.} gives the same expression for $u^{i,N}(\bm x)$. Hence $u^{i,N}=w$ in $B$. Since $B$ is arbitrary, $u^{i,N}\in C^2(\R^{Nd})$ and~\eqref{eq:poisson-general} holds. The gradient bound in~\eqref{eq:erg-corrector-bounds} follows from the Lipschitz bound by letting $\bm y\to\bm x$.
\vspace{2mm}

\textit{Step 3.} (Uniqueness)\\
Let $w\in C^2(\R^{Nd})$ have at most quadratic growth and satisfy $\widetilde\lambda-\mathcal L_\infty w=r^{i,N}$. Then $\mathcal L_\infty(w-u^{i,N})=\widetilde\lambda-\lambda^{i,N}$. We apply It\^o's formula to $(w-u^{i,N})(\bm X^{\infty,\bm x})$ up to the exit time from a ball of radius $R$, take expectations, and let $R\to\infty$. The quadratic growth of $w-u^{i,N}$ and~\eqref{eq:canonical-finite-moments} justify the passage to the limit, and we obtain
$$P_t(w-u^{i,N})(\bm x)=(w-u^{i,N})(\bm x)+(\widetilde\lambda-\lambda^{i,N})\,t,\quad t\geq0.$$
By \eqref{eq:erg-Wasserstein-contraction}, $\delta_{\bm x}P_t\to\mu_\infty^N$ in $W_2$ as $t\to\infty$. Convergence in $W_2$ implies convergence of the integrals of every continuous function with at most quadratic growth, so the left-hand side converges to the finite number $\int(w-u^{i,N})\,d\mu_\infty^N$. Therefore $\widetilde\lambda=\lambda^{i,N}$, and then $(w-u^{i,N})(\bm x)=\int(w-u^{i,N})\,d\mu_\infty^N$, for every $\bm x\in\R^{Nd}$.
\end{proof}

\begin{proposition}\label{prop:erg-forward-BSDE}
Assume \textup{(A1)-(A3)} and, for $s\ge 0$, let $\mu_s$ be the law of the canonical state $\bm X^\infty_s$ started with law
$\mu_0\in\mathcal P_2(\R^{Nd})$. Then, $(\mu_s)_{s\ge 0}$ satisfies the Fokker--Planck equation $\partial_s\mu_s=\mathcal L_\infty^*\mu_s$ in the weak sense, where $\mathcal{L}_\infty^*$ is the adjoint of the generator $\mathcal{L}_\infty$ from \eqref{eq:erg-generator}. In particular, $\mathcal L_\infty^*\mu_\infty^N=0$, in the weak sense for the invariant measure $\mu^N_\infty$.
\end{proposition}

\begin{proof}
The proof is standard. Let $\varphi\in C_c^\infty(\R^{Nd})$ and $s\geq0$. We will show that
\begin{equation}\label{eq:erg-weak-FP}
\int\varphi\,d\mu_s=\int\varphi\,d\mu_0+\int_0^s\int\mathcal L_\infty\varphi\,d\mu_u\,du.
\end{equation}
Let $\bm X^\infty$ be the canonical state with $\mathcal L(\bm X_0^\infty)=\mu_0$. It\^o's formula and~\eqref{eq:erg-generator} give
\[
\varphi(\bm X_s^\infty)=\varphi(\bm X_0^\infty)+\int_0^s\mathcal L_\infty\varphi(\bm X_u^\infty)\,du+\sqrt2\int_0^sD\varphi(\bm X_u^\infty)\cdot\Sigma\,d\bm W_u.
\]
Since $\varphi$ has compact support and $b$ is continuous, $D\varphi$ and $\mathcal L_\infty\varphi$ are bounded. The stochastic integral is therefore a martingale, and taking expectations proves~\eqref{eq:erg-weak-FP}. If $\mu_0=\mu_\infty^N$, then $\mu_s=\mu_\infty^N$ for every $s\geq0$ by Proposition~\ref{prop:erg-contraction}, and~\eqref{eq:erg-weak-FP} reduces to $\int\mathcal L_\infty\varphi\,d\mu_\infty^N=0$, which is the weak form of $\mathcal L_\infty^*\mu_\infty^N=0$.
\end{proof}

\begin{remark}
\label{rem:erg-corrector-adjoint}
Since $-\mathcal L_\infty v_\infty^i(\bm x)=D_xL^i(x^i,b^i(\bm x))+D_{x^i}F^i(\bm x)$ and $-\mathcal L_\infty u^{i,N}(\bm x)=L^i(x^i,b^i(\bm x))+F^i(\bm x)-\lambda^{i,N}$, it is tempting to compare $D_{x^i}u^{i,N}$ and $v_\infty^i$, $i\in [N]$. In general, we have $D_{x^i}u^{i,N}\ne v_\infty^i$, and therefore, unlike in the closed-loop $N$-player game case, $D_{x^i}u^{i,N}$ does not create a stationary Nash equilibrium via the controls
$-D_pH^i(x^i,D_{x^i}u^{i,N})$, $i\in [N]$.\\[0.3ex]
This distinction already occurs for quadratic data satisfying
\textup{(A1)--(A3)}. Take $N=2$, $d=1$, $\sigma^1=\sigma^2=1$,
$G^1=G^2=0$, and
\[
L^i(x,a)=\tfrac12(x^2+a^2),
\qquad F^i(\bm x)=\tfrac12(x^1+x^2)^2,\;\;\;i\in \{1,2\}.
\]
The terminal gradient map is zero and the running gradient map has
the positive semidefinite matrix $\bigl(\begin{smallmatrix}1&1\\1&1\end{smallmatrix}\bigr)$.
Let
\[
P=\begin{pmatrix}p_0&p_1\\p_1&p_0\end{pmatrix},
\qquad p_0=\frac{\sqrt3+1}{2},\quad
p_1=\frac{\sqrt3-1}{2},\qquad
P^2=\begin{pmatrix}2&1\\1&2\end{pmatrix}.
\]
The stable stationary field is $v_\infty(\bm x)=P\bm x$ and the drift is
$b(\bm x)=-P\bm x$. Player one's running cost is
$r^{1,2}(\bm x)=\tfrac12\bm x^\top R\bm x$, where
\[
R=\begin{pmatrix}p_0^2+2&p_0p_1+1\\p_0p_1+1&p_1^2+1\end{pmatrix}.
\]
Its corrector has the form
$u^{1,2}(\bm x)=\tfrac12\bm x^\top Q\bm x+\mathrm{constant}$,
with $Q$ symmetric and $PQ+QP=R$.
Indeed, substitution in the Poisson equation gives this matrix
identity and $\lambda^{1,2}=\operatorname{tr}Q$; centering fixes the additive constant. If $D_{x^1}u^{1,2}=v_\infty^1$, then $Q_{11}=p_0$ and $Q_{12}=p_1$. The $(1,1)$ entry of the matrix identity would give $2(p_0^2+p_1^2)=p_0^2+2.$ Its left-hand side is $4$ and its right-hand side is $3+\sqrt3/2$, a contradiction.
\end{remark}

\subsection{Long-time behavior of equilibrium distributions}
The following pointwise estimate is needed because its constant is
independent of both the observation time $s\ge 0$ and the horizon $T$.

\begin{lemma}\label{lem:erg-pointwise-turnpike}
Assume \textup{(A1)-(A3)}. Start the finite-horizon and canonical infinite-horizon systems $\bm X^T,\bm X^\infty$ from
the same deterministic $\bm x$ at time $0$, with the same Brownian
motions. There is a constant $C$ such that, for all $T>0$ and $s\in [0,T]$,
\begin{equation}\label{eq:erg-pointwise-turnpike}
\Ex\left[
 |\bm X_s^T-\bm X_s^\infty|^2
 +|\bm Y_s^T-\bm Y_s^\infty|^2
 +|\bm\alpha_T(s)-\bm\alpha_\infty(s)|^2
\right]
\leq C(1+|\bm x|^2)e^{-2c(T-s)},
\end{equation}
where $c$ is the rate from the turnpike estimate in Theorem \ref{thm:turnpike}.
\end{lemma}

\begin{proof}
By Proposition~\ref{N-PDE-classical}, $\bm\alpha_T(s)=b_T(s,\bm X_s^T)$, with $b_T$ as in~\eqref{eq:erg-canonical-drifts}. We compare $\bm X^T$ with $\bm X^\infty$ through the drift error
\[
\bm e_s:=b_T(s,\bm X_s^T)-b(\bm X_s^T),\qquad0\leq s\leq T.
\]
Since $v_T$ is continuous on $[0,T]\times\R^{Nd}$, the estimate~\eqref{goal:estimate2} also holds at $t=T$. Together with~\eqref{eq:H-gradient-Lipschitz} and~\eqref{unifL2}, it gives
\begin{equation}\label{eq:erg-drift-error}
\Ex[|\bm e_s|^2]\leq\frac{1}{C_L^2}\Ex\bigl[|v_T(s,\bm X_s^T)-v_\infty(\bm X_s^T)|^2\bigr]\leq C(1+|\bm x|^2)e^{-2c(T-s)},\qquad0\leq s\leq T.
\end{equation}
Here and below, $C$ does not depend on $s$, $T$ or $\bm x$.

Write $\Delta\bm X_s:=\bm X_s^T-\bm X_s^\infty$, and let $\Delta v_s$ be the difference of $v_\infty$ evaluated at $\bm X_s^T$ and at $\bm X_s^\infty$. Set
\[
f(s):=\Ex[|\Delta\bm X_s|^2]\qquad\text{and}\qquad q(s):=\Ex[\Delta\bm X_s\cdot\Delta v_s],\qquad0\leq s\leq T.
\]
The states $\bm X^T$ and $\bm X^\infty$, with $\bm Y=v_\infty(\bm X)$ and $Z^{i,j}=D_{x^j}v^i_\infty(\bm X)\sigma^j$ for $\bm X\in\{\bm X^T,\bm X^\infty\}$, satisfy the equations of Lemma~\ref{lem:pairing} on $[0,T]$, with $\bm\phi^1=\bm e$, $\bm\psi^1=D_{\bm x}v_\infty(\bm X^T)\bm e$ and $\bm\phi^2=\bm\psi^2=0$. For $\bm X^T$, this follows from $d\bm X_s^T=\bigl(b(\bm X_s^T)+\bm e_s\bigr)\,ds+\sqrt2\Sigma\,d\bm W_s$, It\^o's formula and~\eqref{N-PDEinf}. The integrability conditions hold because $v_\infty$ is Lipschitz, $\bm e$ is square integrable by~\eqref{eq:erg-drift-error}, and both states have square-integrable suprema on $[0,T]$. The function $e$ of Lemma~\ref{lem:pairing} is at least $f$, so Lemma~\ref{lem:pairing}(i) gives, for almost every $s\in[0,T]$,
\begin{equation*}
q'(s)\leq-C_Lf(s)+\Ex\bigl[\bm e_s\cdot\Delta v_s+\Delta\bm X_s\cdot D_{\bm x}v_\infty(\bm X_s^T)\bm e_s\bigr]\,.
\end{equation*} Since $v_\infty$ is $K$-Lipschitz, the expectation on the right-hand side is at most $2K\Ex[|\Delta\bm X_s||\bm e_s|]$. Young's inequality then gives
\begin{equation}\label{eq:erg-perturbed-pairing}
q'(s)\leq-\frac{C_L}{2}f(s)+\frac{2K^2}{C_L}\Ex[|\bm e_s|^2].
\end{equation}
By~\eqref{eq:erg-strong-monotonicity} and the Lipschitz bound, $C_Lf(s)\leq q(s)\leq Kf(s)$. With $\theta=C_L/(2K)$, the preceding inequality becomes $q'(s)\leq-\theta q(s)+\frac{2K^2}{C_L}\Ex[|\bm e_s|^2]$. Since $q(0)=0$, the bound~\eqref{eq:erg-drift-error} gives
\[
C_Lf(s)\leq q(s)\leq\frac{2K^2}{C_L}\int_0^se^{-\theta(s-u)}\Ex[|\bm e_u|^2]\,du\leq C(1+|\bm x|^2)e^{-2c(T-s)}\int_0^se^{-(\theta+2c)(s-u)}\,du,
\]
and the last integral is at most $(\theta+2c)^{-1}$. This is the bound for the states in~\eqref{eq:erg-pointwise-turnpike}.

For the adjoints and the controls, the identities $\bm Y_s^T=v_T(s,\bm X_s^T)$ and $\bm Y_s^\infty=v_\infty(\bm X_s^\infty)$, together with~\eqref{eq:H-gradient-Lipschitz}, give
\begin{align*}
|\bm Y_s^T-\bm Y_s^\infty|&\leq|v_T(s,\bm X_s^T)-v_\infty(\bm X_s^T)|+K|\Delta\bm X_s|,\;\;
|\bm\alpha_T(s)-\bm\alpha_\infty(s)|\leq\frac{K_L}{C_L}|\Delta\bm X_s|+\frac{1}{C_L}|\bm Y_s^T-\bm Y_s^\infty|.
\end{align*}
Squaring, taking expectations, and using the bound for the states (bound for $f$) and the second inequality in~\eqref{eq:erg-drift-error} proves~\eqref{eq:erg-pointwise-turnpike}.
\end{proof}

\begin{remark}[Dependence on the number of players]
Under \textup{(U)}, $C_L$ can be bounded below and $K$ above by
positive constants independent of $N$. Hence the monotonicity,
contraction, and averaged-residual estimates admit constants
independent of $N$. Combining the same pointwise turnpike proof with
\eqref{goal:estimate2-uniform} and the population-uniform moment
estimate gives
\[
\Ex\left[
 |\bm X_s^T-\bm X_s^\infty|^2
 +|\bm Y_s^T-\bm Y_s^\infty|^2
 +|\bm\alpha_T(s)-\bm\alpha_\infty(s)|^2
\right]
\leq C_*(N+|\bm x|^2)e^{-2c_*(T-s)}.
\]
\end{remark}

\subsubsection*{Proof of Theorem \ref{thm:mainLTMdistr}}\phantomsection\label{pf:mainLTMdistr} It follows from the corollary below.

\begin{corollary}\label{cor:erg-long-time-limit}
Consider the setup of Theorem \ref{thm:mainLTMdistr}. For a deterministic initial configuration $\bm x\in (\R^d)^N$, there is a constant $C$ independent of $T$ and $s\in [0,T]$ such that
\begin{equation}\label{eq:erg-state-law-turnpike}
W_2\!\left(\mathcal L(\bm X_s^{T,\bm x}),\mu_\infty^N\right)
\leq C(1+|\bm x|)
\left(e^{-\theta s}+e^{-c(T-s)}\right),
\qquad 0\leq s\leq T,
\end{equation}
where $\theta=C_L/(2K)$.
\end{corollary}

\begin{proof}
Start the horizon-$T$ equilibrium $\bm X^T$ and the canonical process $\bm X^\infty$ from $\bm x$ at time zero, driven by the same Brownian motions, as in Lemma~\ref{lem:erg-pointwise-turnpike}. The triangle inequality for $W_2$ gives
\[
W_2\!\left(\mathcal L(\bm X_s^{T,\bm x}),\mu_\infty^N\right)
\leq W_2\!\left(\mathcal L(\bm X_s^{T,\bm x}),\mathcal L(\bm X_s^{\infty,\bm x})\right)
+W_2\!\left(\mathcal L(\bm X_s^{\infty,\bm x}),\mu_\infty^N\right).
\]
The synchronous pair $(\bm X_s^{T,\bm x},\bm X_s^{\infty,\bm x})$ is a coupling of the two laws in the first term. By \eqref{eq:erg-pointwise-turnpike} and $\sqrt{1+|\bm x|^2}\leq1+|\bm x|$, the first term is therefore at most
\[
\left(\Ex\left[|\bm X_s^{T,\bm x}-\bm X_s^{\infty,\bm x}|^2\right]\right)^{1/2}\leq C(1+|\bm x|)e^{-c(T-s)}.
\]
For the second term, $\mathcal L(\bm X_s^{\infty,\bm x})=\delta_{\bm x}P_s$, so \eqref{eq:erg-Wasserstein-contraction} with $\nu=\delta_{\bm x}$ gives the bound $\kappa^{-1/2}e^{-\kappa s/2}W_2(\delta_{\bm x},\mu_\infty^N)$. The only coupling of $\delta_{\bm x}$ and $\mu_\infty^N$ is their product, hence
\[
W_2(\delta_{\bm x},\mu_\infty^N)^2=\int|\bm x-\bm y|^2\,\mu_\infty^N(d\bm y)\leq2|\bm x|^2+2\int|\bm y|^2\,\mu_\infty^N(d\bm y),
\]
and $W_2(\delta_{\bm x},\mu_\infty^N)\leq C(1+|\bm x|)$ because $\mu_\infty^N\in\mathcal P_2(\R^{Nd})$. Since $\kappa/2=C_L/(2K)=\theta$, the second term is at most $C(1+|\bm x|)e^{-\theta s}$. Adding the two bounds proves \eqref{eq:erg-state-law-turnpike}.

Under \textup{(U)}, the canonical moment estimate obtained in the
proof of Theorem~\ref{thm:stationary-limit}, started from $\bm0$,
gives
\[
\sup_{s\geq0}
\Ex[|\bm X_s^{\infty,\bm0}|^2]\leq C_*N.
\]
Since $\delta_{\bm0}P_s\to\mu_\infty^N$ in $W_2$, it follows that
\[
\int|\bm y|^2\,\mu_\infty^N(d\bm y)\leq C_*N.
\]
Consequently,
\[
W_2(\delta_{\bm x},\mu_\infty^N)
\leq C_*(\sqrt N+|\bm x|).
\]
Using the lower monotonicity bound $\overline C_L$ and the
Lipschitz bound $K_*$ in the contraction argument gives
\[
W_2(\delta_{\bm x}P_s,\mu_\infty^N)
\leq C_*(\sqrt N+|\bm x|)e^{-\theta_*s},
\qquad
\theta_*=\frac{\overline C_L}{2K_*}.
\]
The population-uniform pointwise turnpike estimate in the preceding
remark also gives
\[
W_2\!\left(
\mathcal L(\bm X_s^{T,\bm x}),
\mathcal L(\bm X_s^{\infty,\bm x})
\right)
\leq C_*(\sqrt N+|\bm x|)e^{-c_*(T-s)}.
\]
The triangle inequality proves
\eqref{eq:state-law-turnpike-uniform} for $t_0=0$.
The initial-time translation gives the general case.
\end{proof}

\vspace{3mm}

\subsection{Long-time behavior of the value functions}
Set $W_\tau^{i,N}(\bm x):=V_\tau^{i,N}(0,\bm x),\; \tau>0$, and $W_0^{i,N}:=G^i$. Time homogeneity and uniqueness of the finite-horizon Pontryagin system
give, for $0\leq t\leq T$,
\begin{equation}\label{eq:erg-time-consistency}
V_T^{i,N}(t,\bm x)=W_{T-t}^{i,N}(\bm x).
\end{equation}
We start by showing the following (weighted) Lipschitz bound.

\begin{lemma}
\label{lem:W-lip-general}
Suppose that \textup{(A1)-(A3)} hold. For every $i\in[N]$, there is a constant $C_W\ge 0$ such that, for every $\tau\geq0$,
\begin{equation}\label{eq:erg-W-weighted-Lipschitz}
|W_\tau^{i,N}(\bm y)-W_\tau^{i,N}(\bm y')|
\leq C_W(1+|\bm y|+|\bm y'|)|\bm y-\bm y'|.
\end{equation}
\end{lemma}

\begin{proof}
Drive the two horizon-$\tau$ equilibrium systems~\eqref{N-FBSDE} synchronously, starting from $\bm y,\bm y'\in(\R^d)^N$ at time $0$, denote them by $(\bm X,\bm Y)$ and $(\bm X',\bm Y')$, and put $f(s)=\Ex[|\bm X_s-\bm X_s'|^2]$. Lemma~\ref{lem:pairing}(i), integrated over $[s,\tau]$, the monotonicity of $D\bm G$ at $\tau$ and the $K$-Lipschitz bound~\eqref{lipbound} for $v_\tau$ give
\[
C_L\int_s^\tau f(r)\,dr
\leq\Ex[(\bm X_s-\bm X_s')\cdot(\bm Y_s-\bm Y_s')]
\leq Kf(s),\qquad0\leq s\leq\tau.
\]
Set $g(s)=\int_s^\tau f(r)\,dr$ and recall $\kappa=C_L/K$. Then $g'(s)+\kappa g(s)\leq0$ and $g(0)\leq\kappa^{-1}f(0)$, so $g(s)\leq\kappa^{-1}e^{-\kappa s}|\bm y-\bm y'|^2$. By~\eqref{eq:H-gradient-Lipschitz} and~\eqref{lipbound}, the drift $b_\tau(s,\cdot)$ in~\eqref{eq:erg-canonical-drifts} is Lipschitz with constant $\Lambda=(K_L+K)/C_L$, uniformly in $\tau$ and $s$. Since the synchronous noise cancels, $|f'(s)|\leq2\Lambda f(s)$, and therefore $f(s)\leq e^{2\Lambda}f(r)$ whenever $|s-r|\leq1$ and $r\in [0,s]$. If $s\geq1$, averaging this inequality over $r\in[s-1,s]$ gives $f(s)\leq e^{2\Lambda}g(s-1)$. If $s<1$, then $f(s)\leq e^{2\Lambda}f(0)$. Thus, in both cases,
\[
f(s)\leq Ce^{-\kappa s}|\bm y-\bm y'|^2,
\qquad0\leq s\leq\tau,
\]
where $C$ depends only on $\Lambda$ and $\kappa$.

The two equilibrium controls are $b_\tau(s,\bm X_s)$ and $b_\tau(s,\bm X_s')$, so their difference is bounded in $L^2$ by $\Lambda\sqrt{f(s)}$. Using~\eqref{erg:growth}, Cauchy--Schwarz, and the moment bound~\eqref{unifL2} for the states and controls of the two systems, the difference of the running-cost expectations at time $s$ is at most
\[
C(1+|\bm y|+|\bm y'|)\sqrt{f(s)}\leq C(1+|\bm y|+|\bm y'|)e^{-\kappa s/2}|\bm y-\bm y'|.
\]
Integrate this bound on $[0,\tau]$. The Hessian of $G^i$ is bounded, so $G^i$ satisfies the same bound as $F^i$ in~\eqref{erg:growth}, and the same estimate at $s=\tau$ applies to the terminal-cost difference. The resulting constant is independent of $\tau$. For $\tau=0$, \eqref{eq:erg-W-weighted-Lipschitz} holds because $W_0^{i,N}=G^i$ has a bounded Hessian. The proof is complete.
\end{proof}

\subsubsection*{Proof of Theorem \ref{thm:mainLTB}}\phantomsection\label{pf:mainLTB}
The following theorem implies Theorem \ref{thm:mainLTB}. Indeed, let $u^{i,N}$ be the corrector of Proposition~\ref{prop:erg-poisson}, let $\kappa^{i,N}_G$ be the limit~\eqref{eq:erg-terminal-normalization}, and set $w^{i,N}_G:=u^{i,N}+\kappa^{i,N}_G$. By Proposition~\ref{prop:erg-poisson}, $w^{i,N}_G$ is the unique classical solution with at most quadratic growth of the equation in Theorem~\ref{thm:mainLTB} whose integral with respect to $\mu^N_\infty$ is $\kappa^{i,N}_G$, and~\eqref{eq:erg-value-expansion} is the estimate of Theorem~\ref{thm:mainLTB}.

\begin{theorem}
\label{thm:expansion-general}
Assume \textup{(A1)-(A3)} and fix $(t,\bm x)\in [0,T]\times (\R^d)^N$. For every $i\in[N]$, the limit
\begin{equation}\label{eq:erg-terminal-normalization}
\kappa_G^{i,N}:=
\lim_{\tau\to\infty}
\left(\int W_\tau^{i,N}(\bm y)\,d\mu_\infty^N(\bm y)
 -\tau\lambda^{i,N}\right)
\end{equation}
exists.
There are $C_E,c_E>0$, independent of $t,T,\bm x,i$, such that
\begin{equation}\label{eq:erg-value-expansion}
\left|V_T^{i,N}(t,\bm x)-(T-t)\lambda^{i,N}
 -u^{i,N}(\bm x)-\kappa_G^{i,N}\right|
\leq C_E(1+|\bm x|^2)e^{-c_E(T-t)},
\end{equation}
where $u^{i,N}$ is the corrector of Proposition \ref{prop:erg-poisson}.
\end{theorem}

\begin{proof}
Fix $i\in[N]$. For notational ease we write $W_\tau$, $\lambda$, $u$, $r$ and $\kappa_G$ for $W_\tau^{i,N}$, $\lambda^{i,N}$, $u^{i,N}$, $r^{i,N}$ and $\kappa_G^{i,N}$. By Lemma~\ref{lem:W-lip-general}, $W_\tau$ has at most quadratic growth, so the quantity
$$\kappa_G(\tau):=\int W_\tau\,d\mu_\infty^N-\tau\lambda,\quad\tau\geq0,$$
is finite. The first assertion of the theorem is that $\kappa_G(\tau)$ converges as $\tau\to\infty$.
Throughout the proof, $C$ denotes a constant independent of $\tau$, $S$ and $\bm x$. Also, $\theta=C_L/(2K)=\kappa/2$, where $\kappa$ appeared in \eqref{eq:erg-semigroup-cost-rate}.
\vspace{2mm}

\textit{Step 1.} (Splitting the horizon, dynamic programming)\\
Let $S\in [0,\tau]$ and $\bm x\in\R^{Nd}$. Let $\bm X^\tau$ be the equilibrium state of the horizon-$\tau$ game started from $\bm x$ at time $0$, with control $\bm\alpha_\tau$, and let $\bm X^\infty$ be the canonical state started from $\bm x$, with the same Brownian motions. We claim that
\begin{equation}\label{eq:erg-split}
W_\tau(\bm x)=\Ex\left[\int_0^S\Bigl(L^i\bigl(X_s^{i,\tau},\alpha_\tau^i(s)\bigr)+F^i(\bm X_s^\tau)\Bigr)ds\right]+\Ex\bigl[W_{\tau-S}(\bm X_S^\tau)\bigr].
\end{equation}
By Proposition \ref{N-PDE-classical}, $\bm X^\tau$ solves the stochastic differential equation with drift $b_\tau(s,\cdot)$, which is Lipschitz uniformly in $s$, and $\bm\alpha_\tau(s)=b_\tau(s,\bm X_s^\tau)$. As we argued for \eqref{eq:erg-time-consistency}, time homogeneity and the uniqueness part of Proposition~\ref{N-PDE-classical} give $v_\tau(S+q,\cdot)=v_{\tau-S}(q,\cdot)$, hence $b_\tau(S+q,\cdot)=b_{\tau-S}(q,\cdot)$, for $q\in[0,\tau-S]$. By the Markov property of this equation, the conditional law of $(\bm X_{S+q}^\tau)_{q\in[0,\tau-S]}$ given $\mathcal F_S$ is the law of the equilibrium state of the horizon-$(\tau-S)$ game started from $\bm X_S^\tau$. The conditional expectation, given $\mathcal F_S$, of the cost on $[S,\tau]$, including the terminal cost, is therefore $W_{\tau-S}(\bm X_S^\tau)$. This proves~\eqref{eq:erg-split}.
\vspace{2mm}

\textit{Step 2.} (Comparison with the canonical state $\bm X^\infty=\bm X^{\infty,\bm x}$)\\
We estimate the two terms in~\eqref{eq:erg-split}. Since $\alpha_\infty^i(s)=b^i(\bm X_s^\infty)$, we have $r(\bm X_s^\infty)=L^i(X_s^{i,\infty},\alpha_\infty^i(s))+F^i(\bm X_s^\infty)$. The bounds~\eqref{erg:growth}, the Cauchy--Schwarz inequality, the moment bounds~\eqref{unifL2} and~\eqref{eq:erg-canonical-moments}, the linear growth of $b$, and the turnpike estimate in Lemma \ref{lem:erg-pointwise-turnpike} with $T=\tau$ give
$$\Bigl|\Ex\bigl[L^i\bigl(X_s^{i,\tau},\alpha_\tau^i(s)\bigr)+F^i(\bm X_s^\tau)\bigr]-\Ex\bigl[r(\bm X_s^\infty)\bigr]\Bigr|\leq C(1+|\bm x|^2)e^{-c(\tau-s)},\quad s\in [0,\tau].$$
The integral of the right-hand side over $s\in[0,S]$ is at most $c^{-1}C(1+|\bm x|^2)e^{-c(\tau-S)}$. Next, It\^o's formula for $u(\bm X^\infty)$, the Poisson equation~\eqref{eq:poisson-general} and the gradient bound in~\eqref{eq:erg-corrector-bounds} give
\[
\int_0^S\Ex\bigl[r(\bm X_s^\infty)\bigr]\,ds=S\lambda+u(\bm x)-P_Su(\bm x).
\]
By~\eqref{eq:erg-corrector-bounds}, we may take $\varphi=u$ in~\eqref{eq:erg-semigroup-cost-rate} and, since $\int u\,d\mu_\infty^N=0$, we obtain $|P_Su(\bm x)|\leq C(1+|\bm x|^2)e^{-\theta S}$. For the second term in~\eqref{eq:erg-split}, Lemma~\ref{lem:W-lip-general}, the Cauchy--Schwarz inequality, the same moment bounds and Lemma~\ref{lem:erg-pointwise-turnpike} give
$$\bigl|\Ex\bigl[W_{\tau-S}(\bm X_S^\tau)\bigr]-P_SW_{\tau-S}(\bm x)\bigr|\leq C_W\Ex\bigl[(1+|\bm X_S^\tau|+|\bm X_S^\infty|)|\bm X_S^\tau-\bm X_S^\infty|\bigr]\leq C(1+|\bm x|^2)e^{-c(\tau-S)}.$$
On the other hand, \eqref{eq:erg-semigroup-cost-rate} with $\varphi=W_{\tau-S}$ and $C_\varphi=C_W$ gives
\[
\left|P_SW_{\tau-S}(\bm x)-\int W_{\tau-S}\,d\mu_\infty^N\right|\leq C(1+|\bm x|^2)e^{-\theta S}.
\]
Since $\int W_{\tau-S}\,d\mu_\infty^N=\kappa_G(\tau-S)+(\tau-S)\lambda$, inserting these four estimates in~\eqref{eq:erg-split} yields a constant $C_0$, independent of $\tau$, $S$ and $\bm x$, such that
\begin{equation}\label{eq:erg-one-step}
\bigl|W_\tau(\bm x)-\tau\lambda-u(\bm x)-\kappa_G(\tau-S)\bigr|\leq C_0(1+|\bm x|^2)\bigl(e^{-c(\tau-S)}+e^{-\theta S}\bigr),\quad S\in [0,\tau].
\end{equation}
\textit{Step 3.} (Existence of the limit)\\
We integrate~\eqref{eq:erg-one-step} with respect to $\mu_\infty^N(d\bm x)$ and use $\int u\,d\mu_\infty^N=0$. With $C_1:=C_0\bigl(1+\int|\bm y|^2\,\mu_\infty^N(d\bm y)\bigr)$, this gives
\[
|\kappa_G(\tau)-\kappa_G(\tau-S)|\leq C_1\bigl(e^{-c(\tau-S)}+e^{-\theta S}\bigr),\quad S\in [0,\tau].
\]
Set $c_E:=\frac12\min\{c,\theta\}$. For $h>0$ and $\tau\geq2h$, the choice $S=\tau-h\geq h$ gives
\[
|\kappa_G(\tau)-\kappa_G(h)|\leq2C_1e^{-2c_Eh}.
\]
Hence $|\kappa_G(\tau)-\kappa_G(\tau')|\leq4C_1e^{-2c_Eh}$ for all $\tau,\tau'\geq2h$. This means that $\tau\mapsto \kappa_G(\tau)$ is Cauchy and hence it converges as $\tau\to\infty$. Its limit is $\kappa_G$ in~\eqref{eq:erg-terminal-normalization}, and letting $\tau\to\infty$ in the last inequality yields
$$|\kappa_G-\kappa_G(h)|\leq2C_1e^{-2c_Eh},\quad h>0.$$
\textit{Step 4.} (Conclusion)\\
Let $\tau>0$. Taking $S=\tau/2$ in~\eqref{eq:erg-one-step} and $h=\tau/2$ in the last bound of Step~3, we obtain
\[
\bigl|W_\tau(\bm x)-\tau\lambda-u(\bm x)-\kappa_G\bigr|\leq2C_0(1+|\bm x|^2)e^{-c_E\tau}+2C_1e^{-c_E\tau}\leq C_E(1+|\bm x|^2)e^{-c_E\tau},
\]
with $C_E:=2(C_0+C_1)$. By~\eqref{eq:erg-time-consistency}, this is~\eqref{eq:erg-value-expansion} with $\tau=T-t$, for $t<T$. For $t=T$, the estimate follows from the quadratic growth of $G^i$ and $u$ after enlarging $C_E$. Taking the largest $C_E$ over $i\in[N]$ makes the constants independent of $i$.
\end{proof}

\begin{remark}
The terminal cost affects the additive constant $\kappa_G^{i,N}$ and not $\lambda^{i,N}, u^{i,N}$.
\end{remark}

\appendix

\section{Technical proofs}\label{Appendix2}

\begin{lemma}\label{rem:erg-B-admissible}
Let $\bm x\in(\R^d)^N$. For $i\in [N]$, $\beta\in\mathcal{A}^i_\infty(\bm x)$ and $R\geq1$, we have
$$\Ex[|X_R^{i,\beta}|^2]\leq3R\mathcal E^i(\beta)+6\|\sigma^i\|_{\mathrm{HS}}^2.$$
In particular,
$\Ex[|X_R^{i,\beta}|^2]/R^2\to0$ as $R\to\infty$.
Moreover, $\mathcal{A}^i_\infty(\bm x)$ is closed under affine combinations
of its controls with deterministic coefficients whose sum is one.
\end{lemma}

\begin{proof}
For $s\in[R-1,R]$, the state equation and Cauchy--Schwarz give
\[
\Ex[|X_R^{i,\beta}|^2]
\leq3\Ex[|X_s^{i,\beta}|^2]
+3\int_{R-1}^R\Ex[|\beta_u|^2]du
+6\|\sigma^i\|_{\mathrm{HS}}^2.
\]
Integrate over $s\in[R-1,R]$ and bound the resulting state and
control integrals by $R\mathcal E^i(\beta)$. For an affine
combination, the corresponding state is the same affine combination
of the states, since all controls have the same initial state and
Brownian motion. The assertions then follow from elementary squared
norm bounds.
\end{proof}

\begin{lemma}[First-order condition]\label{lem:erg-first-order}
Suppose that Assumptions \textup{(A1)-(A3)} hold. Let $\bm x\in (\R^d)^N$ and $\bm\beta\in\mathcal{A}_\infty(\bm x)$ be a strong average-cost
equilibrium. For $i\in[N]$ and $\gamma\in\mathcal{A}^i_\infty(\bm x)$, set
$h:=\gamma-\beta^i$ and
$\delta X_s:=X_s^{i,\gamma}-X_s^{i,\beta^i}=\int_0^sh_u\,du$.
Then
\begin{equation}\label{erg:FOC}
\lim_{R\to\infty}\frac1R\int_0^R
\Ex\left[
\bigl(D_xL^i(X_s^{i,\bm\beta},\beta_s^i)
+D_{x^i}F^i(\bm X_s^{\bm\beta})\bigr)\cdot\delta X_s
+D_aL^i(X_s^{i,\bm\beta},\beta_s^i)\cdot h_s
\right]ds=0.
\end{equation}
\end{lemma}

\begin{proof}
Let $D_R$ denote the normalized integral inside the limit in \eqref{erg:FOC}. Its
absolute value is uniformly bounded for $R\geq1$, by the linear
growth of the gradients and the Cauchy--Schwarz inequality. For every $\varepsilon\in\R$, the control
$\beta^{i,\varepsilon}=\beta^i+\varepsilon h$ belongs to
$\mathcal{A}_\infty^i(\bm x)$, and its implied state is
$X^{i,\bm\beta}+\varepsilon\delta X$. The bounded Hessians in
\textup{(A1)} imply, by Taylor's formula for the costs,
\[
\left|
A_R^i([\bm\beta^{-i};\beta^{i,\varepsilon}])
-A_R^i(\bm\beta)-\varepsilon D_R
\right|
\leq \frac{C\varepsilon^2}{R}\Ex\left[\int_{0}^R(|\delta X_s|^2+|h_s|^2)ds\right],
\]
for some constant $C>0$ depending on (A1). For $R\geq1$, the right-hand side is bounded from above by
$C\e^2Q:=2C\e^2\bigl(\mathcal E^i(\gamma)+
\mathcal E^i(\beta^i)\bigr)<\infty$. Using \eqref{erg:strong}, we then obtain
$\liminf_R(\varepsilon D_R)\geq-C\varepsilon^2Q$.
For $\varepsilon>0$, this implies
$\liminf_RD_R\geq-C\varepsilon Q$; for $\varepsilon<0$ it implies
$\limsup_RD_R\leq C|\varepsilon|Q$. Letting
$\varepsilon\to0$ from the two sides proves the claim.
\end{proof}

\section{Characterizations of strong average-cost equilibria}\label{char}

In this section of the appendix we present some characterizations of strong average-cost equilibria (see Definition \ref{erg-strong}).
We distinguish uniqueness of the infinite-horizon Pontryagin
system from uniqueness for the average-cost criterion.

\begin{proposition}[Uniqueness of the stable Pontryagin system]
\label{prop:erg-stable-Pontryagin}
For a prescribed deterministic initial configuration $\bm x\in (\R^d)^N$, the
infinite-horizon system
\begin{equation}\label{eq:erg-Pontryagin}
\begin{cases}
dX_s^i=-D_pH^i(X_s^i,Y_s^i)\,ds+\sqrt2\,\sigma^i\,dW_s^i,\\
dY_s^i=\bigl(D_xH^i(X_s^i,Y_s^i)-D_{x^i}F^i(\bm X_s)\bigr)\,ds +\sqrt2\displaystyle\sum_{j=1}^NZ_s^{i,j}\,dW_s^j,\\
\bm X_0=\bm x
\end{cases}
\end{equation}
has exactly one adapted solution satisfying
\begin{equation}\label{eq:erg-stable-class}
\sup_{s\geq0}\Ex\bigl[|\bm X_s|^2+|\bm Y_s|^2\bigr]<\infty,
\qquad
\Ex\left[\int_0^R\|\bm Z_s\|_{\mathrm{HS}}^2\,ds\right]<\infty
\quad\text{for every }R<\infty.
\end{equation}
It is the solution supplied by Theorem~\ref{thm:stationary-limit},
and its control is
\begin{equation}\label{eq:erg-optimal-control}
\alpha_\infty^{i}(s)
=-D_pH^i(X_s^{i,\infty},Y_s^{i,\infty}),
\qquad Y_s^{i,\infty}=v_\infty^i(\bm X_s^\infty).
\end{equation}
Uniqueness of $\bm X,\bm Y$ is up to indistinguishability; uniqueness
of $\bm\alpha,\bm Z$ is up to $ds\otimes d\Pr$-null sets.
\end{proposition}

\begin{proof}
Existence follows from Theorem~\ref{thm:stationary-limit}, the uniform
moment estimate, and the boundedness of $D_{\bm x}v_\infty$. For uniqueness,
let two solutions be driven by the same Brownian motions, and write
their differences as $\Delta\bm X,\Delta\bm Y,\Delta\bm Z$ and
$\Delta\bm a$, where
$a^i=-D_pH^i(X^i,Y^i)$.
 On every interval $[0,R]$, the two solutions satisfy the equations of Lemma~\ref{lem:pairing} with $\bm\phi^k=\bm\psi^k=0$, and $\Delta\bm a$ is the process $\Delta\bm\alpha$ of the lemma; the integrability conditions hold because~\eqref{eq:erg-stable-class} and the linear growth of $D_pH^i$ give $\Ex[\sup_{s\leq R}|\bm X_s|^2]<\infty$. With $q(s):=\Ex[\Delta\bm X_s\cdot\Delta\bm Y_s]$, integrating Lemma~\ref{lem:pairing}(i) gives
\begin{equation}\label{eq:erg-unique-pairing}
q(r)-q(s)
\leq-C_L\int_s^r
\Ex\bigl[|\Delta\bm X_u|^2+|\Delta\bm a_u|^2\bigr]du,
\qquad 0\leq s\leq r.
\end{equation}

Since $q(0)=0$ and $q$ is bounded by the uniform second moments,
\eqref{eq:erg-unique-pairing} implies
\(
\int_0^\infty
\Ex\bigl[|\Delta\bm X_s|^2+|\Delta\bm a_s|^2\bigr]ds<\infty.
\)
Choose $s_n\to\infty$ with
$\Ex[|\Delta\bm X_{s_n}|^2]\to0$. The uniform bound for
$\Ex[|\Delta\bm Y_s|^2]$ gives $q(s_n)\to0$. Since $q$ is
nonincreasing and $q(0)=0$, it follows that $q\equiv0$. Hence both
differences in the preceding integral vanish. Continuity identifies
the states, and $Y^i=-D_aL^i(X^i,a^i)$ identifies the adjoints.
Subtracting the backward equations and using It\^o's isometry on each
finite interval identifies $\bm Z$.
\end{proof}

The time-averaged characterization can also be stated directly in
terms of the canonical Pontryagin rule. This formulation keeps the
opponents' controls frozen as open-loop processes in every unilateral
deviation.

\begin{theorem}[Averaged Pontryagin characterization]
\label{thm:erg-Pontryagin-residual}
A profile $\bm\beta\in\mathcal{A}_\infty(\bm x)$ is a strong average-cost
equilibrium if and only if
\begin{equation}\label{eq:erg-Pontryagin-residual}
\lim_{R\to\infty}\frac1R\Ex\left[\int_0^R
\sum_{i=1}^N
\left|
\beta_s^i+D_pH^i\bigl(X_s^{i,\bm\beta},v_\infty^i(\bm X_s^{\bm\beta})\bigr)
\right|^2ds\right]=0.
\end{equation}
In particular, this condition characterizes equilibrium up to
time-averaged equivalence, whereas
\eqref{eq:erg-optimal-control} selects the unique stable Pontryagin
solution pointwise.
\end{theorem}

\begin{proof}
Let $L_b$ be a global Lipschitz constant for $b$, and put
\[
\bm e_s:=\bm\beta_s-b(\bm X_s^{\bm\beta}),\qquad
\Delta\bm X_s:=\bm X_s^{\bm\beta}-\bm X_s^\infty,
\qquad
\Delta\bm a_s:=\bm\beta_s-\bm\alpha_s^\infty.
\]
If $\bm\beta$ is a strong equilibrium, Proposition~\ref{thm:erg-uniqueness}
and the bound
$|\bm e_s|^2\leq2|\Delta\bm a_s|^2+2L_b^2|\Delta\bm X_s|^2$
give \eqref{eq:erg-Pontryagin-residual}.

For the converse, set
\[
f(s):=\Ex[|\Delta\bm X_s|^2]\qquad\text{and}\qquad q(s):=\Ex\Bigl[\Delta\bm X_s\cdot\bigl(v_\infty(\bm X_s^{\bm\beta})-v_\infty(\bm X_s^\infty)\bigr)\Bigr].
\]
By the definition of $\bm e$, the controlled state solves $d\bm X_s^{\bm\beta}=\bigl(b(\bm X_s^{\bm\beta})+\bm e_s\bigr)\,ds+\sqrt2\Sigma\,d\bm W_s$, and $\bm e$ is square integrable on finite time intervals, because $\mathcal E(\bm\beta)<\infty$ and $b$ has at most linear growth.
Lemma~\ref{lem:pairing}(i), applied on every interval $[0,R]$ to the states $\bm X^{\bm\beta}$ and $\bm X^\infty$, with $\bm Y=v_\infty(\bm X)$ and $Z^{i,j}=D_{x^j}v^i_\infty(\bm X)\sigma^j$ for $\bm X\in\{\bm X^{\bm\beta},\bm X^\infty\}$ and with $\bm\phi^1=\bm e$, $\bm\psi^1=D_{\bm x}v_\infty(\bm X^{\bm\beta})\bm e$ and $\bm\phi^2=\bm\psi^2=0$, followed by Young's inequality as in the proof of~\eqref{eq:erg-perturbed-pairing}, gives
\[
q'(s)\leq-\frac{C_L}{2}f(s)+\frac{2K^2}{C_L}\Ex[|\bm e_s|^2],\quad\text{for almost every }s\geq0.
\]
By Lemma~\ref{lem:erg-strong-monotonicity}, $q(s)\geq0$, while $q(0)=0$. Integrating the preceding inequality over $[0,R]$ therefore gives
\[
\frac{C_L}{2}\int_0^Rf(s)\,ds
\leq\frac{2K^2}{C_L}\int_0^R\Ex[|\bm e_s|^2]ds.
\]
Since $\Delta\bm a_s=b(\bm X_s^{\bm\beta})-b(\bm X_s^\infty)+\bm e_s$ and $b$ is $L_b$-Lipschitz, also
\[
\int_0^R\Ex[|\Delta\bm a_s|^2]ds
\leq2L_b^2\int_0^Rf(s)\,ds
+2\int_0^R\Ex[|\bm e_s|^2]ds.
\]
Division by $R$ proves \eqref{erg:cesaro} under
\eqref{eq:erg-Pontryagin-residual}. The converse implication in
Proposition~\ref{thm:erg-uniqueness} completes the proof.
\end{proof}

\begin{corollary}[Uniqueness of the stationary rule]
\label{cor:erg-stationary-rule-unique}
Let $\widetilde b:\R^{Nd}\to\R^{Nd}$ be globally Lipschitz.
Suppose the strong solution
\[
d\widetilde{\bm X}_s
=\widetilde b(\widetilde{\bm X}_s)ds
+\sqrt2\Sigma\,d\bm W_s,
\qquad\widetilde{\bm X}_0=\bm x,
\]
induces a profile
$\bm\beta_s=\widetilde b(\widetilde{\bm X}_s)$ in
$\mathcal{A}_\infty(\bm x)$ which is a strong average-cost open-loop
equilibrium. Then $\widetilde b=b$ on $\R^{Nd}$, and
the corresponding state and control processes are the canonical
ones. Thus the stationary rule generating such an open-loop
equilibrium is unique.
\end{corollary}

\begin{proof}
Set $g(\bm y):=|\widetilde b(\bm y)-b(\bm y)|^2$.
Global Lipschitz continuity of the two rules implies
\[
|g(\bm y)-g(\bm z)|
\leq C(1+|\bm y|+|\bm z|)|\bm y-\bm z|.
\]
Theorem~\ref{thm:erg-Pontryagin-residual} gives
$R^{-1}\Ex[\int_0^Rg(\widetilde{\bm X}_s)ds]\to0$.
Proposition \ref{thm:erg-uniqueness}, the preceding weighted Lipschitz
bound, and Cauchy--Schwarz with the mean-energy bounds give
\[
\frac1R\Ex\left[\int_0^R
\left|g(\widetilde{\bm X}_s)-g(\bm X_s^\infty)\right|ds\right]\xrightarrow{R\to +\infty}0.
\]
The convergence established above for weighted Lipschitz functions now implies
\[
\int_{\R^{Nd}}g\,d\mu_\infty^N
=\lim_{R\to\infty}\frac1R\Ex\left[\int_0^Rg(\bm X_s^\infty)ds\right]=0.
\]
Since $g$ is nonnegative and continuous and $\mu_\infty^N$ has full
support, $g$ vanishes everywhere. Strong uniqueness for the SDE then
identifies the state processes and their controls.
\end{proof}



\bibliographystyle{alpha}
\bibliography{ref}

\end{document}